\documentclass{article}
\usepackage{amsmath}
\usepackage{graphicx}
\usepackage{enumerate}
\usepackage{natbib}
\usepackage{url} 

\newcommand{\blind}{1}
\usepackage{fullpage}

\usepackage{amsthm,amssymb,mathtools}
\usepackage{altvars,mathextras}

\usepackage{enumitem}
\setlist{leftmargin=*}
\setlist[enumerate]{align=left}

\graphicspath{{figs/}}
\usepackage{caption, subcaption}
\usepackage[algo2e,linesnumbered,ruled]{algorithm2e}

\usepackage[hidelinks]{hyperref}
\usepackage[capitalise]{cleveref}
\crefname{equation}{}{}
\Crefname{equation}{Equation}{Equations}
\crefname{figure}{Fig.}{Figs.}
\crefname{condition}{Condition}{Conditions}
\newtheorem{theorem}{Theorem}[section]
\newtheorem{lemma}[theorem]{Lemma}
\newtheorem{condition}[theorem]{Condition}
\newtheorem{corollary}[theorem]{Corollary}
\newtheorem{proposition}[theorem]{Proposition}
\newtheorem{remark}[theorem]{Remark}

\newcommand{\ipto}{\overset{i.p.}{\to}}
\newcommand{\paquant}{q}
\newcommand{\ep}{\varepsilon}

\usepackage{xcolor}
\newcommand{\edit}[1]{{#1}}

\begin{document}

\title{Selecting the number of components in PCA via random signflips}
\author{%
  David Hong%
  \thanks{Department of Electrical and Computer Engineering, University of Delaware. Email: \href{mailto:hong@udel.edu}{hong@udel.edu}}
  \and
  Yue Sheng%
  \thanks{Graduate Group in Applied Mathematics and Computational Science, University of Pennsylvania}
  \and
  Edgar Dobriban%
  \thanks{Department of Statistics and Data Science, University of Pennsylvania}}
\maketitle

\begin{abstract}

Principal component analysis (PCA) is a foundational tool in modern data analysis,
and a crucial step in PCA is selecting the number of components to keep.
However, classical selection methods
(e.g., scree plots, parallel analysis, etc.)
lack statistical guarantees
in the increasingly common setting
of large-dimensional data with heterogeneous noise,
i.e., where each entry may have a different noise variance.
Moreover, it turns out that these methods,
which are highly effective for homogeneous noise,
can fail dramatically for data with heterogeneous noise.
This paper proposes a new method called signflip parallel analysis (FlipPA) for the setting of approximately symmetric noise:
it compares the data singular values to those of ``empirical null'' matrices
generated by flipping the sign of each entry randomly with probability one-half.
We develop a rigorous theory for FlipPA,
showing that it has nonasymptotic type~I error control
and that it consistently selects the correct rank
for signals rising above the noise floor
in the large-dimensional limit
(even when the noise is heterogeneous).
We also rigorously explain why classical permutation-based parallel analysis
degrades under heterogeneous noise.
Finally, we illustrate that FlipPA compares favorably to state-of-the-art methods
via numerical simulations
and an illustration on data coming from astronomy.

\end{abstract}

\noindent%


\section{Introduction}

Discovering latent low-dimensional structure
in large and noisy datasets is one of the central challenges faced in modern data analysis.
Indeed, examples arise across virtually all of science and engineering,
and unsupervised dimensionality reduction
is a standard element of statistical analysis.
In particular,
Factor Analysis (FA) and Principal Component Analysis (PCA)
remain incredibly popular and successful techniques.
They continue to be integral parts of myriad data analysis pipelines,
being performed routinely in thousands of studies every year.
Applications abound in
psychology and education \citep{horn1965, tran2009parallel},
public health \citep{patil2010factor},
management and marketing \citep{stewart1981factor},
economics and finance \citep{bai2002determining, ahn2013eigenvalue},
genomics \citep{lin2016pnas, yano2019pca},
environmental sensing \citep{subbarao1996factor},
and manufacturing \citep{apley2001factor},
to name just a few.
See, e.g., \citet{anderson1958introduction,jolliffe2002principal,yao2015large},
for further references.

Given measurements of $p$ features (covariates)
over a set of $n$ observations (datapoints),
FA and PCA identify components driving variation in the data.
However, these components do not all capture \emph{meaningful} variation, i.e., signal;
many capture variation simply due to noise.
Hence, an important question is:
how many components capture signals rising above the noise?
Informally stated, we have data
that consists of a low-rank signal $\bmS$ plus noise $\bmN$:
$
\bmX = \bmS + \bmN \in \bbR^{n \times p}
$,
$k \coloneqq \rank \bmS \ll n$.
The problem is to estimate the rank $k$ given only $\bmX$.%
\footnote{A key point is that weak signal components can get buried in the noise,
and as a result the corresponding principal components of $\bmX$ completely fail to capture the signal
\citep[see, e.g.,][]{baik2005phase,johnstone2018pca}.
In this case, the goal is to estimate only the number of components
that actually capture the signal and not the noise,
i.e., the number of signal components rising above the noise.}
This paper tackles this problem
in the increasingly common
setting where the noise can be heterogeneous.
As we will show below,
methods that do not appropriately account for heterogeneity
can dramatically degrade,
and theory developed for homogeneous cases does not directly apply.

\subsection{Related work}
\label{sec:related}

Estimating the rank is well known to significantly impact downstream data analyses,
with the standard textbook \citet{brown2014book} calling it
``the most crucial decision'' in exploratory FA.
Choosing too few factors deprives downstream steps of potentially crucial information,
while choosing too many passes on unnecessary noise.
Moreover, data in many important applications have weak ``emergent'' factors
that produce signal singular values of the same order as those coming from the noise,
making them challenging to identify.%
\footnote{Notably, thresholding methods based on bounding the operator norm of noise (e.g., \citet{chatterjee2015meb})
can be overly conservative. See \cref{sec:weak:signals}.}
Such settings are common, e.g., in behavioral and biological sciences.
As a result, a tremendous amount of work has gone into the development of effective methods for this problem;
we give a brief overview of some related literature here.

Classical and standard methods
include Cattell's scree plot \citep{cattell1966scree, cattell1977scree},
sphericity tests based on likelihood ratios \citep{bartlett1954chisquare, lawley1956test},
the Kaiser-Guttman criterion for correlation matrices \citep{guttman1954snc,kaiser1960factor},
the minimum average partial test \citep{velicer1976partial},
and approaches based on minimum description length \citep{wax1985detection, fishler2002detection}.
A popular choice is parallel analysis \citep{horn1965, buja1992rop},
which computes a parallel set of singular values
(typically by computing singular values after
applying random permutations to each column
of the data matrix),
and selects all the components whose singular values
rise above some quantile of their parallel counterparts.
\citet{owen2016bicross} note that
``there is a large amount of evidence that PA
is one of the most accurate [...] classical methods for determining the number of factors''.
Indeed many works find parallel analysis (PA) to be highly effective;
see, e.g., \citet[Section~1.2]{dobriban2017permutation} and references therein.

More recently, tremendous progress has been made
in developing and analyzing methods for large-dimensional data
by using modern insights from high-dimensional probability and random matrix theory.
For example,
several works have developed methods
based on the asymptotic behavior of the eigenvalues of sample covariance matrices
(including the behavior of the differences or the ratios of consecutive eigenvalues)
under various assumptions
\citep{kapetanios2004factor,kapetanios2010factor,kritchman2009non,onatski2010factor,lam2012factor,ahn2013eigenvalue,passemier2014estimation,li2017identifying}.
Other works have studied methods based on various information criteria
\citep{bai2002determining,nadakuditi2008sample,alessi2010factor},
including very recent work by \citet{bai2018consistency} and \citet{hu2020lla:arxiv:v1}.
Other recently proposed approaches include
bi-cross-validation \citep{owen2016bicross},
double cross validation \citep{zeng2019double},
and a random matrix theoretic thresholding of eigenvalues
from the sample correlation matrix \citep{fan2020eno}.
A few works focus specifically on rank estimation for heterogeneous noise
\citep{ke2021eot,landa2021brt:arxiv:v2},
which we will discuss in greater detail below.
Related to parallel analysis,
\citet{dobriban2017permutation}
analyzed permutation-based parallel analysis
for large-dimensional data,
and \citet{dobriban2018dpa} proposed a derandomized variant.
Finally, we note that rank estimation for large-dimensional data
under various assumptions
remains an incredibly active research field;
see, e.g.,
\citet{cai2017spike,xu2022eigenvalue}
for some recent works
that study other settings
(e.g., diverging spikes, correlated data, etc.).

\subsection{Heterogeneous noise and the need for new methods}

\begin{figure} \centering
  \begin{subfigure}[b]{\linewidth} \centering
    \includegraphics[scale=0.22]{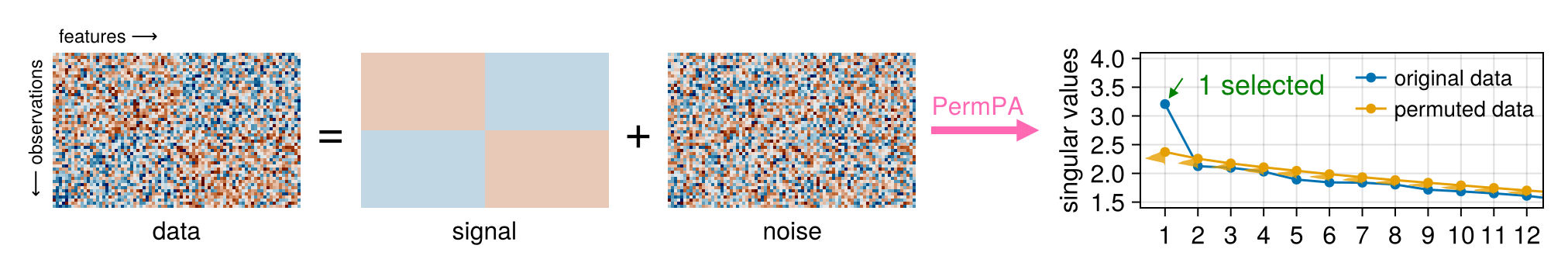}
    \caption{Permutation-based parallel analysis (PermPA) for a rank-one signal in homogeneous noise.}
    \label{fig:illustration:permpa:hom}
  \end{subfigure}

  \begin{subfigure}[b]{\linewidth} \centering
    \includegraphics[scale=0.22]{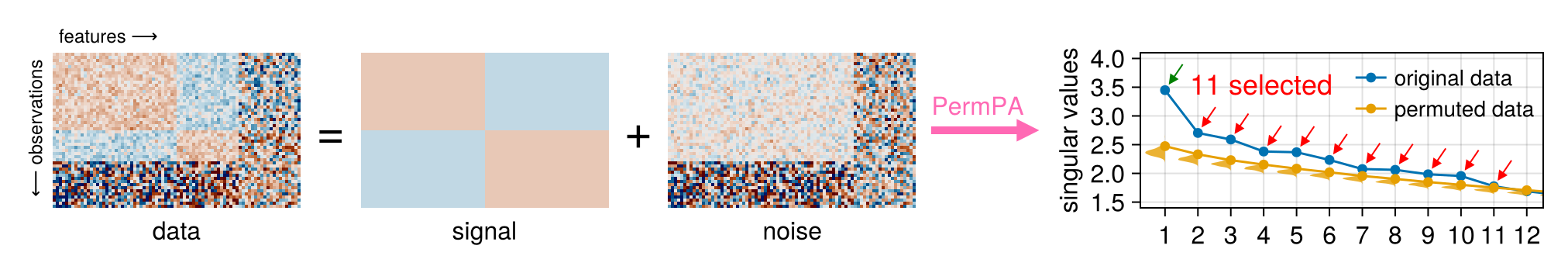}
    \caption{Permutation-based parallel analysis (PermPA) for a rank-one signal in heterogeneous noise.}
    \label{fig:illustration:permpa:het}
  \end{subfigure}

  \caption{Illustrative example of permutation-based parallel analysis (PermPA)
    for a simple rank-one signal in homogeneous and heterogeneous noise.
    PermPA is incredibly effective for homogeneous noise,
    but can substantially over-estimate when the noise is heterogeneous.}
  \label{fig:illustration:permpa}
\end{figure}

In modern applications,
data are now commonly not only large-dimensional
but also heterogeneous.
Heterogeneous noise in particular occurs often in modern settings,
whether due to heteroscedasticity in the features
or due to heterogeneous quality among the observations.
For example,
the noise level in medical imaging varies both within images and from image to image.
Likewise, atmospheric corruptions in astronomical image data
vary both from night to night and from pixel to pixel,
and the quality of environmental sensors can vary from location to location.
These types of
effects all contribute to
noise heterogeneity.

However,
much of the existing work to date has been for data with homogeneous noise,
and methods made for homogeneous noise can dramatically degrade when the noise is heterogeneous.
For example,
as mentioned above,
permutation-based parallel analysis
(which we will refer to as ``PermPA'')
 \citep{buja1992rop}
is highly effective for homogeneous noise
and even remains effective when the data is large-dimensional \citep{dobriban2017permutation}.
It turns out, however, that PermPA can substantially over-estimate the rank
when the noise is heterogeneous,
as shown in the example of \cref{fig:illustration:permpa}.
The entries of the noise in \cref{fig:illustration:permpa:het}
have heterogeneous variances:
the data matrix is moderately noisy on the right side,
least noisy in the upper left,
and most noisy in the lower left.
PermPA correctly selects one component in the homogeneous setting
of \cref{fig:illustration:permpa:hom},
but incorrectly selects eleven components in the heterogeneous setting
of \cref{fig:illustration:permpa:het}.%
\footnote{The over-estimation in \cref{fig:illustration:permpa:het}
is not unique to this particular example
and is a general result
(which we rigorously characterize in \cref{sec:analysis:perm})
that occurs when the noise is heterogeneous.}
Similarly, we will show in experiments that other highly effective methods for homogeneous noise
can substantially degrade for heterogeneous noise.

Consequently, recent works have begun to study
how to properly account for
noise heterogeneity
when carrying out PCA for large data.
Great progress has been made on improving the quality of the estimated components
\citep[see, e.g.,][]{zhang2018hpa:arxiv:v2,leeb2018oss:arxiv:v4,leeb2019mdf:arxiv:v3,hong2016tat,hong2018asymptotic,hong2023owp},
but relatively
fewer works have addressed
how to estimate the \emph{number} of components
in these heterogeneous settings.
For heterogeneity across features,
\citet{leeb2018oss:arxiv:v4}
consider selecting the singular values rising
above the asymptotic operator norm of the noise matrix,
which can be predicted when the noise is whitened
or when the noise variances are well-estimated.
\citet{ke2021eot}
consider noise variances drawn from a Gamma distribution,
and propose setting a threshold based on fitting
the associated Mar{\v{c}}enko-Pastur distribution to the bulk singular values.
\citet{landa2021brt:arxiv:v2} propose rescaling the rows and columns
of the data matrix to biwhiten the data,
then selecting the singular values
rising above the same cut-off as applies
for large-dimensional homogeneous noise.
They show that the biwhitening procedure makes
the formerly heterogeneous noise
behave like homogeneous noise, making this a good choice of cut-off.

Overall, the analysis of large-dimensional data with heterogeneous noise
is an actively developing area.
In this paper,
we develop an elegant and flexible method in the style of
parallel analysis
that is not only effective for homogeneous noise
but also remains effective for heterogeneous noise.
The method is intuitive and simple to implement,
and we provide theoretical guarantees
(including not only asymptotic consistency, but also nonasymptotic type~I error control in selecting incorrect components)
leveraging modern insights
from random matrix theory and randomization-based testing.
In contrast to other approaches,
the proposed signflip parallel analysis (FlipPA) method
allows for entry-wise noise heterogeneity,
does not need exact or distributional knowledge of the variances,
does not modify the data,
and provides nonasymptotic type~I error control.
We show in a broad range of experiments that our method compares favorably to prior approaches
\citep{buja1992rop,bai2018consistency,hu2020lla:arxiv:v1,fan2020eno,ke2021eot,landa2021brt:arxiv:v2}.

\subsection{Organization}

The remainder of the paper is organized as follows.
\Cref{sec:method} describes the proposed signflip parallel analysis (FlipPA) method.
\Cref{sec:theory} rigorously analyzes FlipPA, 
discussing nonasymptotic type~I error control and
establishing asymptotic consistency (\cref{sec:consistency});
\cref{sec:analysis:perm} rigorously explains
why PermPA degrades under heterogeneous noise.
\Cref{sec:simulations} demonstrates the performance of FlipPA
through numerical simulations
and compares with several state-of-the-art methods.
\Cref{sec:exp} demonstrates FlipPA on an empirical dataset
coming from astronomy.
Finally, \cref{sec:discuss} concludes with a discussion
of future directions.


\section{Proposed method: Signflip parallel analysis (FlipPA)}
\label{sec:method}

We propose a new method for rank estimation called \emph{signflip parallel analysis} (FlipPA).
FlipPA is a method in the style of parallel analysis \citep{horn1965, buja1992rop}
that largely retains the excellent battle-tested performance of parallel analysis
under homogeneous noise,
while expanding these benefits to data with heterogeneous noise.
The idea is to generate \emph{parallel} matrices
that look like pure noise,
then select the components whose singular values rise above
data-driven cut-offs determined by their parallel counterparts.
FlipPA generates parallel matrices
from the data matrix
by randomly flipping the signs
of its entries,
i.e., each entry is signflipped
independently with probability $1/2$.
Each parallel matrix is formed using a new independent set of random signflips.
A feature of this approach is that it provides non-asymptotic type~I error control
at a user-specified level (see \cref{sec:theory}).

\begin{algorithm2e}[t]
  \caption{Signflip parallel analysis (FlipPA)}
  \label{alg:pa:signflip}
  \SetKwInOut{Input}{Input}
  \SetKwInOut{Output}{Output}

  \Input{Data $\bmX \in \bbR^{n \times p}$,
    quantile $\paquant \in [0,1]$,
    number of trials $T$,
    threshold $\tau \geq 0$.}
  \Output{Selected number of factors $\htk$.}
  \For{$t = 1$ \KwTo $T$}{
    \label{alg:pa:signflip:line:for}
    Randomly signflip the entries of $\bmX$:
    \begin{equation}
      \label{alg:pa:signflip:flipping}
      \tlX^{(t)}_{ij}
      \overset{\text{ind}}{\sim}
      \begin{cases}
        +X_{ij} ,& \text{with probability } 1/2, \\
        -X_{ij} ,& \text{with probability } 1/2,
      \end{cases}
    \end{equation}
    i.e., $\btlX^{(t)} \gets \bmR^{(t)} \circ \bmX$
    where $\bmR^{(t)} \in \{-1,+1\}^{n \times p}$ has i.i.d. Rademacher entries\;
    \label{alg:pa:signflip:line}

    $\btlsigma^{(t)} \gets \text{singular values of $\btlX^{(t)}$}$\;
    \label{alg:pa:signflip:line:paval}
  }
  \label{alg:pa:signflip:line:endfor}
  $\bmsigma \gets \text{singular values of $\bmX$}$\;
  \label{alg:pa:signflip:line:dataval}
  $\htk \gets$ first $k$ for which
  \begin{subequations}
  \begin{align}
    \sigma_{k+1} &\leq \text{$\paquant$-quantile of} \left(\tlsigma_{k+1}^{(1)},\dots,\tlsigma_{k+1}^{(T)}\right) + \tau
    , &&\text{if using pairwise comparison}
    \label{alg:pa:signflip:pairwise}
    , \\
    \sigma_{k+1} &\leq \text{$\paquant$-quantile of} \left(\tlsigma_1^{(1)},\dots,\tlsigma_1^{(T)}\right) + \tau
    , &&\text{if using upper-edge comparison}
    \label{alg:pa:signflip:upperedge}
    .
  \end{align}
  \end{subequations}
  \label{alg:pa:signflip:line:cmp}
\end{algorithm2e}

\begin{figure} \centering
  \includegraphics[scale=0.22]{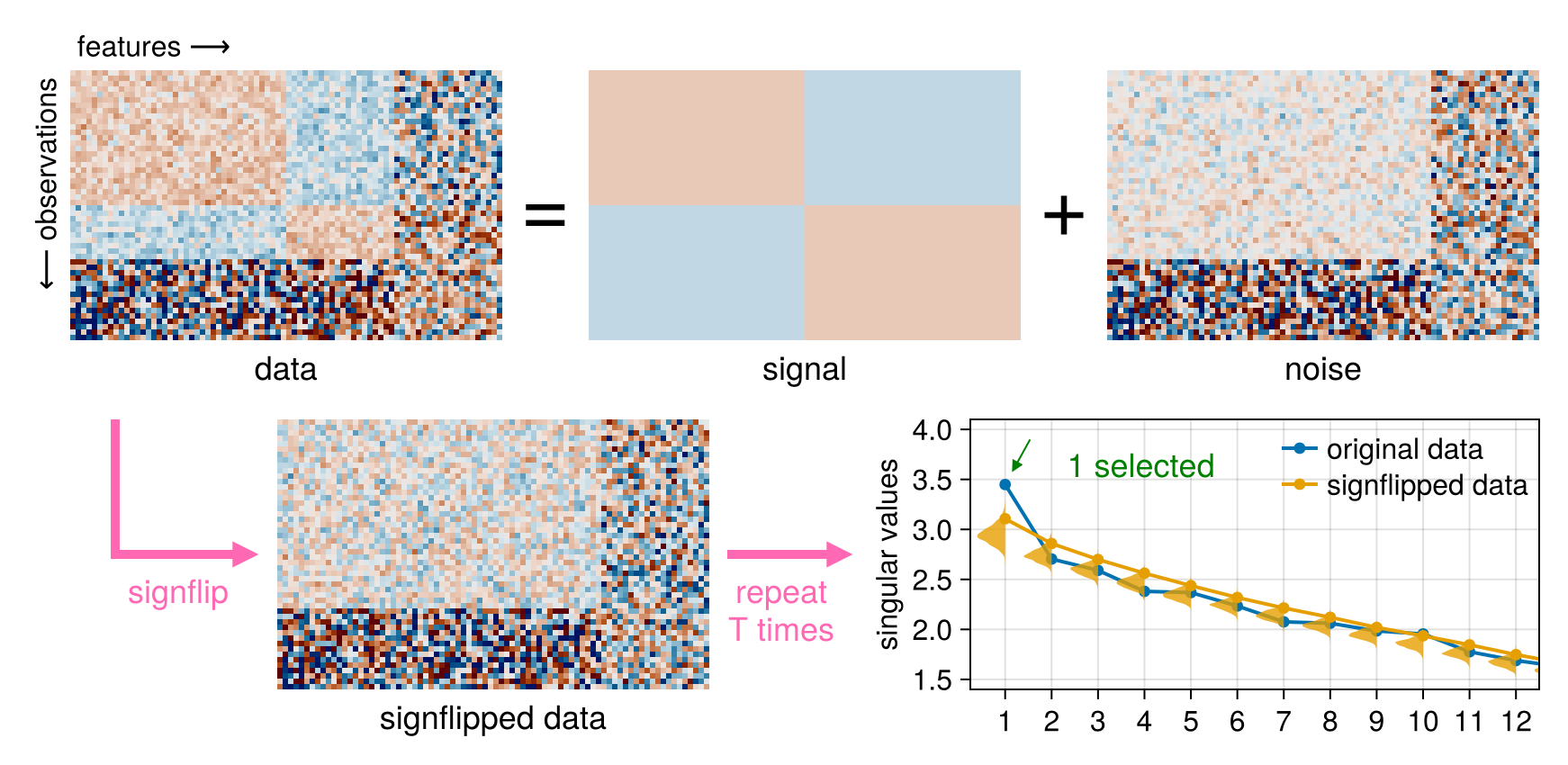}
  \caption{Illustration of signflip parallel analysis (FlipPA) for heterogeneous noise.
    FlipPA selects the data singular values
    that rise above the $\paquant$-quantiles of
    those from the signflipped data,
    using the empirical distributions
    from $T$ signflips
    (shown here as vertical densities).}
  \label{fig:illustration:flippa:het}
\end{figure}

\Cref{alg:pa:signflip} provides a detailed description of FlipPA,
and \cref{fig:illustration:flippa:het} illustrates
the process using the heterogeneous data
from \cref{fig:illustration:permpa:het}.
The first step
(\cref{alg:pa:signflip:line:for,alg:pa:signflip:line,alg:pa:signflip:line:paval,alg:pa:signflip:line:endfor})
is to compute $T \ge 1$ sets of parallel singular values
$\btlsigma^{(1)},\dots,\btlsigma^{(T)}$
from the signflipped data matrices $\btlX^{(1)},\dots,\btlX^{(T)}$ in \cref{alg:pa:signflip:flipping};%
\footnote{Note that
only the singular values $\btlsigma^{(1)},\dots,\btlsigma^{(T)}$ need to be saved;
the memory for storing $\btlX^{(t)}$ can be reused across trials.}
\cref{fig:illustration:flippa:het}
illustrates the signflipped data matrix from one of the trials.
The next and final step
(\cref{alg:pa:signflip:line:dataval,alg:pa:signflip:line:cmp})
is to compute the singular values of the data $\bmX$
and compare them against
a chosen quantile $q$ (e.g., $50\%$, $95\%$, or $99\%$)
computed with their parallel counterparts.
There are two choices for the comparison method:%
\footnote{Note that both comparison methods are sequential:
the selection ends once a data singular value $\sigma_{k}$
falls below its appropriate parallel cut-off,
even if singular values later on rise above again.}
\begin{itemize}
  \item pairwise comparison \cref{alg:pa:signflip:pairwise}
  selects all the leading components
  for which the data singular value $\sigma_{k}$
  rises $\tau$-above the $q$-quantile
  of its \emph{parallel counterparts} $\tlsigma_{k}^{(1)},\dots,\tlsigma_{k}^{(T)}$,

  \item upper-edge comparison \cref{alg:pa:signflip:upperedge}
  selects all the leading components
  for which the data singular value $\sigma_{k}$
  rises $\tau$-above the $q$-quantile
  of the \emph{parallel upper-edges} $\tlsigma_{1}^{(1)},\dots,\tlsigma_{1}^{(T)}$.
\end{itemize}
The scatter plot in \cref{fig:illustration:flippa:het}
shows the singular values $\sigma_{k}$ of the original data $\bmX$,
overlaid with each empirical distribution
for the parallel singular values $\tlsigma_{k}^{(1)},\dots,\tlsigma_{k}^{(T)}$
and the resulting $q$-quantiles.
Both comparison methods correctly select
one component
with $\tau = 0$ here.

While pairwise comparison is the classical selection method \citep{horn1965, buja1992rop},
upper-edge comparison can have some benefits.
For example,
upper-edge comparison never selects more factors than pairwise comparison,
making it more conservative
and less likely to over-estimate the rank.
Moreover, upper-edge comparison
can be computationally cheaper,
since it only requires computing and storing
the first singular value
$\tlsigma_{1}^{(t)}$
of the signflipped data
(i.e., the operator norm $\|\btlX^{(t)}\|$)
in each trial.
At the same time,
upper-edge comparison is also more likely
to under-estimate the rank,
which can be an issue if strong signals
are shadowing weak signals in the data.
The two comparison methods also agree in many cases,
especially when the data dimensions are large.
When they disagree,
which comparison method to choose will depend on the application and the salient priorities.
We suggest using upper-edge comparison as a default if over-selection of the number of components or a large computational cost are important issues.

As a parallel analysis method,
FlipPA has the same structure
as the permutation-based parallel analysis (PermPA) of \citet{buja1992rop}.
The crucial difference is that
PermPA creates parallel matrices
by randomly permuting the entries of each column of the data matrix
while FlipPA uses random entrywise signflips.
Thus, FlipPA immediately shares many of the practical benefits of PermPA,
being intuitive and simple to implement.%
\if1\blind%
\footnote{A Julia package implementing the method is available at:
\url{https://github.com/dahong67/FlipPA.jl}}
\else{ }\fi
Moreover, it is a flexible method.
As we will demonstrate in \cref{sec:sim:block:dep},
it can be straightforwardly modified,
e.g., to handle noise whose entries have blockwise dependence.

The remaining sections
investigate
the performance of FlipPA
theoretically (\cref{sec:theory}),
via simulations (\cref{sec:simulations}),
and via an illustration on data from astronomy (\cref{sec:exp}).


\section{Theoretical Properties}
\label{sec:theory}

This section studies the theoretical properties 
of the proposed FlipPA method.
Throughout this section,
we consider data of the following signal-plus-noise form:
\begin{equation}
  \label{eq:signal:plus:noise}
  \bmX = \underbracket{\sum_{i = 1}^{k} \theta_i \bmu_i \bmz_i^\top}_{\eqqcolon \bmS}
  +
  \bmN
  \in \bbR^{n \times p}
  ,
\end{equation}
where the signal $\bmS$ is of rank at most $k$,
and
$\theta_i\ge0$, $\bmu_i \in \bbR^n$, and $\bmz_i \in \bbR^p$, for $i\in\{1,\ldots, k\}$ may be deterministic or random and are not necessarily orthogonal, as discussed below.
Furthermore, the signal component coefficients $\theta_i$ and the number of terms $k$ may both change with $n$ and $p$.
The noise matrix $\bmN$ satisfies the following condition.

\begin{condition}[Noise with independent symmetric entries]
  \label{assump:noise:sym:ind}
  The noise matrix $\bmN$
  has independent entries
  with symmetric distributions,
  i.e., $ N_{ij} =_d -N_{ij}$ for all $i, j$.
\end{condition}

Importantly, this does not require the noise entries to be identically distributed.
Moreover, as discussed below in \cref{rem:noise:sym},
the noise entries need not have symmetric distributions
for FlipPA to be consistent.
Another slightly subtle point is that
the precise output of FlipPA
can depend on how one computes sample quantiles;
all the results in this section hold
for sample quantiles
computed using any of the nine commonly used definitions
enumerated in \cite{hyndman1996sqi}.

Consider first the nonasymptotic type~I error rate of FlipPA
under the null hypothesis of $H_0 : k = 0$.
There is no signal here,
and the type~I error rate of FlipPA
is $\Pr_{H_0}\bigl[\htk > 0 \bigr]$,
where $\htk$ is the output of FlipPA.
Combining some general arguments for randomization tests \citep{hemerik2017etw}
with some careful but straightforward additional analysis
(to connect sample quantiles to order statistics)
yields the following \lcnamecref{prop:type-I}.

\begin{proposition}[Type~I error control]
  \label{prop:type-I}
  Suppose the signal-plus-noise model~\cref{eq:signal:plus:noise}
  satisfies \cref{assump:noise:sym:ind}.
  Then the type~I error rate of FlipPA, for both upper-edge and pairwise comparison methods,
  is bounded above as $
  \Pr_{H_0}\bigl[\htk > 0\bigr]
  \leq 1 - \lfloor \paquant T \rfloor/(T+1)
  $.
\end{proposition}

For completeness, we provide a proof in \cref{thm:false:positive:proof}.
Note that the type~I error bound
holds for any chosen threshold $\tau \geq 0$ and
takes values in increments of $1/(T+1)$,
starting at $1/(T+1)$ and going up.
These bounds provide a simple way to control the type~I error of FlipPA.
Namely, to obtain a type~I error rate below level $\alpha$,
simply choose the number of trials $T$
and the quantile $\paquant$ so that
$1 - \lfloor \paquant T \rfloor/(T+1) \leq \alpha$,
e.g., one can choose
$\paquant = 1$
and
$T = \lceil 1/\alpha \rceil - 1$,
yielding a type~I error rate bounded as
$ \Pr_{H_0}\bigl[\htk > 0 \bigr]
  \leq 1/\lceil 1/\alpha \rceil\leq \alpha$.
One may choose a larger $T$ to reduce the variability in FlipPA due to the random signflips,
and a correspondingly smaller $\paquant$
to control type~I error at level $\alpha$.
As usual,
the choice of $q$ tunes the trade-off
between the type~I error of FlipPA and its power,
though interestingly, this trade-off can essentially vanish
when $n$ and $p$ are large
(see \cref{sec:sim:roc}).

The following subsections present our main theoretical results.
\Cref{sec:consistency} considers the alternative hypothesis
(i.e., nonzero components)
and proves asymptotic consistency
under suitable conditions.
\Cref{sec:analysis:perm} provides a rigorous explanation
of why permutation-based parallel analysis degrades
under heterogeneous noise.

\subsection{Consistency of FlipPA}
\label{sec:consistency}

Here,
we consider the asymptotic consistency of FlipPA
in the large-dimensional limit
where $n,p \to \infty$.
This asymptotic regime
captures many modern datasets,
for which the number of features is comparable to the number of observations.
For practical purposes,
this asymptotic consistency means that
FlipPA provides accurate rank estimates
with high probability when there are sufficiently many features and many observations.

\begin{figure} \centering
  \includegraphics[scale=0.21]{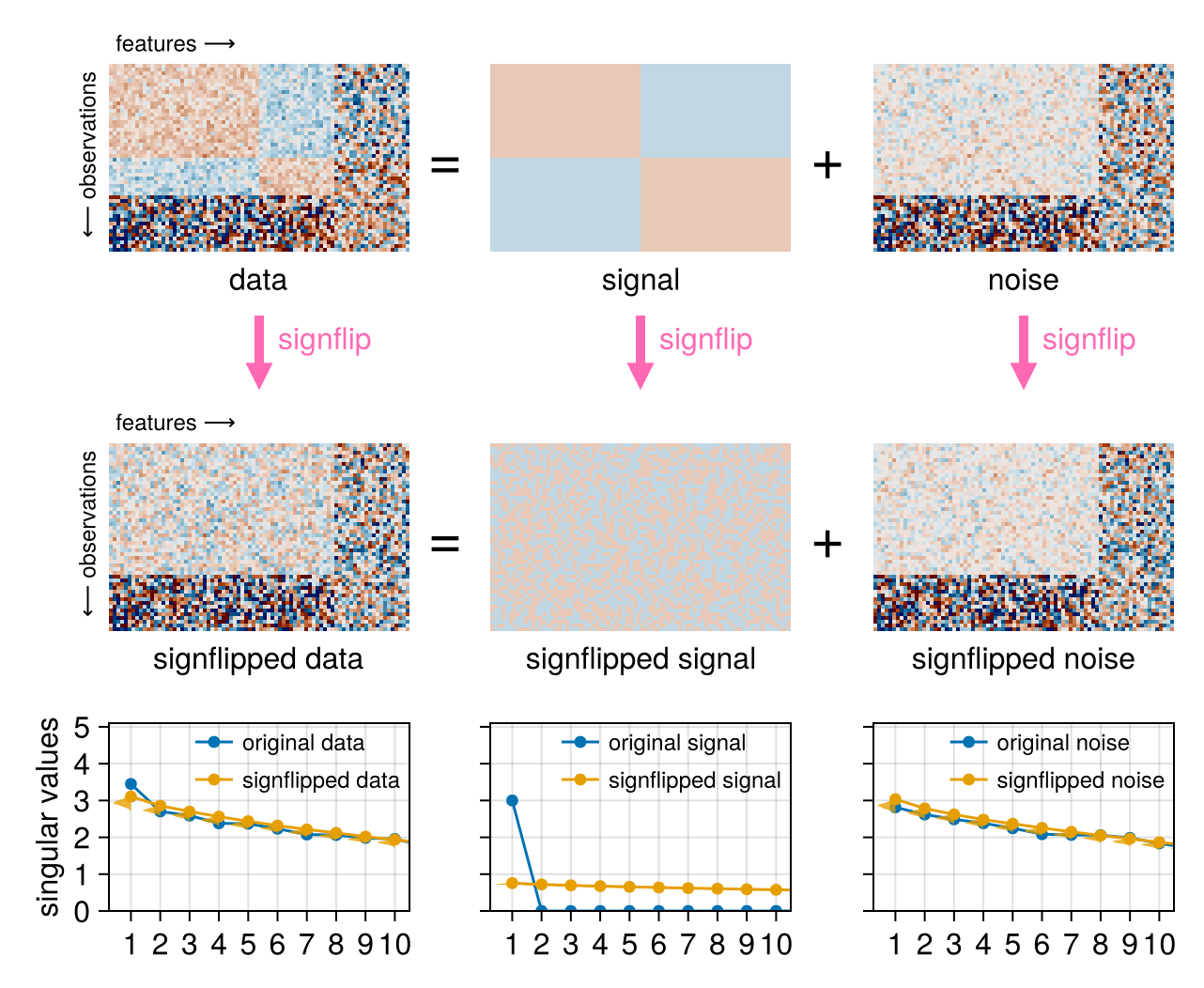}
  \caption{Example illustrating how
    FlipPA scrambles the signal
    while preserving heterogeneous noise,
    resulting in singular values similar to those of the noise.}
  \label{fig:consistency}
\end{figure}

The key insight is
the observation that the signflipping operation in FlipPA
scrambles the signal in the data while preserving heterogeneous noise,
producing matrices that look like the noise,
as illustrated in \cref{fig:consistency}.
Specifically, the signflipped signal has a much reduced operator norm,
while the signflipped noise is similar to the original noise.
As a result,
the singular values of the signflipped data
are very similar to those of the original noise (which we do not observe directly)
and hence lead to a good choice for cut-off in estimating the rank.%
\footnote{Similar observations were also
used to derive properties of PermPA
in \citet{dobriban2017permutation}
and to derive general properties of randomization-based tests
in \citet{dobriban2021consistency}.
However, the results from those works cannot be used to derive the results presented here.}
The remainder of this section develops this intuition
into a rigorous theory.
To simplify the presentation,
here we will focus on FlipPA
with the upper-edge comparison method.
Similar statements can be shown
for the pairwise comparison method
under similar conditions.
\edit{Moreover, we focus here on settings
where signflipping causes the signal to vanish in operator norm,
leading to an accurate estimate of the noise level;
see \cref{sec:relative:strength} for an analysis of settings with strong components that do not vanish.}

To carry out the rigorous analysis,
we consider the following two conditions.

\begin{condition}[Well-defined asymptotic noise upper-edge]
  \label{assump:noise:upper:edge}
  The noise has a well-defined asymptotic upper-edge $\brsigma$,
  i.e.,
  $\| \bmN \| \ipto \brsigma$
  where $\ipto$ denotes convergence in probability
  and $\brsigma$ is a constant (independent of $n$ and $p$).
\end{condition}

This first condition
simply ensures a well-defined asymptotic noise limit.
Essentially it means that the variability of the noise operator norm is small
\emph{relative to its value},
since one can always re-scale the data by the order of the noise operator norm
without loss of generality.
This condition is satisfied in many natural models
studied in random matrix theory \cite[see, e.g.,][]{benaych2012singular},
including ones with heterogeneous noise
(see, e.g., \cite{husson2020ld,ajanki2016ufg,ajanki2019qve}).

\begin{condition}[Asymptotically perceptible signal]
  \label{assump:signal:perceptible}
  The signal is asymptotically perceptible,
  i.e., there exists some $\ep>0$ such that
  $\Pr[ \sigma_k(\bmX) > \| \bmN \| + \ep ] \to 1$.
\end{condition}

This second condition
ensures that the signal rises above the noise
and is perceptible as a result.
As is well known \citep[see, e.g.,][]{baik2005phase, johnstone2018pca},
signals that are too weak
can get buried in the noise
for large-dimensional data,
and the corresponding components
can consequently completely fail
to recover
them.
Our goal is to identify only the components above the noise level.

We are now ready to state the main consistency results;
see \cref{sec:consistency:proofs} for their proofs.

\begin{theorem}[Asymptotic consistency of FlipPA]
  \label{thm:consistency:component}
  Suppose the signal-plus-noise model~\cref{eq:signal:plus:noise}
  satisfies \cref{assump:noise:sym:ind,assump:noise:upper:edge,assump:signal:perceptible}.
  Then FlipPA using the upper-edge comparison method
  with threshold $\tau \in (0,\ep)$
  is consistent,
  i.e.,
  $\Pr\big[ \htk = k \big] \to 1$,
  as long as
  the signal components are delocalized as follows:
  \begin{equation}
    \label{eq:consistency:component}
    \bbE \left\{
      \sum_{i = 1}^{k}
      \theta_i \|\bmu_i\|_2 \|\bmz_i\|_2
      \cdot
      \left[
        \frac{\|\bmu_i\|_\infty/\|\bmu_i\|_2 + \|\bmz_i\|_\infty/\|\bmz_i\|_2}{2}
      \right]
    \right\}
    \to 0
    .
  \end{equation}
\end{theorem}

\begin{remark}[Asymptotic consistency without symmetry]
  \label{rem:noise:sym}
  For simplicity,
  we have stated \cref{thm:consistency:component}
  under the condition that the noise entries have symmetric distributions
  (\cref{assump:noise:sym:ind}).
  This condition can be weakened to drop symmetry
  as long as the signflipped noise $\bmR \circ \bmN$
  shares the same
  asymptotic upper-edge $\brsigma$ as $\bmN$ (from \cref{assump:noise:upper:edge}),
  i.e.,
  as long as $\| \bmR \circ \bmN \| \ipto \brsigma$.
  This weaker condition is satisfied in many natural models
  studied in random matrix theory
  due to universality,
  including ones with heterogeneous noise
  (see, e.g., \cite{ajanki2016ufg}).
  The same applies to \cref{thm:consistency:component:rates,thm:consistency} below.
\end{remark}

\cref{thm:consistency:component}
follows from the more general \cref{thm:consistency} below,
and its proof is provided in \cref{thm:consistency:component:proof}.
The terms
$\|\bmu_i\|_\infty/\|\bmu_i\|_2$ and $\|\bmz_i\|_\infty/\|\bmz_i\|_2$
in \cref{eq:consistency:component}
capture how localized the $i$-th signal component
is in $\bmu_i$ and $\bmz_i$, respectively.
If the signal is completely concentrated on a single entry,
e.g., if $\bmu_i = [1,0,\dots,0]^\top$,
then $\|\bmu_i\|_\infty/\|\bmu_i\|_2$ achieves its maximum possible value of~$1$,
and likewise for $\|\bmz_i\|_\infty/\|\bmz_i\|_2$.
On the other hand,
if the signal is completely spread out,
e.g., if $\bmu_i = [1,\dots,1]^\top$,
then $\|\bmu_i\|_\infty/\|\bmu_i\|_2$ achieves its minimum possible value of $1/\sqrt{n}$,
and again likewise for $\|\bmz_i\|_\infty/\|\bmz_i\|_2$.
Thus, roughly speaking, \cref{thm:consistency:component}
states that FlipPA is consistent as long as
each component delocalizes fast enough
to overcome the overall signal strength of $\theta_i \|\bmu_i\|_2 \|\bmz_i\|_2$.
This delocalization makes the signal incoherent,
which FlipPA is able to exploit to obtain an accurate estimate of the noise operator norm,
and consequently an accurate estimate of the rank.
The following \lcnamecref{thm:consistency:component:rates}
shows how this condition translates
to growth and delocalization rates for the signal components.%
\footnote{To simplify the presentation,
here we suppose the rank $k$
and the signal component coefficients $\theta_1,\dots,\theta_k$
are deterministic,
and that the component vectors $\bmu_1,\dots,\bmu_k$ and $\bmz_1,\dots,\bmz_k$
are all jointly independent.}

\begin{corollary}[Asymptotic consistency of FlipPA in terms of growth rates]
  \label{thm:consistency:component:rates}
  Suppose the signal-plus-noise model \cref{eq:signal:plus:noise}
  satisfies \cref{assump:noise:sym:ind,assump:noise:upper:edge,assump:signal:perceptible},
  where
  \begin{itemize}
    \item the signal rank $k$ and component coefficients $\theta_1,\dots,\theta_k$ are deterministic
    and grow at the following rates:
      $k = O(m^{\nu_1}\log^{\nu_2}m)$
      and
      $\displaystyle \max_{i=1,\dots,k} \! \theta_i = O(m^{\beta_1} \log^{\beta_2} m)$
      where $m = \min(n,p)$
      for some $\nu_1,\nu_2,\beta_1,\beta_2 \!\ge\! 0$,

    \item the component vectors $\bmu_1,\dots,\bmu_k$ and $\bmz_1,\dots,\bmz_k$
    are all jointly independent,
    are normalized such that
    $\bbE \|\bmu_i\|_2 = \bbE \|\bmz_i\|_2 = 1$ for $i \in \{1,\dots,k\}$,
    and delocalize at the following rates:
    $\displaystyle \max_{i=1,\dots,k} \, \bbE \|\bmu_i\|_\infty = O(n^{-\alpha_1} \log^{-\alpha_2} n)$
    and
    $\displaystyle \max_{i=1,\dots,k} \, \bbE \|\bmz_i\|_\infty = O(p^{-\alpha_1} \log^{-\alpha_2} p)$
    for some $\alpha_1,\alpha_2 \in \bbR$.
    \end{itemize}
  Then FlipPA using the upper-edge comparison method
  with threshold $\tau \in (0,\ep)$
  is consistent,
  i.e.,
  $\Pr\big[ \htk = k \big] \to 1$,
  as long as:
  a) $\alpha_1 > \nu_1 + \beta_1$,
  or
  b) $\alpha_1 = \nu_1 + \beta_1$ and $\alpha_2 > \nu_2+\beta_2$.
\end{corollary}

\Cref{thm:consistency:component:rates}
is proven in \cref{thm:consistency:component:rates:proof}, and it
covers many important settings.
For example,
in the commonly studied settings where the component vectors $\bmu_i$ and $\bmz_i$
either:
a) have i.i.d.~sub-Gaussian entries,
or
b) are uniformly distributed on the unit spheres in $\bbR^n$ and $\bbR^p$, respectively,
it follows that%
\footnote{This can be verified using, e.g.,
\citet[Exercise 2.5.10 and Theorem 3.4.6]{vershynin2018hdp}.}
$\max_{i=1,\dots,k} \, \bbE \|\bmu_i\|_\infty = O(n^{-1/2} \log^{1/2} n)$
and
$\max_{i=1,\dots,k} \, \bbE \|\bmz_i\|_\infty = O(p^{-1/2} \log^{1/2} p)$.
Hence, in these settings, $\alpha_1 = 1/2$ and $\alpha_2 = -1/2$,
and FlipPA remains consistent even
if the signal rank and component coefficients both grow,
as long as their combined rate is small enough, i.e.,
 $k \cdot \max_{i=1,\dots,k}  \theta_i= O(m^{\zeta_1}\log^{\zeta_2}m)$
      with
$\zeta_1 < 1/2$
or $\zeta_1=1/2$ with $\zeta_2 < -1/2$.
For instance, this allows $k = m^{1/4}$ and $\max_{i=1,\dots,k}  \theta_i = m^{1/4-\delta}$ for any fixed $\delta>0$.
Moreover, it is common in many applications for the rank to be fixed or bounded,
i.e., $\nu_1 = \nu_2 = 0$,
in which case the condition simply reduces
to
      $\max_{i=1,\dots,k} \, \theta_i = O(m^{\beta_1} \log^{\beta_2} m)$
      with
$\beta_1 < 1/2$
or $\beta_1 = 1/2$ with $\beta_2 < -1/2$.
Overall,
FlipPA remains consistent even if the signal rank and component coefficients grow,
as long as the component vectors delocalize sufficiently quickly
to compensate.

We conclude this section with a discussion
of the more general sufficient conditions for FlipPA consistency,
of which
\cref{thm:consistency:component} (and consequently \cref{thm:consistency:component:rates})
are actually special cases.
The following \lcnamecref{thm:consistency}
shows that FlipPA remains consistent
under an even weaker form of signal delocalization
than that of \cref{eq:consistency:component}.
This weaker form of delocalization is stated
in terms of properties of the signal matrix $\bmS$ directly,
rather than the component coefficients $\theta_i$ and vectors $\bmu_i$ and $\bmz_i$.
In particular, we define the delocalization coefficients
\begin{align}
  \label{eq:delocalization:coeffs}
  \rho_2(\bmS)
  &
  \coloneqq
  \begin{Vmatrix} & \bmS \\ \bmS^\top & \end{Vmatrix}_{2,\infty}
  \sqrt[4]{\log (n+p)}
  ,
  &
  \rho_\infty(\bmS)
  &
  \coloneqq
  \max_{i=1,\dots,n+p} \begin{Vmatrix} & \bmS \\ \bmS^\top & \end{Vmatrix}_{\infty,(i)} \sqrt{\log i}
  ,
\end{align}
where $\|\cdot\|_{2,\infty}$ denotes the maximum row $\ell_2$ norm,
and
$\|\cdot\|_{\infty,(i)}$ denotes the $i$-th largest row $\ell_\infty$ norm,
i.e.,
for any matrix $\bmA \in \bbR^{m \times \ell}$,
$\|\bmA\|_{\infty,(1)} \geq \cdots \geq \|\bmA\|_{\infty,(m)}$
sorts the row $\ell_\infty$ norms $\|\bmA_{1:}\|_\infty, \dots, \|\bmA_{m:}\|_\infty$
in descending order.
We are now ready to state the theorem.

\begin{theorem}[Asymptotic consistency of FlipPA: general condition]
  \label{thm:consistency}
  Suppose the signal-plus-noise model \cref{eq:signal:plus:noise}
  satisfies \cref{assump:noise:sym:ind,assump:noise:upper:edge,assump:signal:perceptible}.
  Then FlipPA using the upper-edge comparison method
  with threshold $\tau \in (0,\ep)$
  is consistent,
  i.e.,
  $\Pr\bigl[ \htk = k \bigr] \to 1$,
  as long as
  $\bbE \| \bmS \|_{2,\infty} \to 0$,
  $\bbE \| \bmS^{\top} \|_{2,\infty} \to 0$, and
  $\min\big[ \bbE \, \rho_2(\bmS), \bbE \, \rho_\infty(\bmS) \big] \to 0$.
\end{theorem}

The proof involves showing  
that the noise is preserved
(so that the signflipped noise looks like the true noise),
while the signal is destroyed by signflipping
(so that the signflipped data looks like the signflipped noise).
See \cref{pf:thm:consistency}.
We leverage recent advances in a line of work in random matrix theory
on nonhomogeneous random matrices \citep{seginer2000ten, latala2005seo, schuett2013ote, bandeira2016sharp, van2017structured, latala2018tdf}.
The delocalization conditions
$\bbE \| \bmS \|_{2,\infty} \to 0$
and
$\bbE \| \bmS^{\top} \|_{2,\infty} \to 0$
ensure
that the signal is not
localized on any particular row or column.
The condition
$\min\big[ \bbE \, \rho_2(\bmS), \bbE \, \rho_\infty(\bmS) \big] \to 0$
ensures
that the decay of row and column norms is sufficiently fast
with respect to either
of the delocalization coefficients
$\rho_2$ or $\rho_\infty$
in \cref{eq:delocalization:coeffs}.
The coefficient $\rho_2$
captures a dimension-dependent form of decay;
$\bbE \, \rho_2(\bmS) \to 0$
is equivalent to
saying that both
$\bbE\|\bmS\|_{2,\infty}$ and $\E\|\bmS^\top\|_{2,\infty}$
decay at a rate of $o(\log^{-1/4}(n+p))$.
The coefficient $\rho_\infty$
captures a dimension-free form of decay,
and is satisfied
in many natural circumstances,
as described
by the following \lcnamecref{thm:decay:suff}
(see \cref{thm:decay:suff:proof} for its proof).

\begin{theorem}[Sufficient conditions for $\rho_\infty$ decay]
  \label{thm:decay:suff}
  We have $\bbE \, \rho_\infty(\bmS) \to 0$,
  under either of the following conditions:
  (1) there is $t \geq 2$ such that $\bbE \|\bmS\|_{t,t} \to 0$;
(2) the rank grows slowly enough that $\bbE \left\{\rank^{1/2}(\bmS) \sqrt{\|\bmS\|_{2,\infty} \cdot \|\bmS^\top\|_{2,\infty}} \right\} \to 0$.
\end{theorem}

The mixed norm $\|\cdot\|_{t,t}$ denotes
the $\ell_t$ norm of the vector made of the row $\ell_t$ norms,
i.e., $\|\bmS\|_{t,t} = \| ( \|\bmS_{1:}\|_t, \dots, \|\bmS_{n:}\|_t ) \|_t$,
which is equivalent to the $\ell_t$ norm of the vectorized matrix.
The row-wise and column-wise delocalization conditions
$\bbE \| \bmS \|_{2,\infty} \to 0$
and
$\bbE \| \bmS^{\top} \|_{2,\infty} \to 0$
imply the second sufficient condition in \cref{thm:decay:suff}
when the rank of $\bmS$ is bounded,
and hence are already sufficient in that case.
Overall,
\cref{thm:consistency,thm:decay:suff} show that
FlipPA is consistent as long as the signal is sufficiently delocalized
in an even weaker sense than in \cref{thm:consistency:component,thm:consistency:component:rates}.


\subsection{Degradation of PermPA under heterogeneous noise}
\label{sec:analysis:perm}

This section gives a rigorous explanation
of why PermPA degrades
under heterogeneous noise.
The key insight is that
randomly permuting the entries of each column,
as is done by PermPA,
has the effect of homogenizing the noise variance in each column.
More precisely,
let $\bmpi_1,\dots,\bmpi_p \in \{0,1\}^{n \times n}$
be the $p$ PermPA permutation matrices
(drawn i.i.d. uniformly at random)
for a single trial,
and let
$\bmN_\pi \coloneqq [\bmpi_1 \bmN_{:1}; \cdots; \bmpi_p \bmN_{:p}] \in \bbR^{n \times p}$
be the resulting parallel noise matrix
generated by PermPA.
Then it immediately follows that for all $i,j$, 
\begin{align}
  \label{eq:perm:homvar}
  \bbE |(\bmN_\pi)_{ij}|^2
  &
  =
  \bbE |(\bmpi_j \bmN_{:j})_i|^2
  =
  \bbE |(\bmpi_j)_{i:} \bmN_{:j}|^2
  =
  \bbE \left[ \bmN_{:j}^\top (\bmpi_j)_{i:}^\top (\bmpi_j)_{i:} \bmN_{:j} \right]
  \\& \nonumber
  =
  \bbE \left[
    \bmN_{:j}^\top \,
    \bbE \left[ (\bmpi_j)_{i:}^\top (\bmpi_j)_{i:} \right]
    \bmN_{:j}
  \right]
  =
  \bbE \left[ \bmN_{:j}^\top \left(\frac{1}{n} I_n \right) \bmN_{:j} \right]
  =
  \frac{1}{n} \sum_{m=1}^n \bbE |N_{mj}|^2
  .
\end{align}
Namely,
the variances of all entries within a given column of $\bmN_\pi$ 
equal the
average of the variances in that column of $\bmN$.
PermPA homogenizes the noise variances,
and hence one might expect the parallel singular values
to look more like
those of a noise matrix with the homogenized variances \cref{eq:perm:homvar}
rather than those of the actual noise $\bmN$.

\begin{figure} \centering
  \includegraphics[scale=0.05]{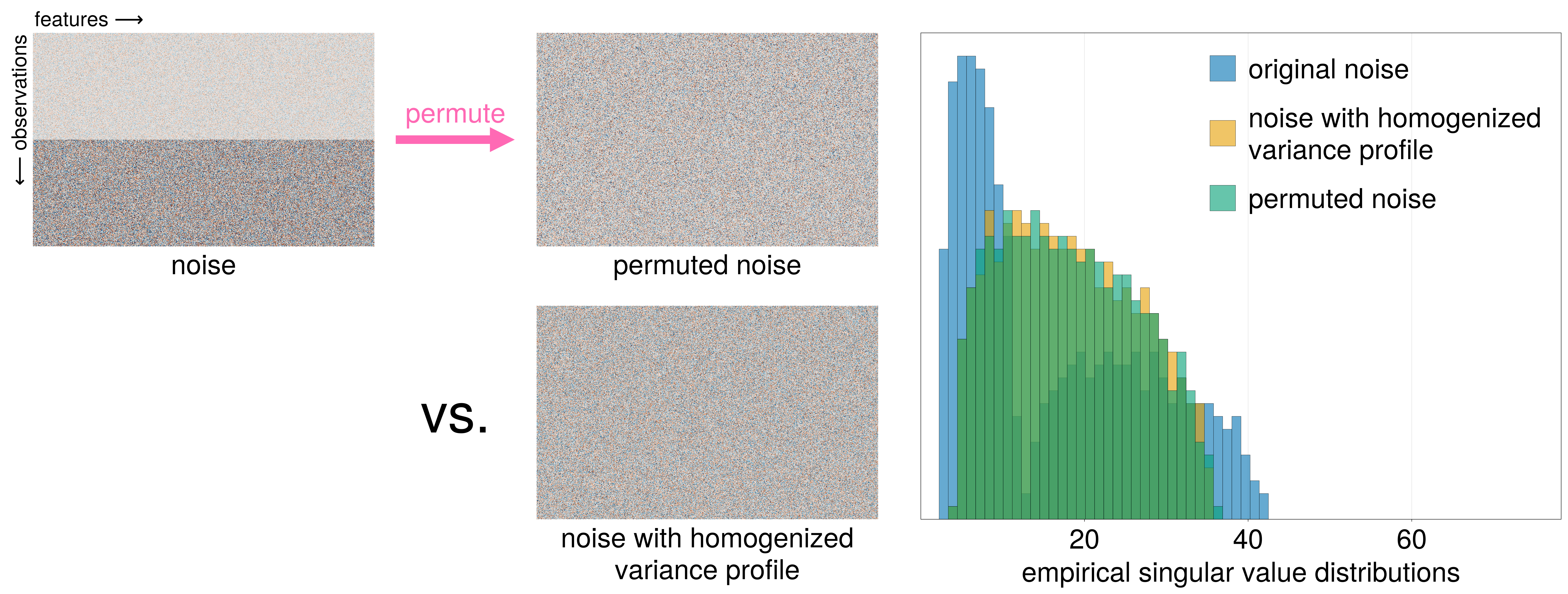}
  \caption{Example illustrating
    the homogenization effect of permutations
    on a $500 \times 800$ noise matrix
    with heterogeneous variances:
    the first half of the observations have a noise variance of $0.1$
    and the second half have a noise variance of $0.9$.
  }
  \label{fig:perm:homvar}
\end{figure}

Indeed, this distortion of the noise singular values
is exactly what can happen,
as shown in \cref{fig:perm:homvar}.
The empirical singular value distribution of the permuted noise
is very different from that of the original noise
and instead closely matches that of
a noise matrix with independent Gaussian entries
generated using the homogenized variance profile \cref{eq:perm:homvar}.
The consequence of such a distortion is that PermPA
incorrectly estimates the noise floor,
which can then lead to incorrect rank estimates.
For example, in \cref{fig:perm:homvar},
running PermPA directly on the original noise (which has zero signal components)
would incorrectly select many components
since the permuted noise has substantially smaller leading singular values
than the original noise.

The following \lcnamecref{thm:perm:homvar} makes
this intuition rigorous;
see \cref{thm:perm:homvar:proof} for the proof.

\begin{theorem}[PermPA homogenizes heterogeneous noise]
  \label{thm:perm:homvar}
  Let $\bmN \in \bbR^{n \times p}$ have independent entries
  with zero mean and uniformly bounded fourth moments,
  let $\bmN_\pi$ be the corresponding permuted noise,
  and let $\bbrN$ be the homogenized noise matrix defined as
  $\brN_{ij} \overset{ind}{\sim} \clN(0, v_j)$,
   where $v_j \coloneqq (1/n) \sum_{m=1}^n \bbE |N_{mj}|^2$.
     Then as $n,p \to \infty$ with $p/n \to \gamma > 0$,
  the empirical singular value distribution of $\bmN_\pi/\sqrt{n}$
  converges to the same almost sure weak limit
  as $\bbrN/\sqrt{n}$,
  as long as
  the empirical distribution of
  the column-wise variances $v_1,\dots,v_p$
  converges to some deterministic distribution $H$.
  In particular,
  their empirical singular value distributions both converge to the generalized Mar{\v{c}}enko-Pastur distribution
  defined by $H$.
\end{theorem}

Note that proving \cref{thm:perm:homvar} is more nontrivial than may initially appear.
Even though $\bmN_\pi$ has the same variance profile as $\bbrN$,
the entries of $\bmN_\pi$ are \emph{not independent};
one must show that the dependence induced by the permutations
has no impact in the limit.
This requires showing a generalized Mar{\v{c}}enko-Pastur law under relaxed independence conditions, related to, but different from those of \citet{hui2010lsd,wei2016tls,bryson2019marchenko}, 
and then showing that the permuted matrix satisfies these conditions by controlling the variance of certain quadratic forms.


\section{Simulation Studies}
\label{sec:simulations}

This section demonstrates the empirical performance of FlipPA
with numerical simulations.%
\if1\blind%
\footnote{Codes for reproducing the figures in this section are available online at:
\edit{\url{https://gitlab.com/dahong/rank-selection-via-random-signflips}}}
\else{ }\fi
We compare with the following
state-of-the-art methods
for high-dimensional data:
\begin{itemize}
  \item PermPA \citep{buja1992rop}:
  traditional permutation-based parallel analysis,
  which is highly effective for homoscedastic noise
  but not for heteroscedastic noise.
  We use $T = 100$ trials,
  a quantile of $\paquant = 1.0$,
  and a threshold of $\tau = 0$
  with the upper-edge comparison;
  this choice of comparison and quantile
  discourages overestimation.

  \item BCF and GIC \citep{bai2018consistency,hu2020lla:arxiv:v1}:
  the Bai-Choi-Fujikoshi (BCF) method
  studied in \citet{bai2018consistency}
  and the generalized information criterion (GIC)
  proposed by \citet{hu2020lla:arxiv:v1}
  are penalization-based model selection approaches;
  they each select the rank that minimizes
  an associated selection criterion.
  For BCF,
  we use the presentation
  from Equations~(4.5)--(4.6) of \cite{hu2020lla:arxiv:v1}.
  For GIC,
  we use Equation~(4.9) of \cite{hu2020lla:arxiv:v1}
  with $\gamma$ chosen as suggested in Remark 6 of \cite{hu2020lla:arxiv:v1}.

  \item ACT \citep{fan2020eno}:
  the adjusted correlation thresholding (ACT) method
  uses the eigenvalues of the correlation matrix
  instead of the covariance matrix
  to address observed variables with heterogeneous scales.
  The method has a maximum rank parameter, which we set to $100$.
  We use the implementation provided online
  in their supplementary materials:
  \url{https://www.tandfonline.com/doi/suppl/10.1080/01621459.2020.1825448}.

  \item BEMA0 and BEMA \citep{ke2021eot}:
  two methods based on
  bulk eigenvalue matching analysis (BEMA).
  BEMA0 is designed for homoscedastic noise;
  it estimates the noise variance
  by fitting the sample covariance matrix eigenvalues
  to a standard scaled Mar{\v{c}}enko-Pastur distribution,
  then uses the estimated noise variance
  to determine an appropriate cut-off
  from the Tracy-Widom distribution.
  BEMA is a generalization designed to
  allow for heteroscedasticity in the noise across features;
  instead of estimating a single noise variance for the data,
  it assumes that the feature-wise noise variances
  are drawn from a Gamma distribution
  and estimates the parameters of this Gamma distribution
  from the data
  to determine an appropriate cut-off.
  We use the R implementations (with the default values)
  provided online by the authors: \url{https://github.com/ZhengTracyKe/BEMA}.

  \item Biwhitening \citep[Section~4]{landa2021brt:arxiv:v2}:
  a method that addresses heteroscedasticity
  in nonnegative data matrices
  by rescaling the rows and columns of the data
  to whiten the noise along both axes;
  the necessary rescaling is determined
  via the Sinkhorn-Knopp algorithm (we use $1000$ iterations).
  The rank is estimated by
  counting how many singular values
  the biwhitened matrix has
  that exceed the homoscedastic cut-off of $\sqrt{n} + \sqrt{p}$.
  Since we consider data that may not be nonnegative,
  the main algorithms in \cite{landa2021brt:arxiv:v2}
  cannot be directly applied.
  We adapt them here by using the entrywise square of the data matrix
  as the variance estimate
  needed in Equation~(4.2) of \cite{landa2021brt:arxiv:v2}.
\end{itemize}
For FlipPA,
we use the same parameters as PermPA,
i.e., we use $T = 100$ trials,
a quantile of $\paquant = 1.0$,
and a threshold of $\tau = 0$.
For simplicity,
we focus on the upper-edge comparison method;
the results for the pairwise comparison method are similar.
Overall, we find that FlipPA compares favorably to the existing methods.

\subsection{Homogeneous noise variances}
\label{sec:sim:hom}

This section considers an $n \times p$ data matrix $\bmX$
with a rank-$1$ signal
generated as $\bmX = \theta \bmu \bmz^\top + \bmN$,
where $N_{ij} \overset{iid}{\sim} \clN(0, v/n)$,
$\bmu \in \bbR^n$ and $\bmz \in \bbR^p$
are drawn uniformly from the respective unit spheres,
$\theta$ is swept from zero to two,
and we take $v = 1$ without loss of generality.
As $\theta$ increases,
the signal transitions from being buried in the noise
to rising above it.
The noise is homogeneous,
so we expect the existing methods to also work well.

\begin{figure} \centering
  \includegraphics[scale=0.13]{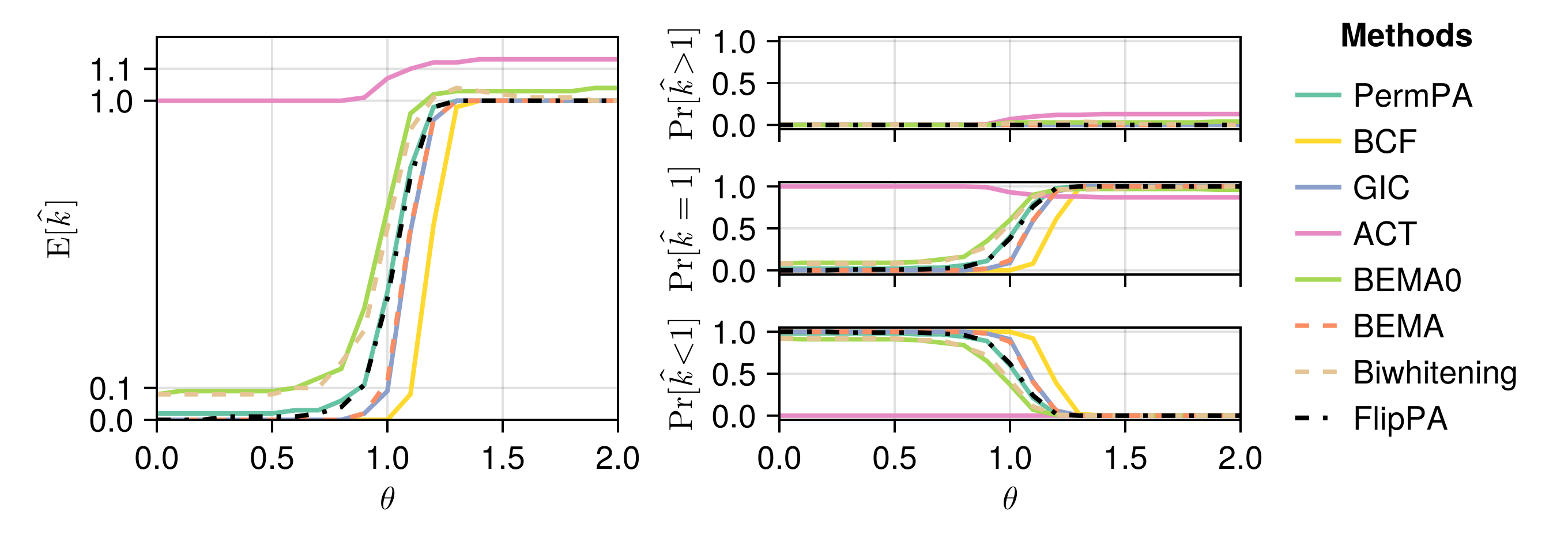}
  \caption{Performance across $100$ runs of each method for
    a rank-one signal in homogeneous noise,
    where the signal strength $\theta$ increases
    from zero (in the noise)
    to two (above the noise).}
  \label{fig:sim:hom}
\end{figure}

\Cref{fig:sim:hom} shows the resulting performance of each method
across $100$ runs,
where $n = 500$ and $p = 300$.
The left plot shows the average selected rank  $\bbE\bigl[\htk\bigr]$ across the runs;
the second column of plots shows the proportion $\Pr\bigl[\htk > 1\bigr]$
of runs resulting in over-estimation,
correct estimation (i.e., $\Pr\bigl[\htk = 1\bigr]$),
and under-estimation (i.e., $\Pr\bigl[\htk < 1\bigr]$).

As expected, all the methods were highly effective
  once $\theta$ was large enough,
  i.e., once the signal rose above the noise,
  and they estimated the rank correctly in most of the runs.
Interestingly, ACT also correctly estimated a rank of one
  for small $\theta$;
  using the eigenvalues of the correlation matrix may have helped here.
  However, ACT also incorrectly estimated a rank of one for $\theta = 0$,
  where the data is pure noise and the correct rank is zero.

The behavior at $\theta=0$ illustrates the type~I error properties of the methods
  since the data is pure noise in that case.
  PermPA rejected the pure noise null hypothesis
  in 2 of the 100 trials
  (achieving an empirical type~I error rate of 2\%);
  BEMA0 and Biwhitening rejected the null
  in 8 of the 100 trials
  (achieving empirical type~I error rates of 8\%);
  and ACT rejected the null
  in all of the 100 trials
  (achieving an empirical type~I error rate of 100\%).
  The other methods
  (BCF, GIC, BEMA, and FlipPA)
  never rejected the null,
  but FlipPA seemed to have the most power among them;
  for $\theta > 0$,
  FlipPA correctly rejected the null (and selected one component)
  more often than BCF, GIC, and BEMA.

We also considered a higher dimensional setting
($n = 60, p = 5000$)
in \cref{sec:sim:highdim};
\cref{fig:sim:highdim:hom} provides the analogue to \cref{fig:sim:hom}.
With the exceptions of BCF, GIC, and ACT
(all of which dramatically over-estimated the rank),
the remaining methods performed similarly:
PermPA, BEMA0, Biwhitening, and FlipPA
were all highly effective
once $\theta$ was large enough.
The empirical type~I error rates (seen at $\theta = 0$)
for this setting
were better for BEMA0 (4\%) and Biwhitening (3\%),
and were fairly similar for PermPA (1\%) and FlipPA (1\%).
They were dramatically worse for
BCF (100\%) and GIC (100\%),
matching ACT (100\%) here.
Only PermPA and FlipPA seem to have controlled the type~I error
for both data sizes here,
highlighting the rigorously guaranteed nonasymptotic type~I error control
of these randomization-based methods.
This benefit can make them appealing in practice.
Finally, we also considered \edit{settings with multiple signal components} in \cref{sec:sim:shadowing}
to investigate
the parallel analysis phenomenon known as ``shadowing'',
where strong signals can cause weak signals to be missed.
As expected, both parallel analysis methods (PermPA and FlipPA)
exhibited shadowing\edit{.
Notably, however, this only occurred when the strong signal was much stronger than the rest,
as predicted by the theoretical analysis in \cref{sec:relative:strength}}.

\subsection{Block-structured noise variance profiles}
\label{sec:sim:block}

To study how the various methods perform
for heterogeneous noise,
this section considers the same experiment as \cref{sec:sim:hom},
but now the noise matrix is generated
with a block-structured noise variance profile
as follows:
\begin{equation}
  \label{eq:sim:block:noise}
  N_{ij} \overset{ind}{\sim} \clN(0, V_{ij}/n)
  \quad \text{where} \quad
  \bmV
  =
  \begin{bmatrix}
    (1-\Delta) \cdot \bm1_{n_1 \times p_1} &
    1          \cdot \bm1_{n_1 \times p_2} \\
    (1+\Delta) \cdot \bm1_{n_2 \times p_1} &
    1          \cdot \bm1_{n_2 \times p_2}
  \end{bmatrix}
  ,
\end{equation}
and $\Delta$ is swept from zero to one.
As $\Delta$ increases from zero to one,
the noise transitions from being homoscedastic (when $\Delta = 0$)
to being increasingly heteroscedastic (as $\Delta$ grows).

\begin{figure} \centering
  \includegraphics[scale=0.13]{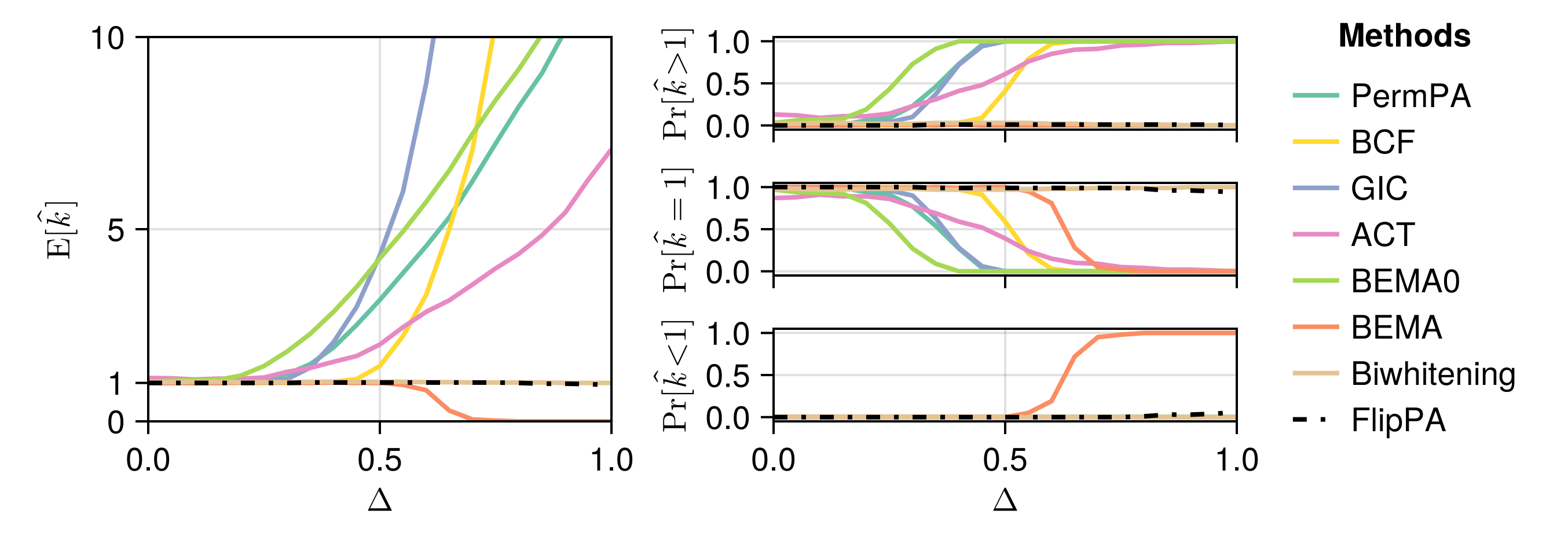}
  \caption{Performance across $100$ runs of each method for
    a rank-one signal in noise having a
    block-structured noise variance profile
    \cref{eq:sim:block:noise},
    where the amount of heteroscedasticity increases
    as $\Delta$ increases
    from zero (homoscedastic noise)
    to one (maximal heteroscedasticity).}
  \label{fig:sim:block}
\end{figure}

\Cref{fig:sim:block} shows the resulting performance of each method
across $100$ runs,
where $\theta = 1.5$, $n = 500$, $p = 300$, $n_1 = n_2 = 250$, $p_1 = 240$, and $p_2 = 60$.
As before, the left plot shows the average selected rank $\bbE\bigl[\htk\bigr]$ across the runs;
the second column of plots shows the proportion $\Pr\bigl[\htk > 1\bigr]$
of runs resulting in over-estimation,
correct estimation (i.e., $\Pr\bigl[\htk = 1\bigr]$),
and under-estimation (i.e., $\Pr\bigl[\htk < 1\bigr]$).

All methods were highly effective for small $\Delta$;
  they all properly accounted for the high-dimensionality of the data
  in this homoscedastic regime.
  However, only FlipPA and Biwhitening correctly estimated the rank across
  the entire sweep,
  with only a few runs where they either over-estimated or under-estimated.
   BEMA correctly estimated the rank until $\Delta \approx 0.5$;
  for $\Delta > 0.5$ it began to under-estimate the rank.
  The block-structured variance profile is outside of
  the assumptions under which BEMA was designed, since
  the noise variances are not uniform across the observations.
  The remaining methods significantly over-estimated the rank
  as $\Delta$ grew, because they do not properly account
  for the heteroscedasticity.

We also considered a higher dimensional setting
($n = 60, p = 5000$)
in \cref{sec:sim:highdim};
\cref{fig:sim:highdim:block} provides the analogue to \cref{fig:sim:block}.
With the exceptions of BCF, GIC, and ACT
(all of which dramatically over-estimated the rank),
the remaining methods again performed similarly:
PermPA, BEMA0, Biwhitening, and FlipPA
were all effective
when $\Delta$ was small,
but only FlipPA and Biwhitening remained effective
in the heteroscedastic case of larger~$\Delta$.
Finally, we also considered settings
with asymmetric noise distributions in \cref{sec:sim:asymmetric}.
Overall, the methods performed similarly
to the symmetric noise case considered here.

\subsection{Noise with dependent entries}
\label{sec:sim:block:dep}

To illustrate the flexibility of FlipPA,
here we consider the experiment from \cref{sec:sim:block}
but add blockwise dependence among the noise entries.
Namely,
suppose the data matrix $\bmX$ has a rank-one signal
as in \cref{sec:sim:hom}
but now the noise matrix $\bmN$
is generated with blockwise dependence
on top of the block-structured noise variance profile
of \cref{eq:sim:block:noise}.
Specifically, $\bmN = \bmSigma_1^{1/2} \btlN \bmSigma_2^{1/2}$,
where $\btlN$ has a block-structured variance profile
similar to \cref{eq:sim:block:noise}:
\begin{equation*}
  \tlN_{ij} \overset{ind}{\sim} \clN(0, V_{ij}/n)
  \quad \text{with} \quad
  \bmV  =
  \begin{bmatrix}
    0.25 \cdot \bm1_{n_1 \times p_1} &
    1    \cdot \bm1_{n_1 \times p_2} \\
    1.75 \cdot \bm1_{n_2 \times p_1} &
    1    \cdot \bm1_{n_2 \times p_2}
  \end{bmatrix}
  .
\end{equation*}
Further, $\bmSigma_1$ and $\bmSigma_2$ induce
blockwise dependence across both observations and features
as
\begin{align*}
  \bmSigma_1
  &
  =
  \bmI_{n/b_1} \otimes (\gamma \cdot \bm1_{b_1 \times b_1} + (1-\gamma) \cdot \bmI_{b_1})
  , &
  \bmSigma_2
  &
  =
  \bmI_{p/b_2} \otimes (\gamma \cdot \bm1_{b_2 \times b_2} + (1-\gamma) \cdot \bmI_{b_2})
  ,
\end{align*}
and $\gamma$ is swept from zero to $1/2$.
As $\gamma$ increases from zero to $1/2$,
the noise transitions from entrywise independence (when $\gamma = 0$)
to increasing blockwise dependence (as $\gamma$ grows),
i.e.,  the entries in each block of $b_1 \times b_2$ entries
become dependent.
Blockwise dependence can arise
in settings where observations and features may form groups,
e.g.,
blockwise dependence across both features and observations
can arise in genomics data
from a combination of
linkage disequilibrium
and
kinship among subjects \citep[see, e.g.,][]{conomos2016mfe}.
The blocks in these settings may be estimated
by one of various application-specific techniques developed in these domains.
We focus here on estimating the signal rank given the blocks.

In addition to the methods in the previous comparisons,
we add a straightforward variant of FlipPA
called BlockFlipPA
that flips the sign of all entries in a block together,
i.e.,
the signflip matrix $\bmR = \btlR \otimes \bm1_{b_1 \times b_2}$, where
  $\tlR_{ij} \overset{iid}{\sim}\pm1$, with probability 1/2.
This method is simply FlipPA applied at the block level%
\footnote{Note that the blocks of entries need not be regular and consecutive as in this case.}
so it preserves the distribution of the noise $\bmN$
and is a valid group-invariance based method,
just like FlipPA is for noise with independent entries.
Indeed,
BlockFlipPA has essentially the same guarantees
for noise with independent blocks
as FlipPA has for noise with independent entries;
see \cref{sec:blockflippa:consistency}.

\begin{figure} \centering
  \includegraphics[scale=0.13]{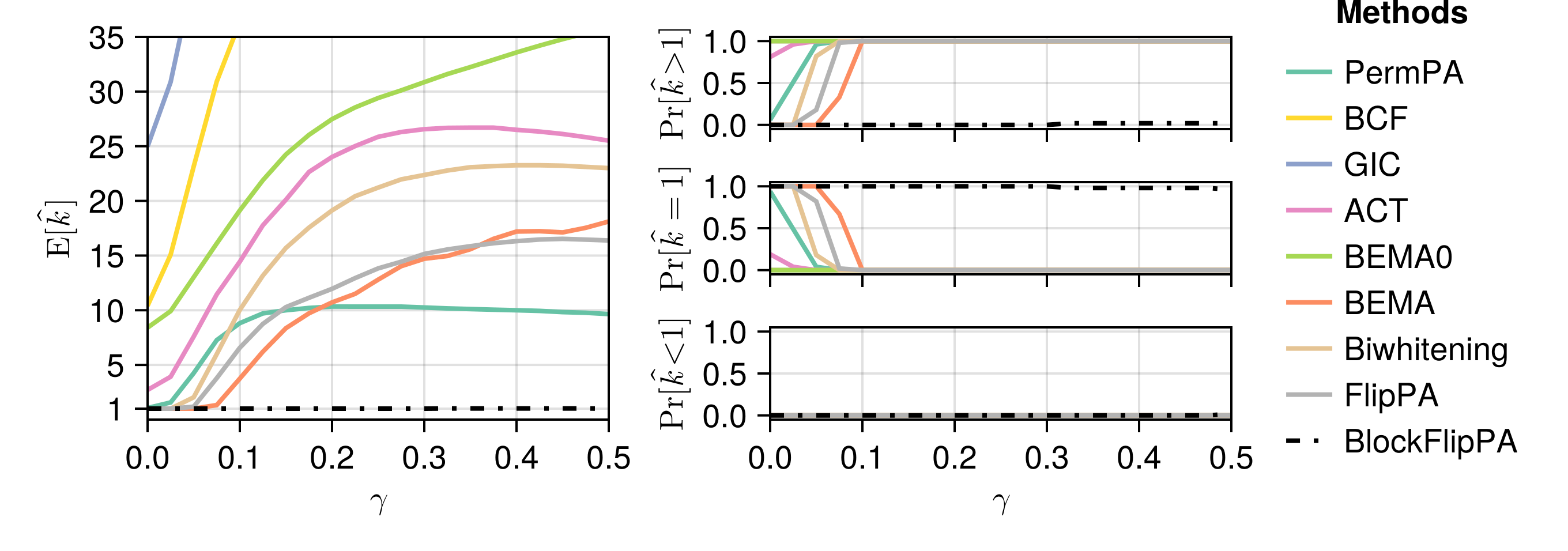}
  \caption{Performance across $100$ runs of each method for
    a rank-one signal in noise having a
    block-structured noise variance profile
    with blockwise dependence,
    where the amount of dependence increases
    as $\gamma$ increases
    from zero (independent entries)
    to $1/2$ (increasing blockwise dependence).}
  \label{fig:sim:block:dep}
\end{figure}

\Cref{fig:sim:block:dep} shows the resulting performance of each method
across $100$ runs,
where $\theta = 6$, $n = 500$, $p = 300$, $n_1 = n_2 = 250$, $p_1 = 240$, $p_2 = 60$, $b_1 = 25$, $b_2 = 15$.
As before,
the left plot shows the average selected rank $\bbE\bigl[\htk\bigr]$ across the runs;
the second column of plots shows the proportion $\Pr\bigl[\htk > 1\bigr]$
of runs resulting in over-estimation,
correct estimation (i.e., $\Pr\bigl[\htk = 1\bigr]$),
and under-estimation (i.e., $\Pr\bigl[\htk < 1\bigr]$).
When $\gamma$ is small,
  several methods are fairly effective.
However, only BlockFlipPA correctly estimates the rank across
  the entire sweep;
  the rest overestimate it.
Notably, this straightforward modification of FlipPA
allows it to remain effective
even in the presence of significant dependence,
highlighting the flexibility of the approach.

\begin{figure} \centering
  \includegraphics[scale=0.13]{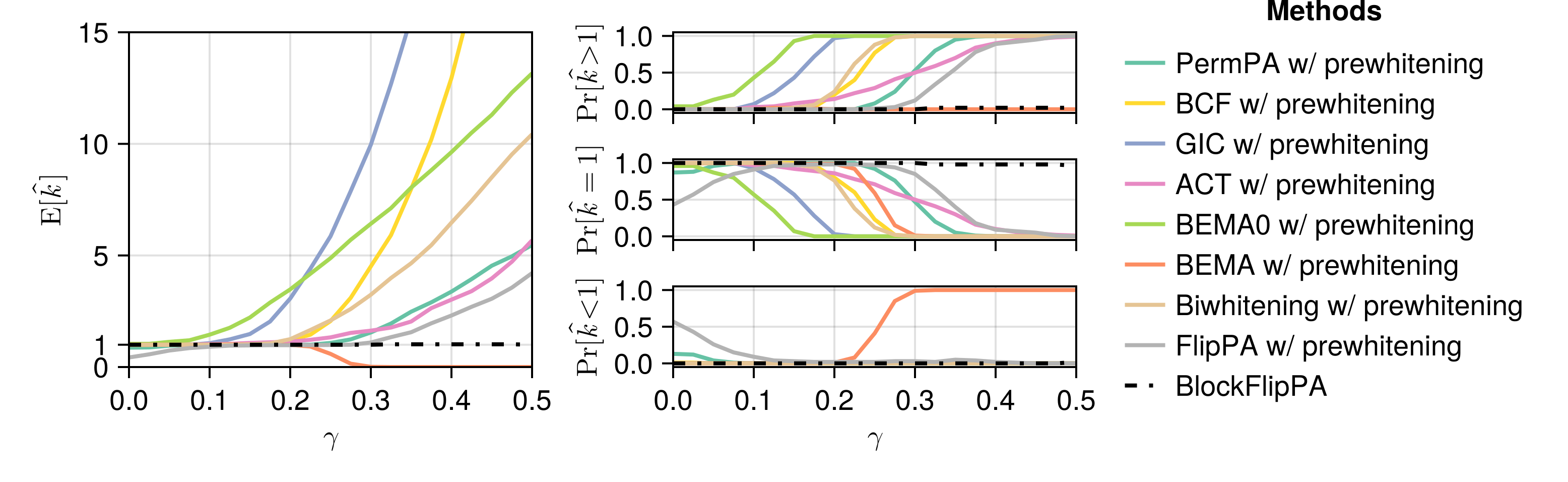}
  \caption{Performance
    for the same setting as \cref{fig:sim:block:dep}
    across $100$ runs of each method
    with prewhitening using $\bhtSigma_1,\bhtSigma_2$.}
  \label{fig:sim:block:dep:prewhiten:mask}
\end{figure}

It is perhaps less obvious how one would modify
the other methods to account for the dependence.
A common approach might be to instead prewhiten the data,
then apply one of the other methods,
in hopes that the prewhitening will sufficiently
remove the dependence.
For example,
one might form the estimates\footnote{
Note that the estimates $\bhtSigma_1,\bhtSigma_2$
use the same knowledge of the blockwise dependence structure
as BlockFlipPA;
they do so by zeroing out all entries
except those corresponding to within-block dependence.
We also considered forming the estimates $\bhtSigma_1$ and $\bhtSigma_2$
without zeroing out these entries;
this form of prewhitening did not generally perform better.}
$\bhtSigma_1
  =
  (\bmI_{n/b_1} \otimes \bm1_{b_1 \times b_1})
  \circ
  (\bmX \bmX^\top)$, 
$\bhtSigma_2
=
  (\bmI_{p/b_2} \otimes \bm1_{b_2 \times b_2})
  \circ
  (\bmX^\top \bmX)$,
then apply the methods to the prewhitened data
$\btlX = \bhtSigma_1^{-1/2} \bmX \bhtSigma_2^{-1/2}$.
Doing so improved their performance at small to moderate $\gamma$,
as shown in \cref{fig:sim:block:dep:prewhiten:mask},
but only to a point.
BlockFlipPA was again the only method
that correctly estimated the rank
across the entire sweep.
For $\gamma > 0.3$,
BEMA with prewhitening under-estimated the rank
and the other methods with prewhitening over-estimated the rank.

We also considered a higher dimensional setting
($n = 60, p = 5000$)
in \cref{sec:sim:highdim};
\cref{fig:sim:highdim:block:dep,fig:sim:highdim:block:dep:prewhiten:mask}
provide the analogues to
\cref{fig:sim:block:dep,fig:sim:block:dep:prewhiten:mask}.
Similar to the lower dimensional setting here,
only BlockFlipPA was effective across the entire sweep.
Moreover, prewhitening did not significantly improve the methods,
perhaps due to poor estimates of the prewhitening operators
in the higher dimensional setting.
To summarize,
this experiment demonstrates
the ease of adapting FlipPA to allow for dependence among the entries
and highlights the flexibility of the approach.


\section{Empirical Data Illustration from Astronomy}
\label{sec:exp}
\label{sec:exp:quasar:spectra}

This section illustrates the behavior of FlipPA on an empirical dataset
coming from astronomy.%
\if1\blind%
\footnote{Codes for reproducing the figures in this section are available online at:
\edit{\url{https://gitlab.com/dahong/rank-selection-via-random-signflips}}}
\else{ }\fi
We consider
quasar spectra obtained
from the 16th data release of the Sloan Digital Sky Survey (SDSS) \citep{ahumada2020t1d}
using the associated DR16Q quasar catalog \citep{lyke2020tsd}.
The spectra are each composed of flux measurements at various wavelengths,
and each spectrum comes with associated variance estimates for the flux measurements.
To illustrate FlipPA and the methods considered in \cref{sec:simulations},
we selected and preprocessed (i.e., interpolated, centered, and normalized)
a subset of the data
in a manner similar to Section~8 of \cite{hong2023owp};
see \cref{sec:exp:quasar:spectra:preprocess} for details.

\begin{figure} \centering
  \hfill
  \begin{subfigure}[t]{0.55\linewidth} \centering
    \includegraphics[scale=0.22]{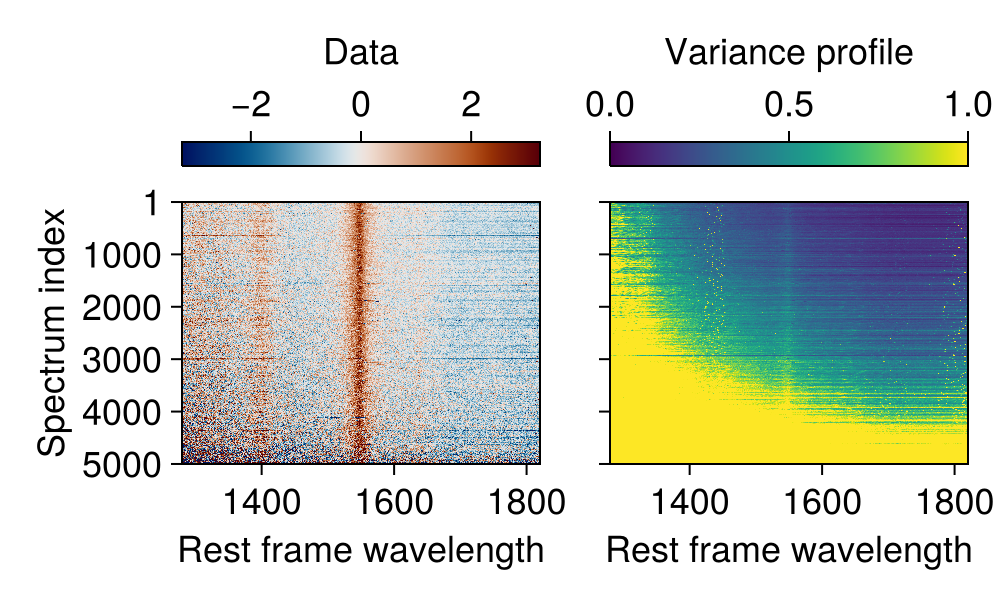}
    \caption{Data $\bmY \in \bbR^{n \times p}$ and variance profile $\bmV \in \bbR^{n \times p}$.}
    \label{fig:exp:quasar:spectra:data}
  \end{subfigure}
  \hfill
  \begin{subfigure}[t]{0.4\linewidth} \centering
    \includegraphics[scale=0.22]{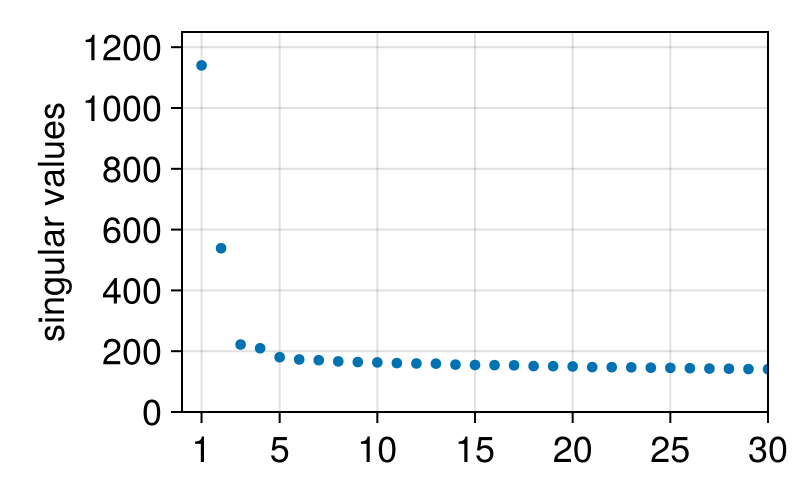}
    \caption{Scree plot for the data matrix $\bmY$.}
    \label{fig:exp:quasar:spectra:scree}
  \end{subfigure}
  \hfill\strut
  \caption{Quasar spectra dataset with associated variance profile
    and scree plot.}
  \label{fig:exp:quasar:spectra:data:scree}
\end{figure}

\Cref{fig:exp:quasar:spectra:data} shows the resulting
data matrix $\bmY \in \bbR^{5000 \times 1081}$
of $5000$ spectra across $1081$ wavelengths,
as well as the corresponding variance profiles
$\bmV \in \bbR^{5000 \times 1081}$.
As can be seen from the variance profile,
the data is heteroscedastic
along both observations (here, spectra) and features (here, wavelengths).
\Cref{fig:exp:quasar:spectra:scree}
shows the associated scree plot,
i.e., a scatter plot of the singular values of $\bmY$.

\begin{table} \centering
  \caption{Rank selections for the quasar spectroscopy data}
  \label{tbl:exp:quasar:spectra}
  \begin{tabular}{l|c@{\hspace{0.2\linewidth}}l|c}
    PermPA (pairwise)            & 133   &  BEMA0               & 168 \\
    PermPA (upper-edge)          & 1     &  BEMA                & 4   \\
    BCF                          & 1080  &  Biwhitening         & 7   \\
    GIC                          & 1080  &  FlipPA (pairwise)   & 4   \\
    ACT ($\text{max rank}=1000$) & 107   &  FlipPA (upper-edge) & 2
  \end{tabular}
\end{table}

\Cref{tbl:exp:quasar:spectra}
shows the ranks selected
by FlipPA and the methods considered in \cref{sec:simulations}.
For PermPA and FlipPA,
we include both the pairwise comparison and upper-edge comparison options.
For ACT, we set the max rank parameter here to be $1000$
because it selected up to the max rank when we set it to $100$.

We observe
that PermPA with pairwise comparisons, BCF, GIC, ACT, and BEMA0
    all selected over 100 components for this dataset.
    Comparing with the scree plot in \cref{fig:exp:quasar:spectra:scree},
    such large selections seem unlikely to be correct.
    Notably, BCF and GIC selected all the components.
In stark contrast to PermPA with pairwise comparisons,
    PermPA with upper-edge comparisons selected only one component
    for this dataset.
    The reason is likely that permutations were not able to sufficiently
    suppress the largest singular value,
    leading to shadowing of smaller components in this case.

The remaining methods are roughly similar,
    and seem fairly reasonable when compared with the scree plot
    in \cref{fig:exp:quasar:spectra:scree}.
    The selection of 2 components by FlipPA with upper-edge comparisons
    corresponds to selecting the two strong components that are substantially above the rest,
    while the selection of 4 components by BEMA and FlipPA with pairwise comparisons
    corresponds to including the next two weaker components
    that nonetheless appear to rise slightly above the rest.
    Biwhitening selects several more components,
    which is likely because the biwhitening operation
    changes the data singular values
    and can cause components that were previously ``buried''
    in the noise to rise above it.

\begin{figure} \centering
  \includegraphics[scale=0.25]{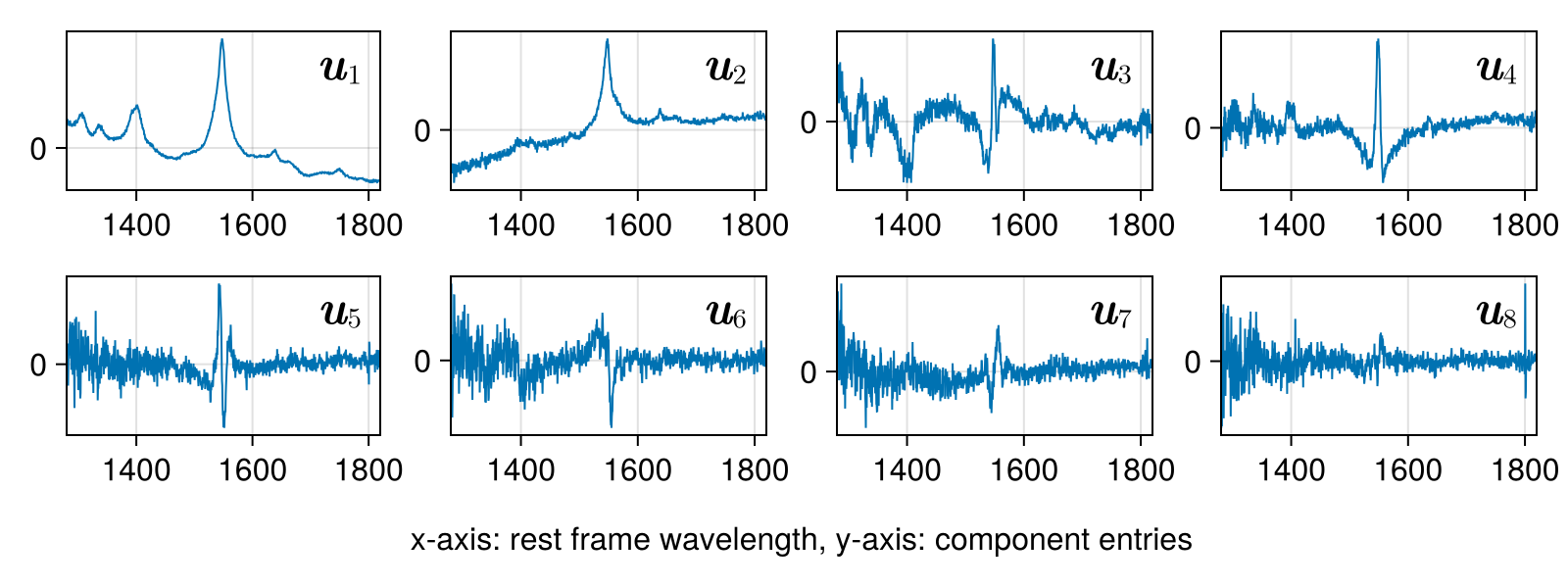}
  \caption{Leading eigenspectra computed from the data matrix $\bmY$.}
  \label{fig:exp:quasar:spectra:eigenspectra}
\end{figure}

 \Cref{fig:exp:quasar:spectra:eigenspectra} shows eigenspectra computed
    by taking the leading right singular vectors of the data matrix $\bmY$.
    The rank selections by FlipPA, BEMA, and Biwhitening
    seem to correspond to relatively ``clean'' eigenspectra,
    providing some further evidence that they are
    selecting meaningful signal components
    rather than components buried in the noise.
Overall, this example
illustrates the potential benefits of FlipPA
for rank selection in practice.


\section{Discussion}
\label{sec:discuss}

This paper proposed a new method (FlipPA)
for tackling the problem of rank estimation
from data with heterogeneous noise.
In future work, there are opportunities to sharpen the theoretical analysis of FlipPA;
for example,
we suspect that the delocalization conditions in \cref{sec:consistency}
can likely be weakened.
Moreover, it could be interesting to better understand the robustness of FlipPA
to violations of the assumptions studied here,
especially when they are rare (e.g., due to outliers).
Another unexplored aspect
is the interplay of the matrix size
and the number of trials to run.
When the matrix is large,
it can be computationally expensive to carry out many trials,
but the singular values also tend to be less variable,
so FlipPA may not require many trials.

\section*{Acknowledgements}

The authors thank
Johannes Alt,
Andreas Buja,
Johannes Heiny,
Boris Landa,
Matthew McKay,
and
Jeffrey Fessler
for helpful and stimulating discussions.
They are also grateful to seminar participants at
RMCDAW 2019
and
Berkeley
for feedback on early versions.
The quasar spectra were provided by the Sloan Digital Sky Survey (SDSS);
see \cref{ack:sdss}.
\edit{AI tools were used to help develop some of the results; see \cref{ack:ai}.}

\section*{Funding}

This work was partially supported by
NSF BIGDATA grant IIS 1837992, 
NSF CAREER grant DMS 2046874,
the Dean's Fund for Postdoctoral Research of the Wharton School,
and
NSF Mathematical Sciences Postdoctoral Research Fellowship DMS 2103353.

\bibliographystyle{plainnat-abbrev}
\bibliography{references}

\newpage
\appendix


\section{Proof of \cref{prop:type-I}}
\label{thm:false:positive:proof}

We use the same overall approach as the proof of \cite[Theorem 2]{hemerik2017etw}.
Note first that
for both pairwise and upper-edge comparison methods
\begin{align}
    \label{eq:false:positive:iff}
    \htk > 0
    \quad \iff \quad
    &
    \sigma_{1} > \text{$\paquant$-quantile of } \left(
        \tlsigma_1^{(1)},\dots,\tlsigma_1^{(T)}
    \right)
    + \tau
    \\ \nonumber
    \quad \implies \quad
    &
    \sigma_{1} > \text{$\paquant$-quantile of } \left(
        \tlsigma_1^{(1)},\dots,\tlsigma_1^{(T)}
    \right)
    .
\end{align}
Now, whenever
$\sigma_{1} > \text{$\paquant$-quantile of } ( \tlsigma_1^{(1)},\dots,\tlsigma_1^{(T)} )$
then it must also be the case that
\begin{align}
    \label{eq:false:positive:implies}
    &
    \#\left\{t \in \{1,\dots,T\} : \tlsigma_1^{(t)} < \sigma_1 \right\}
    \\&\qquad \nonumber
    \geq
    \#\left\{t \in \{1,\dots,T\} :
        \tlsigma_1^{(t)} \leq \text{$\paquant$-quantile of } \left(
            \tlsigma_1^{(1)},\dots,\tlsigma_1^{(T)}
        \right)
    \right\}
    \geq \lfloor \paquant T \rfloor
    ,
\end{align}
where the second inequality holds
for sample quantiles
computed using any of the nine commonly used definitions
enumerated in \cite{hyndman1996sqi}.%
\footnote{Specifically,
the second inequality corresponds to Property P2 in \cite{hyndman1996sqi},
with the correction posted online by the authors \citep{hyndman1996sqi:errata}.}
Moreover, note that
\begin{equation}
    \label{eq:false:positive:expand}
    \left\{t \in \{1,\dots,T\} : \tlsigma_1^{(t)} < \sigma_1 \right\}
    =
    \left\{t \in \{0,\dots,T\} : \sigma_1(\bmR^{(t)} \circ \bmX) < \sigma_1(\bmX) \right\}
\end{equation}
where $\bmR^{(0)} = \bm1_{n \times p}$
and $\bmR^{(1)},\dots,\bmR^{(T)} \overset{iid}{\sim} \operatorname{Unif}(\{-1, 1\}^{n \times p})$
are the $T$ signflip matrices in FlipPA,
since $\sigma_1 = \sigma_1(\bmX)$
and $\tlsigma_1^{(t)} = \sigma_1(\bmR^{(t)} \circ \bmX)$ for $t = 1,\dots,T$.
Combining \cref{eq:false:positive:iff,eq:false:positive:implies,eq:false:positive:expand}
yields
\begin{align}
    \label{eq:false:positive:bound}
    &
    \Pr[\htk > 0]
    =
    \Pr\left[
        \sigma_{1} > \text{$\paquant$-quantile of } \left(
            \tlsigma_1^{(1)},\dots,\tlsigma_1^{(T)}
        \right)
        + \tau
    \right]
    \\&\quad \nonumber
    \leq
    \Pr\left[
        \sigma_{1} > \text{$\paquant$-quantile of } \left(
            \tlsigma_1^{(1)},\dots,\tlsigma_1^{(T)}
        \right)
    \right]
    \\&\quad \nonumber
    \leq
    \Pr\left[
        \#\left\{t \in \{1,\dots,T\} : \tlsigma_1^{(t)} < \sigma_1 \right\}
        \geq \lfloor \paquant T \rfloor
    \right]
    \\&\quad \nonumber
    =
    \Pr\left[
        \#\left\{t \in \{0,\dots,T\} : \sigma_1(\bmR^{(t)} \circ \bmX) < \sigma_1(\bmX) \right\}
        \geq \lfloor \paquant T \rfloor
    \right]
    ,
\end{align}
where these probabilities
are with respect to
the randomness in
$\bmX$ and $\bmR^{(1)},\dots,\bmR^{(T)}$,
all under the null hypothesis $H_0$ of no signal.

So, it remains to bound the final probability in \cref{eq:false:positive:bound},
for which we will exploit the invariance properties of the signflips.
Namely, note that for any $j \in \{0,\dots,T\}$,
\begin{equation}
    \label{eq:false:positive:exchangability}
    \{\{ \bmR^{(0)}, \dots, \bmR^{(T)} \}\} =_d \{\{ \bmR^{(0)} \circ \bmR^{(j)}, \dots, \bmR^{(T)} \circ \bmR^{(j)} \}\}
    ,
\end{equation}
where the double braces here denote a multiset,
i.e., a set with multiplicities for its elements.
Additionally,
\cref{assump:noise:sym:ind} implies that
under $H_0$, for any $j \in \{0,\dots,T\}$,
\begin{equation}
    \label{eq:false:positive:invariance}
    \bmX =_d \bmR^{(j)} \circ \bmX
    .
\end{equation}
Thus, for any $j \in \{0,\dots,T\}$,
\begin{align}
    \label{eq:false:positive:symmetric:form}
    &
    \Pr\left[
        \#\left\{t \in \{0,\dots,T\} : \sigma_1(\bmR^{(t)} \circ \bmX) < \sigma_1(\bmX) \right\}
        \geq \lfloor \paquant T \rfloor
    \right]
    \\&\quad \nonumber
    =
    \Pr\left[
        \#\left\{t \in \{0,\dots,T\} : \sigma_1(\bmR^{(t)} \circ \bmR^{(j)} \circ \bmX) < \sigma_1(\bmX) \right\}
        \geq \lfloor \paquant T \rfloor
    \right]
    \\&\quad \nonumber
    =
    \Pr\left[
        \#\left\{t \in \{0,\dots,T\} : \sigma_1(\bmR^{(t)} \circ \bmR^{(j)} \circ \bmR^{(j)} \circ \bmX) < \sigma_1(\bmR^{(j)} \circ \bmX) \right\}
        \geq \lfloor \paquant T \rfloor
    \right]
    \\&\quad \nonumber
    =
    \Pr\left[
        \#\left\{t \in \{0,\dots,T\} : \sigma_1(\bmR^{(t)} \circ \bmX) < \sigma_1(\bmR^{(j)} \circ \bmX) \right\}
        \geq \lfloor \paquant T \rfloor
    \right]
    ,
\end{align}
where the first equality uses \cref{eq:false:positive:exchangability},
the second equality uses \cref{eq:false:positive:invariance},
and the final equality follows by simplifying.
Thus,
averaging \cref{eq:false:positive:symmetric:form}
across $j \in \{0,\dots,T\}$ yields
\begin{align}
    \label{eq:false:positive:averaged}
    &
    \Pr\left[
        \#\left\{t \in \{0,\dots,T\} : \sigma_1(\bmR^{(t)} \circ \bmX) < \sigma_1(\bmX) \right\}
        \geq \lfloor \paquant T \rfloor
    \right]
    \\&\quad \nonumber
    =
    \frac{1}{T+1}
    \sum_{j = 0}^T
    \Pr\left[
        \#\left\{t \in \{0,\dots,T\} : \sigma_1(\bmR^{(t)} \circ \bmX) < \sigma_1(\bmR^{(j)} \circ \bmX) \right\}
        \geq \lfloor \paquant T \rfloor
    \right]
    \\&\quad \nonumber
    =
    \frac{1}{T+1}
    \sum_{j = 0}^T
    \bbE\left[
        \bm1 \left(
            \#\left\{t \in \{0,\dots,T\} : \sigma_1(\bmR^{(t)} \circ \bmX) < \sigma_1(\bmR^{(j)} \circ \bmX) \right\}
            \geq \lfloor \paquant T \rfloor
        \right)
    \right]
    \\&\quad \nonumber
    =
    \frac{1}{T+1} \;
    \bbE\left[
        \sum_{j = 0}^T
        \bm1 \left(
            \#\left\{t \in \{0,\dots,T\} : \sigma_1(\bmR^{(t)} \circ \bmX) < \sigma_1(\bmR^{(j)} \circ \bmX) \right\}
            \geq \lfloor \paquant T \rfloor
        \right)
    \right]
    \\&\quad \nonumber
    =
    \frac{1}{T+1} \;
    \bbE \left[
    \;\;
        \#\Big\{
    \;
            j \in \{0,\dots,T\} :
    \right.
    \\&\qquad\qquad\qquad\qquad\qquad \nonumber
    \left.
            \#\left\{t \in \{0,\dots,T\} : \sigma_1(\bmR^{(t)} \circ \bmX) < \sigma_1(\bmR^{(j)} \circ \bmX) \right\}
            \geq \lfloor \paquant T \rfloor
    \;
        \Big\}
    \;\;
    \right]
    \\&\quad \nonumber
    \leq
    \frac{1}{T+1} \;
    \bbE [ (T + 1) - \lfloor \paquant T \rfloor ]
    = 1 - \frac{\lfloor \paquant T \rfloor}{T+1}
    ,
\end{align}
where $\bm1(\cdot)$ denotes the indicator function,
which equals one when the argument is true and zero otherwise.
Combining \cref{eq:false:positive:bound,eq:false:positive:averaged}
concludes the proof.
\qed

\section{Proofs for \cref{sec:consistency}}
\label{sec:consistency:proofs}

This \lcnamecref{sec:consistency:proofs}
provides proofs for the consistency guarantees in \cref{sec:consistency}.

\subsection{Useful lemmas}

This section states some lemmas that will be useful
for proving the results from \cref{sec:consistency}.
The first provides a simple expression for
the $\|\cdot\|_{2,\infty}$ norm of a rank-one matrix.
\begin{lemma}[$\|\cdot\|_{2,\infty}$ norm of a rank-one matrix]
    \label{lem:2infty:rankone}
    Let $\bmu \in \bbR^n$ and $\bmz \in \bbR^p$ be arbitrary.
    Then $\| \bmu \bmz^\top \|_{2,\infty} = \|\bmu\|_\infty \|\bmz\|_2$.
\end{lemma}

The proof of \cref{lem:2infty:rankone} is straightforward;
we provide it here for the reader's convenience.
\begin{proof}[Proof of \cref{lem:2infty:rankone}]
    Since the $i$-th row of $\bmu \bmz^\top$
    is $(\bmu \bmz^\top)_{i:} = u_i \cdot \bmz$,
    we immediately have
    \begin{equation*}
        \left\|\bmu \bmz^\top \right\|_{2,\infty}
        = \max_{i=1,\dots,n} \left\|(\bmu \bmz^\top)_{i:} \right\|_{2}
        = \max_{i=1,\dots,n} \left\|u_i \cdot \bmz \right\|_{2}
        = \max_{i=1,\dots,n} | u_i | \cdot \left\| \bmz \right\|_{2}
        = \| \bmu \|_\infty \| \bmz \|_{2}
        ,
    \end{equation*}
    which completes the proof.
\end{proof}

The next lemma provides an elegant relationship between
the entrywise matrix norm $\|\cdot\|_{t,t}$
and a corresponding maximum row norm $\|\cdot\|_{t/2,\infty}$.
\begin{lemma}[Bound for $\|\cdot\|_{t,t}$ norm in terms of $\|\cdot\|_{t/2,\infty}$]
    \label{lem:entrywise:norm:columnwise}
    Let $\bmA \in \bbR^{n \times p}$ and $t \geq 2$ be an arbitrary even number.
    Then $\|\bmA\|_{t,t} \leq \rank^{1/2}(\bmA) \sqrt{\|\bmA\|_{t/2,\infty} \cdot \|\bmA^\top\|_{t/2,\infty}}$.
\end{lemma}
\begin{proof}[Proof of \cref{lem:entrywise:norm:columnwise}]
    Note first that
    \begin{equation*}
        \|\bmA\|_{t,t}^t
        = \sum_{i=1}^n \sum_{j=1}^p \left| A_{i,j} \right|^t
        = \sum_{i=1}^n \sum_{j=1}^p \left| A_{i,j}^{t/2} \right|^2
        = \|\bmB\|_F^2
        ,
    \end{equation*}
    where $\|\cdot\|_F$ denotes the Frobenius norm
    and $\bmB \in \bbR^{n \times p}$ has entries $B_{i,j} \coloneqq A_{i,j}^{t/2}$,
    i.e.,
    \begin{equation*}
        \bmB \coloneqq \underbrace{\bmA \circ \cdots \circ \bmA}_{\text{$t/2$ times}}
        ,
    \end{equation*}
    where $\circ$ denotes the Hadamard product.
    Now, recall that for any matrix $\bmB$,
    \begin{align*}
        \|\bmB\|_F^2 &\leq \rank(\bmB) \|\bmB\|_2^2
        , &
        \|\bmB\|_2 \leq \sqrt{\|\bmB\|_1 \|\bmB\|_\infty}
        ,
    \end{align*}
    where $\|\cdot\|_2$, $\|\cdot\|_1$, and $\|\cdot\|_\infty$
    denote the $\ell_2$, $\ell_1$, and $\ell_\infty$
    operator norms, respectively.
    Thus, we have
    \begin{equation}
        \label{eq:entrywise:norm:columnwise:bound1}
        \|\bmA\|_{t,t}^t
        = \|\bmB\|_F^2
        \leq \rank(\bmB) \|\bmB\|_1 \|\bmB\|_\infty
        .
    \end{equation}
    Next, note that
    \begin{align}
        \label{eq:entrywise:norm:columnwise:l1}
        \|\bmB\|_1
        &
        = \max_{j=1,\dots,p} \, \sum_{i=1}^n |B_{ij}|
        = \max_{j=1,\dots,p} \, \sum_{i=1}^n \left| A_{ij}^{t/2} \right|
        = \max_{j=1,\dots,p} \, \sum_{i=1}^n |A_{ij}|^{t/2}
        = \|\bmA^\top\|_{t/2, \infty}^{t/2}
        , \\
        \label{eq:entrywise:norm:columnwise:linfty}
        \|\bmB\|_\infty
        &
        = \max_{i=1,\dots,n} \, \sum_{j=1}^p |B_{ij}|
        = \max_{i=1,\dots,n} \, \sum_{j=1}^p \left| A_{ij}^{t/2} \right|
        = \max_{i=1,\dots,n} \, \sum_{j=1}^p |A_{ij}|^{t/2}
        = \|\bmA\|_{t/2, \infty}^{t/2}
        ,
    \end{align}
    and
    \begin{equation}
        \label{eq:entrywise:norm:columnwise:rank}
        \rank(\bmB) = \rank(\,\underbracket{\bmA \circ \cdots \circ \bmA}_{\text{$t/2$ times}}\,)
        \leq [\rank(\bmA)]^{t/2}
        .
    \end{equation}
    Substituting \cref{eq:entrywise:norm:columnwise:l1,eq:entrywise:norm:columnwise:linfty,eq:entrywise:norm:columnwise:rank}
    into \cref{eq:entrywise:norm:columnwise:bound1}
    and taking the $t$-th root
    completes the proof.
\end{proof}

The final lemma collects and combines some recent results
on the operator norms of nonhomogeneous random matrices
\citep{seginer2000ten, latala2005seo, schuett2013ote, bandeira2016sharp, van2017structured, latala2018tdf}
and writes the result in our notation.
\begin{lemma}[Operator norms of signflipped matrices]
    \label{lem:opnorm:bound}
    Let $\bmA \in \bbR^{n \times p}$ be arbitrary
    and let $\bmR \sim \operatorname{Unif}(\{-1, 1\}^{n \times p})$
    be a Rademacher random matrix.
    Then
    $\bbE \|\bmA \circ \bmR\|
     \lesssim
     \|\bmA\|_{2,\infty}
     +
     \|\bmA^\top\|_{2,\infty}
     +
     \min\big[ \rho_2(\bmA), \rho_\infty(\bmA) \big]$.
\end{lemma}

\begin{proof}[Proof of \cref{lem:opnorm:bound}]
    Consider the two symmetric $(n+p) \times (n+p)$ matrices
    \begin{align*}
        \btlA &\coloneqq \begin{bmatrix} & \bmA \\ \bmA^\top & \end{bmatrix}
        , &
        \btlR &\coloneqq \begin{bmatrix} & \bmR \\ \bmR^\top & \end{bmatrix}
        .
    \end{align*}
    Several existing results provide
    dimension-dependent bounds for this case,
    e.g., writing \cite[Corollary 4.7]{bandeira2016sharp}
    in our notation and simplifying yields
    \begin{equation}
        \label{eq:opnorm:bound:dimdep}
        \bbE \|\btlA \circ \btlR\|
        \lesssim
        \| \btlA \|_{2,\infty}
        \sqrt[4]{\log (n+p)}
        =
        \begin{Vmatrix} & \bmA \\ \bmA^\top & \end{Vmatrix}_{2,\infty}
        \sqrt[4]{\log (n+p)}
        =
        \rho_2(\bmA)
        .
    \end{equation}

    Next, we develop a bound using the dimension-free bounds of \cite{latala2018tdf},
    which apply to nonhomogeneous Gaussian matrices.
    For this,
    let $\btlG$ be a symmetric  $(n+p) \times (n+p)$ Gaussian random matrix,
    i.e., $\tlG_{ij} \overset{iid}{\sim} \clN(0,1)$ for $i \geq j$.
    Then writing \cite[Theorem 1.1]{latala2018tdf} in our notation
    and simplifying
    yields the following bound for $\btlA \circ \btlG$:
    \begin{align}
        \label{eq:opnorm:bound:dimfree:gauss}
        \bbE \|\btlA \circ \btlG\|
        &
        \lesssim
        \| \btlA \|_{2,\infty}
        +
        \max_{i=1,\dots,n+p} \| \btlA \|_{\infty,(i)} \sqrt{\log i}
        \\& \nonumber
        =
        \begin{Vmatrix} & \bmA \\ \bmA^\top & \end{Vmatrix}_{2,\infty}
        +
        \max_{i=1,\dots,n+p} \begin{Vmatrix} & \bmA \\ \bmA^\top & \end{Vmatrix}_{\infty,(i)} \sqrt{\log i}
        \\& \nonumber
        =
        \max\left(
            \| \bmA \|_{2,\infty}, \| \bmA^\top \|_{2,\infty}
        \right)
        +
        \rho_\infty(\bmA)
        .
    \end{align}
    We next use the following
    relationship between $\btlA \circ \btlG$ and $\btlA \circ \btlR$
    (also used in \cite{latala2005seo}):
    \begin{align}
        \label{eq:opnorm:bound:gaussrad}
        \bbE \|\btlA \circ \btlG\|
        &
        =
        \bbE \big\| \btlA \circ \btlR \circ |\btlG| \big\|
        =
        \bbE_{\btlR} \left[ \bbE_{\btlG} \big\| \btlA \circ \btlR \circ |\btlG| \big\| \right]
        \\& \nonumber
        \geq
        \bbE_{\btlR} \left\| \bbE_{\btlG} \big[ \btlA \circ \btlR \circ |\btlG| \big] \right\|
        =
        \bbE_{\btlR} \left\| \btlA \circ \btlR \circ \bbE_{\btlG} \big[ |\btlG| \big] \right\|
        =
        \sqrt{\frac{2}{\pi}} \;
        \bbE \|\btlA \circ \btlR\|
        ,
    \end{align}
    where $|\btlG| \in \bbR^{(n+p) \times (n+p)}$ is the elementwise absolute value of $\btlG$,
    i.e., $|\btlG|_{ij} = |\tlG_{ij}|$.
    The first equality holds because
    $\btlA \circ \btlG =_d \btlA \circ \btlR \circ |\btlG|$,
    the second equality is the law of total expectation,
    the inequality follows by applying Jensen's inequality to the operator norm $\|\cdot\|$,
    and the final two equalities hold by linearity of the expectation and the fact that $\bbE|\tlG_{ij}| = \sqrt{2/\pi}$.

    Combining \cref{eq:opnorm:bound:dimfree:gauss,eq:opnorm:bound:gaussrad}
    yields the following dimension-free bound for $\btlA \circ \btlR$:
    \begin{equation}
        \label{eq:opnorm:bound:dimfree}
        \bbE \|\btlA \circ \btlR\|
        \lesssim
        \bbE \|\btlA \circ \btlG\|
        \lesssim
        \max\left(
            \| \bmA \|_{2,\infty}, \| \bmA^\top \|_{2,\infty}
        \right)
        +
        \rho_\infty(\bmA)
        ,
    \end{equation}
    and finally combining \cref{eq:opnorm:bound:dimdep,eq:opnorm:bound:dimfree}
    and simplifying yields
    \begin{align*}
        \bbE \|\btlA \circ \btlR\|
        &
        \lesssim
        \min\Big[
            \rho_2(\bmA),
            \max\left(
                \| \bmA \|_{2,\infty}, \| \bmA^\top \|_{2,\infty}
            \right)
            +
            \rho_\infty(\bmA)
        \Big]
        \\&
        \leq
        \| \bmA \|_{2,\infty} + \| \bmA^\top \|_{2,\infty}
        + \min\big[ \rho_2(\bmA), \rho_\infty(\bmA) \big]
        .
    \end{align*}
    Since the nonzero eigenvalues of $\btlA \circ \btlR$
    consist of positive and negative copies of the singular values of $\bmA \circ \bmR$
    \citep[see, e.g.,][Theorem 4.2]{stewart1990mpt},
    it follows that $\bbE \|\bmA \circ \bmR\| = \bbE \|\btlA \circ \btlR\|$,
    which concludes the proof.
\end{proof}

\subsection{Proof of \cref{thm:consistency:component}}
\label{thm:consistency:component:proof}

\Cref{thm:consistency:component} can be proved using similar techniques
as in the proof of \cref{thm:consistency},
but doing so can make it less obvious
that \cref{thm:consistency:component}
is actually a special case of \cref{thm:consistency}.
So, here we will prove \cref{thm:consistency:component}
by showing that it follows from \cref{thm:consistency}.
Namely, we show that
the signal matrix
$\bmS = \sum_{i = 1}^{k} \theta_i \bmu_i \bmz_i^\top$
satisfies the conditions of \cref{thm:consistency},
i.e.,
\begin{equation*}
    \bbE \| \bmS \|_{2,\infty} \to 0
    , \quad
    \bbE \| \bmS^{\top} \|_{2,\infty} \to 0
    , \quad
    \text{and} \quad
    \min\big[ \bbE \, \rho_2(\bmS), \bbE \, \rho_\infty(\bmS) \big] \to 0
    .
\end{equation*}
We begin with the first two conditions.
Note that
\begin{align}
    \label{eq:consistency:component:proof:S}
    \bbE \| \bmS \|_{2,\infty}
    &
    =
    \bbE \left\| \sum_{i = 1}^{k} \theta_i \bmu_i \bmz_i^\top \right\|_{2,\infty}
    \leq
    \bbE \left\{
        \sum_{i = 1}^{k} \left\|\theta_i \bmu_i \bmz_i^\top \right\|_{2,\infty}
    \right\}
    =
    \bbE \left\{
        \sum_{i = 1}^{k} \theta_i \|\bmu_i\|_\infty \|\bmz_i\|_2
    \right\}
    , \\
    \label{eq:consistency:component:proof:St}
    \bbE \| \bmS^{\top} \|_{2,\infty}
    &
    =
    \bbE \left\| \sum_{i = 1}^{k} \theta_i \bmz_i \bmu_i^\top \right\|_{2,\infty}
    \leq
    \bbE \left\{
        \sum_{i = 1}^{k} \left\|\theta_i \bmz_i \bmu_i^\top \right\|_{2,\infty}
    \right\}
    =
    \bbE \left\{
        \sum_{i = 1}^{k} \theta_i \|\bmu_i\|_2 \|\bmz_i\|_\infty
    \right\}
    .
\end{align}
where the inequality in each line follows from the triangle inequality
and the final equality in each line follows from \cref{lem:2infty:rankone}.
Thus, we have
\begin{align}
    \label{eq:consistency:component:proof:results:1:2}
    &
    \max\left(
        \bbE \| \bmS \|_{2,\infty},
        \bbE \| \bmS^{\top} \|_{2,\infty}
    \right)
    \leq
    \bbE \| \bmS \|_{2,\infty} + \bbE \| \bmS^{\top} \|_{2,\infty}
    \\& \nonumber \qquad
    \leq
    \bbE \left\{
        \sum_{i = 1}^{k} \theta_i \|\bmu_i\|_\infty \|\bmz_i\|_2
    \right\}
    +
    \bbE \left\{
        \sum_{i = 1}^{k} \theta_i \|\bmu_i\|_2 \|\bmz_i\|_\infty
    \right\}
    \\& \nonumber \qquad
    =
    2 \cdot
    \bbE \left\{
        \sum_{i = 1}^{k}
        \theta_i \|\bmu_i\|_2 \|\bmz_i\|_2
        \cdot
        \left[
            \frac{\|\bmu_i\|_\infty/\|\bmu_i\|_2 + \|\bmz_i\|_\infty/\|\bmz_i\|_2}{2}
        \right]
    \right\}
    \to 0
    ,
\end{align}
where
the first inequality holds because $\bbE \| \bmS \|_{2,\infty}$ and $\bbE \| \bmS^{\top} \|_{2,\infty}$ are nonnegative,
the second inequality follows from \cref{eq:consistency:component:proof:S,eq:consistency:component:proof:St},
and the final limit follows from the delocalization condition \cref{eq:consistency:component}
in the theorem.
So, we are done with the first two conditions of \cref{thm:consistency}.

To show that the final condition of \cref{thm:consistency} is satisfied,
we will show that $\bbE \|\bmS\|_{4,4} \to 0$ then apply \cref{thm:decay:suff}.
Note first that
\begin{align*}
    \bbE \| \bmS \|_{4,4}
    &
    =
    \bbE \left\| \sum_{i = 1}^{k} \theta_i \bmu_i \bmz_i^\top \right\|_{4,4}
    \leq
    \bbE \left\{
        \sum_{i = 1}^{k} \left\| \theta_i \bmu_i \bmz_i^\top \right\|_{4,4}
    \right\}
    \\&
    \leq
    \bbE \left\{
        \sum_{i = 1}^{k}
        \rank^{1/2}( \theta_i \bmu_i \bmz_i^\top )
        \sqrt{
            \| \theta_i \bmu_i \bmz_i^\top \|_{2,\infty}
            \cdot
            \| \theta_i \bmz_i \bmu_i^\top \|_{2,\infty}
        }
    \right\}
    ,
\end{align*}
where the first inequality follows from the triangle inequality,
and the second inequality follows from \cref{lem:entrywise:norm:columnwise}.
Noting that $\rank( \theta_i \bmu_i \bmz_i^\top ) = 1$, applying \cref{lem:2infty:rankone}, and simplifying
then yields
\begin{align*}
    \bbE \| \bmS \|_{4,4}
    &
    =
    \bbE \left\{
        \sum_{i = 1}^{k}
        \theta_i \|\bmu_i\|_2 \|\bmz_i\|_2
        \cdot
        \sqrt{
            \frac{\|\bmu_i\|_\infty}{\|\bmu_i\|_2}
            \cdot
            \frac{\|\bmz_i\|_\infty}{\|\bmz_i\|_2}
        }
    \right\}
    \\&
    \leq
    \bbE \left\{
        \sum_{i = 1}^{k}
        \theta_i \|\bmu_i\|_2 \|\bmz_i\|_2
        \cdot
        \left[
            \frac{\|\bmu_i\|_\infty/\|\bmu_i\|_2 + \|\bmz_i\|_\infty/\|\bmz_i\|_2}{2}
        \right]
    \right\}
    \to 0
    ,
\end{align*}
where the inequality follows from the AM-GM inequality,
and the limit follows from the delocalization condition \cref{eq:consistency:component}
in the theorem.
Finally, applying \cref{thm:decay:suff} yields
\begin{equation}
    \label{eq:consistency:component:proof:results:3}
    \min\big[ \bbE \, \rho_2(\bmS), \bbE \, \rho_\infty(\bmS) \big]
    \leq
    \bbE \, \rho_\infty(\bmS)
    \to 0
    ,
\end{equation}
so we are done with the final condition of \cref{thm:consistency}.

The proof concludes by combining \cref{thm:consistency}
with \cref{eq:consistency:component:proof:results:1:2,eq:consistency:component:proof:results:3}.
\qed

\subsection{Proof of \cref{thm:consistency:component:rates}}
\label{thm:consistency:component:rates:proof}

We prove \cref{thm:consistency:component:rates}
by showing that 
the signal matrix
$\bmS = \sum_{i = 1}^{k} \theta_i \bmu_i \bmz_i^\top$
satisfies the delocalization condition \cref{eq:consistency:component}
from \cref{thm:consistency:component}.
Note first that
\begin{align*}
    &
    \bbE \left\{
        \sum_{i = 1}^{k}
        \theta_i \|\bmu_i\|_2 \|\bmz_i\|_2
        \cdot
        \left[
            \frac{\|\bmu_i\|_\infty/\|\bmu_i\|_2 + \|\bmz_i\|_\infty/\|\bmz_i\|_2}{2}
        \right]
    \right\}
    \\&
    =
    \bbE \left\{
        \sum_{i = 1}^{k}
        \theta_i
        \left[
            \frac{\|\bmu_i\|_\infty \|\bmz_i\|_2 + \|\bmu_i\|_2 \|\bmz_i\|_\infty}{2}
        \right]
    \right\}
    =
    \sum_{i = 1}^{k}
    \theta_i
    \left[
        \frac{\bbE \|\bmu_i\|_\infty \bbE \|\bmz_i\|_2 + \bbE \|\bmu_i\|_2 \bbE \|\bmz_i\|_\infty}{2}
    \right]
    ,
\end{align*}
because the signal rank and component coefficients are deterministic
and the component vectors are jointly independent.
Substituting the rates yields
\begin{align*}
    &
    \sum_{i = 1}^{k}
    \theta_i
    \left[
        \frac{\bbE \|\bmu_i\|_\infty \bbE \|\bmz_i\|_2 + \bbE \|\bmu_i\|_2 \bbE \|\bmz_i\|_\infty}{2}
    \right]
    \\& \qquad
    \lesssim
    \sum_{i = 1}^{k}
    m^{\beta_1} \log^{\beta_2} m
    \left[
        \frac{n^{-\alpha_1} \log^{-\alpha_2} n + p^{-\alpha_1} \log^{-\alpha_2} p}{2}
    \right]
    \\& \qquad
    =
    k
    (m^{\beta_1} \log^{\beta_2} m)
    \left[
        \frac{n^{-\alpha_1} \log^{-\alpha_2} n + p^{-\alpha_1} \log^{-\alpha_2} p}{2}
    \right]
    \\& \qquad
    \lesssim
    (m^{\nu_1}\log^{\nu_2}m)
    (m^{\beta_1} \log^{\beta_2} m)
    \left[
        \frac{n^{-\alpha_1} \log^{-\alpha_2} n + p^{-\alpha_1} \log^{-\alpha_2} p}{2}
    \right]
    \\& \qquad
    =
    \frac{1}{2} \Big[
        (m^{\nu_1 + \beta_1}\log^{\nu_2 + \beta_2}m) (n^{-\alpha_1} \log^{-\alpha_2} n)
        +
        (m^{\nu_1 + \beta_1}\log^{\nu_2 + \beta_2}m) (p^{-\alpha_1} \log^{-\alpha_2} p)
    \Big]
    .
\end{align*}
Now, since $m = \min(n,p) \geq 1$, $\nu_1+\beta_1 \geq 0$, and $\nu_2+\beta_2 \geq 0$,
it follows that
\begin{align*}
    &
    (m^{\nu_1 + \beta_1}\log^{\nu_2 + \beta_2}m) (n^{-\alpha_1} \log^{-\alpha_2} n)
    \\& \qquad
    \leq
    (n^{\nu_1 + \beta_1}\log^{\nu_2 + \beta_2}n) (n^{-\alpha_1} \log^{-\alpha_2} n)
    =
    n^{\nu_1 + \beta_1 - \alpha_1}\log^{\nu_2 + \beta_2 - \alpha_2}n
    , \\
    &
    (m^{\nu_1 + \beta_1}\log^{\nu_2 + \beta_2}m) (p^{-\alpha_1} \log^{-\alpha_2} p)
    \\& \qquad
    \leq
    (p^{\nu_1 + \beta_1}\log^{\nu_2 + \beta_2}p) (p^{-\alpha_1} \log^{-\alpha_2} p)
    =
    p^{\nu_1 + \beta_1 - \alpha_1}\log^{\nu_2 + \beta_2 - \alpha_2}p
    .
\end{align*}
Thus we have
\begin{align*}
    &
    \bbE \left\{
        \sum_{i = 1}^{k}
        \theta_i \|\bmu_i\|_2 \|\bmz_i\|_2
        \cdot
        \left[
            \frac{\|\bmu_i\|_\infty/\|\bmu_i\|_2 + \|\bmz_i\|_\infty/\|\bmz_i\|_2}{2}
        \right]
    \right\}
    \\& \qquad
    \lesssim
    n^{\nu_1 + \beta_1 - \alpha_1}\log^{\nu_2 + \beta_2 - \alpha_2}n
    +
    p^{\nu_1 + \beta_1 - \alpha_1} \log^{\nu_2 + \beta_2-\alpha_2} p
    ,
\end{align*}%
which converges to zero
if either
$\alpha_1 > \nu_1 + \beta_1$,
or
$\alpha_1 = \nu_1 + \beta_1$ and $\alpha_2 > \nu_2+\beta_2$.
\qed

\subsection{Proof of \cref{thm:consistency}}
\label{pf:thm:consistency}

By \cref{assump:signal:perceptible},
there exists $\ep > 0$
for which
$\Pr[ \sigma_k > \| \bmN \| + \ep ] \to 1$.
Now note that for the upper-edge comparison method
\begin{equation*}
    \htk = k
    \quad \iff \quad
    \sigma_{k+1}
    \leq
    \text{$\paquant$-quantile of } \left(\tlsigma_1^{(1)},\dots,\tlsigma_1^{(T)}\right)
    + \tau
    <
    \sigma_k
    ,
\end{equation*}
which occurs for $\tau \in (0,\ep)$ whenever
the following conditions are simultaneously met:
\begin{gather*}
    \sigma_k > \| \bmN \| + \ep
    \quad \text{and} \quad
    \sigma_{k+1} \leq \| \bmN \|
    \\
    \text{and} \quad
    \| \bmN \| - \tau \leq \tlsigma_1^{(t)} \leq \| \bmN \| + (\ep-\tau)
    \text{ for $t = 1,\dots,T$}
    ,
\end{gather*}
since in that case we have
\begin{alignat*}{3}
    \sigma_{k+1}
    &
    \leq \|\bmN\|
    &&
    \leq \min&&\left(\tlsigma_1^{(1)},\dots,\tlsigma_1^{(T)}\right) + \tau
    \leq \text{$\paquant$-quantile of } \left(\tlsigma_1^{(1)},\dots,\tlsigma_1^{(T)}\right) + \tau
    ,
    \\
    \sigma_{k}
    &
    > \|\bmN\| + \ep
    &&
    \geq \max&&\left(\tlsigma_1^{(1)},\dots,\tlsigma_1^{(T)}\right) + \tau
    \geq \text{$\paquant$-quantile of } \left(\tlsigma_1^{(1)},\dots,\tlsigma_1^{(T)}\right) + \tau
    .
\end{alignat*}
As a result, we have the simple bound
\begin{align}
    \label{eq:consistency:unionbound}
    \Pr[\htk = k]
    &
    \geq
    \Pr\left[
        \begin{gathered}
            \sigma_k > \| \bmN \| + \ep \; \text{and} \; \sigma_{k+1} \leq \| \bmN \|
            \\
            \text{and} \;
            \| \bmN \| - \tau \leq \tlsigma_1^{(t)} \leq \| \bmN \| + (\ep-\tau)
            \text{ for $t = 1,\dots,T$}
        \end{gathered}
    \right]
    \\& \nonumber
    \geq
    1
    +
    \left(
        \Pr[\sigma_k > \| \bmN \| + \ep] - 1
    \right)
    +
    \left(
        \Pr[\sigma_{k+1} \leq \| \bmN \|] - 1
    \right)
    \\&\quad \nonumber
    +
    \sum_{t=1}^{T} \left(
        \Pr[\| \bmN \| - \tau \leq \tlsigma_1^{(t)} \leq \| \bmN \| + (\ep-\tau)] - 1
    \right)
    ,
\end{align}
where the second inequality comes from applying a union bound to the complement of the event.
Now, by \cref{assump:signal:perceptible},
we already have
\begin{equation}
    \label{eq:consistency:limit1}
    \Pr[ \sigma_k > \| \bmN \| + \ep ] \to 1
    .
\end{equation}
Moreover, it follows from Weyl's inequality
\citep[see, e.g.,][Corollary 4.10]{stewart1990mpt}
that
$
\sigma_{k+1}
= \sigma_{k+1}(\bmS + \bmN)
\leq \sigma_{k+1}(\bmS) + \| \bmN \|
= \| \bmN \|
$
since $\rank(\bmS) \leq k$,
so
\begin{equation}
    \label{eq:consistency:limit1b}
    \Pr[ \sigma_{k+1} \leq \| \bmN \| ] = 1
    .
\end{equation}
Hence,
it remains to show that
$\Pr[\| \bmN \| - \tau \leq \tlsigma_1^{(t)} \leq \| \bmN \| + (\ep-\tau)] \to 1$
for $t = 1,\dots,T$.
Namely, we need to show that the signflipped data singular value $\tlsigma_1^{(t)}$
is close to the noise operator norm $\| \bmN \|$.
To do so, we first decompose the difference of the two
using triangle inequality:
\begin{align}
    \label{eq:consistency:split}
    \big| \tlsigma_1^{(t)} - \| \bmN \| \big|
    &
    = \big| \sigma_1(\bmR^{(t)} \circ \bmX) - \sigma_1(\bmN) \big|
    \\& \nonumber
    \leq
    \underbrace{
        \big| \sigma_1(\bmR^{(t)} \circ \bmX) - \sigma_1(\bmR^{(t)} \circ \bmN) \big|
    }_{\text{signal destruction}}
    +
    \underbrace{
        \big| \sigma_1(\bmR^{(t)} \circ \bmN) - \sigma_1(\bmN) \big|
    }_{\text{noise preservation}}
    ,
\end{align}
where $\bmR^{(t)} \sim \operatorname{Unif}(\{-1, 1\}^{n \times p})$
is the signflip matrix for trial $t$.
The first term captures destruction of the signal by signflipping
(so that the signflipped data looks like the signflipped noise)
and the second term captures preservation of the noise
(so that the signflipped noise looks like the true noise).
We now analyze each term separately,
starting with the noise preservation:
\begin{enumerate}
    \item \textbf{noise preservation $\big| \sigma_1(\bmR^{(t)} \circ \bmN) - \sigma_1(\bmN) \big|$:}
    For this term,
    note that
    \begin{equation}
        \label{eq:consistency:noiseflipnoise}
        \big| \sigma_1(\bmR^{(t)} \circ \bmN) - \sigma_1(\bmN) \big|
        \leq
        \underbrace{
        \big| \sigma_1(\bmR^{(t)} \circ \bmN) - \brsigma \big|
        }_{=_d \; | \sigma_1(\bmN) - \brsigma |}
        +
        \big|\brsigma - \sigma_1(\bmN) \big|
        \ipto 0
        ,
    \end{equation}
    because $\sigma_1(\bmN) = \| \bmN \| \ipto \brsigma$ by \cref{assump:noise:upper:edge}
    and $\bmR^{(t)} \circ \bmN =_d \bmN$ by \cref{assump:noise:sym:ind}.

    \item \textbf{signal destruction $\big| \sigma_1(\bmR^{(t)} \circ \bmX) - \sigma_1(\bmR^{(t)} \circ \bmN) \big|$:}
    For this term,
    note that Weyl's inequality
    \citep[see, e.g.,][Corollary 4.10]{stewart1990mpt}
    yields
    \begin{equation*}
        \big| \sigma_1(\bmR^{(t)} \circ \bmX) - \sigma_1(\bmR^{(t)} \circ \bmN) \big|
        \leq
        \| \bmR^{(t)} \circ \bmS \|
        ,
    \end{equation*}
    since $\bmR^{(t)} \circ \bmX = \bmR^{(t)} \circ (\bmS + \bmN) = \bmR^{(t)} \circ \bmS + \bmR^{(t)} \circ \bmN$.
    Now, using \cref{lem:opnorm:bound} we have
    \begin{align*}
        &
        \bbE \big| \sigma_1(\bmR^{(t)} \circ \bmX) - \sigma_1(\bmR^{(t)} \circ \bmN) \big|
        \leq
        \bbE \| \bmR^{(t)} \circ \bmS \|
        \\&\qquad
        \lesssim
        \bbE \|\bmS\|_{2,\infty}
        +
        \bbE \|\bmS^\top\|_{2,\infty}
        +
        \min\big[ \bbE \, \rho_2(\bmS), \bbE \, \rho_\infty(\bmS) \big]
        \to 0
        ,
    \end{align*}
    since each of these terms go to zero by assumption;
    note that the expectations in the final line are only with respect to $\bmS$.
    So, we have shown that
    $\big| \sigma_1(\bmR^{(t)} \circ \bmX) - \sigma_1(\bmR^{(t)} \circ \bmN) \big|$ converges to zero in $L^1$,
    and thus we have
    \begin{equation}
        \label{eq:consistency:flipdataflipnoise}
        \big| \sigma_1(\bmR^{(t)} \circ \bmX) - \sigma_1(\bmR^{(t)} \circ \bmN) \big| \ipto 0
        .
    \end{equation}
\end{enumerate}
Combining \cref{eq:consistency:split,eq:consistency:noiseflipnoise,eq:consistency:flipdataflipnoise}
yields
\begin{equation*}
    \big| \tlsigma_1^{(t)} - \| \bmN \| \big|
    \leq
    \big| \sigma_1(\bmR^{(t)} \circ \bmX) - \sigma_1(\bmR^{(t)} \circ \bmN) \big|
    +
    \big| \sigma_1(\bmR^{(t)} \circ \bmN) - \sigma_1(\bmN) \big|
    \ipto 0
    ,
\end{equation*}
i.e., $\tlsigma_1^{(t)} \ipto \| \bmN \|$,
and so for all $t = 1,\dots,T$
\begin{equation}
    \label{eq:consistency:limit2}
    \Pr[\| \bmN \| - \tau \leq \tlsigma_1^{(t)} \leq \| \bmN \| + (\ep-\tau)] \to 1
    .
\end{equation}
The proof concludes by combining \cref{eq:consistency:unionbound,eq:consistency:limit1,eq:consistency:limit1b,eq:consistency:limit2}.
\qed

\subsection{Proof of \cref{thm:decay:suff}}
\label{thm:decay:suff:proof}

Borrowing an argument from \cite{latala2018tdf},
note that for any matrix $\bmA \in \bbR^{m \times \ell}$,
$t \geq 2$,
and $i \in \{1,\dots,m\}$,
we have
\begin{equation*}
    i \cdot \|\bmA\|_{\infty,(i)}^t
    \leq
    \sum_{j=1}^i \|\bmA\|_{\infty,(j)}^t
    \leq
    \sum_{j=1}^{m} \|\bmA\|_{\infty,(j)}^t
    =
    \sum_{j=1}^{m} \|\bmA_{j:}\|_{\infty}^t
    \leq \sum_{j=1}^{m} \sum_{s=1}^{\ell} |A_{js}|^t
    = \|\bmA\|_{t,t}^t
    ,
\end{equation*}
so $\|\bmA\|_{\infty,(i)} \leq i^{-1/t} \, \|\bmA\|_{t,t}$.
Thus, for any $t \geq 2$
\begin{align*}
    \rho_\infty(\bmS)
    &
    =
    \max_{i=1,\dots,n+p} \begin{Vmatrix} & \bmS \\ \bmS^\top & \end{Vmatrix}_{\infty,(i)} \sqrt{\log i}
    \leq
    \max_{i=1,\dots,n+p} \; i^{-1/t} \, \begin{Vmatrix} & \bmS \\ \bmS^\top & \end{Vmatrix}_{t,t} \sqrt{\log i}
    \\&
    =
    \left[ \max_{i=1,\dots,n+p} \frac{\sqrt{\log i}}{i^{1/t}} \right]
    \begin{Vmatrix} & \bmS \\ \bmS^\top & \end{Vmatrix}_{t,t}
    =
    \left[ \max_{i=1,\dots,n+p} \frac{\sqrt{\log i}}{i^{1/t}} \right]
    2^{1/t}
    \,
    \|\bmS\|_{t,t}
    \\&
    \leq
    \underbrace{
        \left[ \max_{i=1,\dots,\infty} \frac{\sqrt{\log i}}{i^{1/t}} \right]
    }_{< \infty}
    2^{1/t}
    \,
    \|\bmS\|_{t,t}
    \lesssim \|\bmS\|_{t,t}
    ,
\end{align*}
and so
$\bbE \, \rho_\infty(\bmS) \to 0$
if $\bbE \|\bmS\|_{t,t} \to 0$,
which establishes the first condition of \cref{thm:decay:suff}.
Setting $t = 4$ and applying \cref{lem:entrywise:norm:columnwise}
yields
\begin{equation}
    \label{eq:rhoinf:bound}
    \rho_\infty(\bmS)
    \lesssim \|\bmS\|_{4,4}
    \leq
    \rank^{1/2}(\bmS) \sqrt{\|\bmS\|_{2,\infty} \cdot \|\bmS^\top\|_{2,\infty}}
    ,
\end{equation}
so
$\bbE \, \rho_\infty(\bmS) \to 0$
if
$\bbE \left\{\rank^{1/2}(\bmS) \sqrt{\|\bmS\|_{2,\infty} \cdot \|\bmS^\top\|_{2,\infty}} \right\} \to 0$,
which establishes the second condition of \cref{thm:decay:suff}
and completes the proof.
\qed

\edit{
\section{Consistency with strong signal components}
\label{sec:relative:strength}

\Cref{sec:consistency} established the asymptotic consistency of FlipPA in settings where signflipping causes the signal to vanish in operator norm relative to the noise level.
In these settings, FlipPA produces accurate estimates of the asymptotic noise upper-edge and selects the correct rank.
This includes settings with weak ``emergent'' components that produce signal singular values of the same order as those coming from the noise.
Moreover, as shown in \cref{thm:consistency:component:rates}, this also includes settings where the signal singular values diverge at rates up to $\max_{i=1,\dots,k} \, \theta_i = O(m^{\beta_1} \log^{\beta_2} m)$
with
$\beta_1 < 1/2$
or $\beta_1 = 1/2$ with $\beta_2 < -1/2$,
where $m=\min(n,p)$.
It does not, however, include settings with stronger signals whose singular values diverge at rates of $\max_{i=1,\dots,k} \, \theta_i \asymp \sqrt{np}$ as studied, e.g., in \citet{bai2002determining}.

This \lcnamecref{sec:relative:strength} considers settings with strong components that do not vanish after signflipping.
The key is that even though these components are not completely destroyed after signflipping, they are still significantly reduced.
Hence, weaker signals are still correctly selected by FlipPA as long as they are sufficiently strong relative to the strong signals.
Indeed, \cref{sec:sim:shadowing} demonstrates this phenomenon in numerical experiments.
Here, we provide a precise theoretical characterization.
To simplify the presentation, we will take the signal decomposition
\begin{equation}
    \label{eq:relative:signflip:signal:svd}
    \bmS = \sum_{i = 1}^{k} \theta_i \bmu_i \bmz_i^\top
    = \bmU \operatorname{diag}(\bmtheta) \bmZ^\top
    \in \bbR^{n \times p}
    ,
\end{equation}
in the signal-plus-noise model \cref{eq:signal:plus:noise} to be a compact singular value decomposition, where
\begin{itemize}
    \item $\bmtheta = (\theta_1,\dots,\theta_k) \in \bbR^{k}$ is the vector of positive singular values $\theta_1 \geq \dots \geq \theta_k > 0$,
    \item $\bmU = [\bmu_1,\dots,\bmu_k] \in \bbR^{n \times k}$ is the matrix of orthonormal left singular vectors, and
    \item $\bmZ = [\bmz_1,\dots,\bmz_k] \in \bbR^{p \times k}$ is the matrix of orthonormal right singular vectors.
\end{itemize}
Moreover, we will suppose here that $\bmS$ is deterministic;
similar results hold for random $\bmS$.

Our first result quantifies the interplay of the signal strengths, noise preservation, and signal destruction.
As in the proof of \cref{thm:consistency} (in \cref{pf:thm:consistency}),
the crucial comparison is between the weaker signal components and the much smaller operator norm $\|\bmR\circ\bmS\|$ of the whole signal after signflipping.
The following \lcnamecref{thm:relative:selection} provides a scale-free version of this statement (\cref{thm:consistency} assumed the data was scaled, or re-scaled, to produce a noise operator norm of constant order).
Namely, here the signflipped signal is allowed to be non-vanishing, and it simply becomes part of the empirical null generated by FlipPA.
A component is selected once its observed singular value exceeds the noise upper-edge by the combination of the noise-preservation error, the algorithmic threshold, and the signflipped-signal size, which we call the signflipped-signal scale.
Here $O_{\Pr}(\cdot)$ denotes boundedness in probability, and all deterministic sequences may depend on both $n$ and $p$.

\begin{theorem}[Asymptotic consistency in terms of the signflipped-signal scale]
  \label{thm:relative:selection}
  Suppose the signal-plus-noise model~\cref{eq:signal:plus:noise} satisfies
  \begin{align}
    \label{eq:relative:condition}
    \|\bmR\circ\bmS\| &= O_{\Pr}(a_n),
    &
    \big|\|\bmR\circ\bmN\|-\|\bmN\|\big| &= O_{\Pr}(b_n),
  \end{align}
  for some deterministic positive sequences $a_n$ and $b_n$,
  where $\bmR \sim \operatorname{Unif}(\{-1, 1\}^{n \times p})$ is a Rademacher random matrix (independent of $\bmS$ and $\bmN$).
  Suppose further that
  \begin{equation}
    \label{eq:relative:observed:gap}
    \Pr[\sigma_k(\bmX)>\|\bmN\|+\varepsilon_n] \to 1
    ,
  \end{equation}
  for some deterministic positive sequence $\varepsilon_n > 0$.
  Then FlipPA using the upper-edge comparison method
  with threshold $\tau = \tau_n$
  is consistent,
  i.e.,
  $\Pr\big[ \htk = k \big] \to 1$,
  as long as $\tau_n < \varepsilon_n$ eventually and $a_n+b_n=o\{\min(\tau_n,\varepsilon_n-\tau_n)\}$.
\end{theorem}

\Cref{thm:relative:selection} establishes consistency%
\footnote{For simplicity, we focus here on consistent, i.e., exact, selection $\htk = k$. Selection of all perceptible components can be established in a similar fashion under the weaker condition that $a_n+b_n=o\{\varepsilon_n-\tau_n\}$.}
of FlipPA as long as the observed singular value $\sigma_k(\bmX)$ separates from the noise edge $\|\bmN\|$ by more than the sum of the signflipped-signal scale $a_n$, the noise-preservation error $b_n$, and the algorithmic threshold $\tau_n$.
The settings studied in \cref{sec:consistency} correspond to taking $a_n+b_n=o(1)$ with a fixed positive gap $\varepsilon_n=\varepsilon$ and a fixed threshold $\tau\in(0,\varepsilon)$.
The proof is essentially a refined variant of \cref{pf:thm:consistency}.

\begin{proof}[Proof of \cref{thm:relative:selection}]
    As in \cref{pf:thm:consistency}, we have that FlipPA eventually correctly selects $\htk = k$ whenever the following conditions are simultaneously met:
    \begin{gather*}
        \sigma_k > \| \bmN \| + \varepsilon_n
        \quad \text{and} \quad
        \sigma_{k+1} \leq \| \bmN \|
        \\
        \text{and} \quad
        \| \bmN \| - \tau_n \leq \tlsigma_1^{(t)} \leq \| \bmN \| + (\varepsilon_n-\tau_n)
        \text{ for $t = 1,\dots,T$}
        ,
    \end{gather*}
    since $\tau_n < \varepsilon_n$ eventually,
    which leads via a union bound to the simple bound
    \begin{align}
        \label{eq:relative:selection:unionbound}
        \Pr[\htk = k]
        &
        \geq
        1
        +
        \left(
            \Pr[\sigma_k > \| \bmN \| + \varepsilon_n] - 1
        \right)
        +
        \left(
            \Pr[\sigma_{k+1} \leq \| \bmN \|] - 1
        \right)
        \\&\quad \nonumber
        +
        \sum_{t=1}^{T} \left(
            \Pr[\| \bmN \| - \tau_n \leq \tlsigma_1^{(t)} \leq \| \bmN \| + (\varepsilon_n-\tau_n)] - 1
        \right)
        .
    \end{align}
    Similar to \cref{pf:thm:consistency},
    it follows immediately from \cref{eq:relative:observed:gap}
    and from Weyl's inequality that
    \begin{align}
        \label{eq:relative:selection:noisebounds}
        \Pr[ \sigma_k > \| \bmN \| + \varepsilon_n ] &\to 1
        , &
        \Pr[ \sigma_{k+1} \leq \| \bmN \| ] &= 1
        ,
    \end{align}
    so it remains to show that
    $\Pr[\| \bmN \| - \tau_n \leq \tlsigma_1^{(t)} \leq \| \bmN \| + (\varepsilon_n-\tau_n)] \to 1$
    for $t = 1,\dots,T$.
    Considering the same decomposition as in \cref{eq:consistency:split} yields
    \begin{equation*}
        \big| \tlsigma_1^{(t)} - \| \bmN \| \big|
        \leq
        \underbrace{
            \big| \sigma_1(\bmR^{(t)} \circ \bmX) - \sigma_1(\bmR^{(t)} \circ \bmN) \big|
        }_{\text{signal destruction}}
        +
        \underbrace{
            \big| \sigma_1(\bmR^{(t)} \circ \bmN) - \sigma_1(\bmN) \big|
        }_{\text{noise preservation}}
        ,
    \end{equation*}
    where we now have from \cref{eq:relative:condition} that these terms are bounded as
    \begin{align*}
        \big| \sigma_1(\bmR^{(t)} \circ \bmX) - \sigma_1(\bmR^{(t)} \circ \bmN) \big|
        &
        \leq
        \| \bmR^{(t)} \circ \bmS \|
        =
        O_{\Pr}(a_n)
        , \\
        \big| \sigma_1(\bmR^{(t)} \circ \bmN) - \sigma_1(\bmN) \big|
        &
        =
        \big| \|\bmR^{(t)} \circ \bmN\| - \|\bmN\| \big|
        =
        O_{\Pr}(b_n)
        ,
    \end{align*}
    since $\bmR^{(t)} =_d \bmR$.
    Thus, it follows that $\big| \tlsigma_1^{(t)} - \| \bmN \| \big| = O_{\Pr}(a_n + b_n)$
    and consequently
    \begin{equation}
        \label{eq:relative:selection:signfliplim}
        \Pr[\| \bmN \| - \tau_n \leq \tlsigma_1^{(t)} \leq \| \bmN \| + (\varepsilon_n-\tau_n)]
        \geq
        \Pr\Big[\big| \tlsigma_1^{(t)} - \| \bmN \| \big| \leq \min(\tau_n,\varepsilon_n-\tau_n)\Big]
        \to 1
        ,
    \end{equation}
    as long as $a_n+b_n=o\{\min(\tau_n,\varepsilon_n-\tau_n)\}$.
    Combining \cref{eq:relative:selection:unionbound,eq:relative:selection:noisebounds,eq:relative:selection:signfliplim} completes the proof.
\end{proof}

Our next result quantifies the signflipped-signal scale $a_n$ for dense bounded-rank signals
in terms of the strongest signal singular value $\theta_1$ and the coherence of the singular vectors, given by
\begin{equation}
    \label{eq:relative:coherence}
    \mu_U \coloneqq \frac{n}{k}\|\bmU\|_{2,\infty}^2,
    \qquad
    \mu_Z \coloneqq \frac{p}{k}\|\bmZ\|_{2,\infty}^2
    .
\end{equation}
It shows that $a_n$ is controlled by $\theta_1$ divided by the square root of the minimum dimension $m=\min(n,p)$ when the components are incoherent.

\begin{theorem}[Signflipped-signal scale for incoherent components]
    \label{thm:relative:signflip:scale}
    Suppose the signal matrix $\bmS \in \bbR^{n \times p}$ has the compact SVD given in \cref{eq:relative:signflip:signal:svd},
    and let $\bmR \sim \operatorname{Unif}(\{-1, 1\}^{n \times p})$ be a Rademacher random matrix.
    Then, we have
    \begin{equation}
        \label{eq:relative:signflip:maxscale}
        \|\bmR \circ \bmS\|
        = O_{\Pr}\left\{
            \theta_1
            \left(
                \sqrt{\frac{k\mu_U}{n}}
                +
                \sqrt{\frac{k\mu_Z}{p}}
                +
                \frac{k(\mu_U\mu_Z)^{1/4}}{(np)^{1/4}}
            \right)
        \right\}
        .
    \end{equation}
    In particular,
    when the signal has bounded rank and incoherent components,
    i.e., when $k = O(1)$, $\mu_U = O(1)$, and $\mu_Z = O(1)$,
    then we have
    \begin{equation}
        \label{eq:relative:signflip:simple}
        \|\bmR \circ \bmS\|
        = O_{\Pr} (\theta_1 / \sqrt{m})
        ,
    \end{equation}
    where $m = \min(n,p)$.
\end{theorem}

\begin{proof}[Proof of \cref{thm:relative:signflip:scale}]
    Combining \cref{lem:opnorm:bound} with the bound \cref{eq:rhoinf:bound} from \cref{thm:decay:suff:proof} and recalling that $\rank(\bmS) = k$ here yields
    \begin{align}
        \label{eq:relative:signflip:scale:initbound}
        \bbE \| \bmR \circ \bmS \|
        \lesssim
        \|\bmS\|_{2,\infty}
        +
        \|\bmS^\top\|_{2,\infty}
        +
        \sqrt{k} \sqrt{\|\bmS\|_{2,\infty} \cdot \|\bmS^\top\|_{2,\infty}}
        .
    \end{align}
    Note next that for any matrices $\bmA$ and $\bmB$
    \begin{equation*}
        \|\bmA \bmB\|_{2,\infty}
        =
        \max_{i=1,\dots,n} \|(\bmA \bmB)_{i:}\|_2
        \leq
        \max_{i=1,\dots,n} \|\bmA_{i:}\|_2 \|\bmB\|
        =
        \|\bmA\|_{2,\infty} \|\bmB\|
        ,
    \end{equation*}
    so we have by appropriately grouping the terms in the SVD $\bmS = \bmU \operatorname{diag}(\bmtheta) \bmZ^\top$ that
    \begin{align}
        \label{eq:relative:signflip:scale:2infbound:S}
        \|\bmS\|_{2,\infty}
        &
        \leq
        \|\bmU\|_{2,\infty} \|\operatorname{diag}(\bmtheta) \bmZ^\top\|
        =
        \theta_1 \|\bmU\|_{2,\infty}
        =
        \theta_1 \sqrt{\frac{k\mu_U}{n}}
        , \\
        \label{eq:relative:signflip:scale:2infbound:St}
        \|\bmS^\top\|_{2,\infty}
        &
        \leq
        \|\bmZ\|_{2,\infty} \|\operatorname{diag}(\bmtheta) \bmU^\top\|
        =
        \theta_1 \|\bmZ\|_{2,\infty}
        =
        \theta_1 \sqrt{\frac{k\mu_Z}{p}}
        .
    \end{align}
    Combining \cref{eq:relative:signflip:scale:initbound,eq:relative:signflip:scale:2infbound:S,eq:relative:signflip:scale:2infbound:St}, using Markov's inequality, and simplifying completes the proof.
\end{proof}

Combining \cref{thm:relative:selection,thm:relative:signflip:scale} yields the following consistency guarantee for FlipPA in terms of the relative strength of the signals that applies in the presence of strong signals.

\begin{theorem}[Asymptotic consistency of FlipPA under strong signals]
    \label{cor:relative:dense:factors}
    Consider the signal-plus-noise model \cref{eq:signal:plus:noise}, where
    \begin{itemize}
        \item the signal matrix $\bmS$ has the compact SVD given in \cref{eq:relative:signflip:signal:svd} with bounded rank and incoherent components, i.e., $k = O(1)$, $\mu_U = O(1)$, and $\mu_Z = O(1)$, and
        \item the noise matrix $\bmN$ is preserved by signflipping at a rate of $\big|\|\bmR\circ\bmN\|-\|\bmN\|\big| = O_{\Pr}(b_n)$ for some deterministic positive sequence $b_n$.
    \end{itemize}
    Suppose further that the $k$-th component is perceptible in the sense that there exists some constant $c > 0$ for which
    \begin{equation}
        \label{eq:relative:dense:factors:eps}
        \Pr[\sigma_k(\bmX)>\|\bmN\|+c\theta_k] \to 1
        ,
    \end{equation}
    and the threshold $\tau_n$ is chosen to satisfy
    \begin{align}
        \label{eq:relative:dense:factors:tau}
        b_n + \tau_n &= O(\theta_1/\sqrt{m})
        , &
        b_n + \theta_1/\sqrt{m} &= o(\tau_n)
        .
    \end{align}
    Then FlipPA using the upper-edge comparison method
    with threshold $\tau = \tau_n$
    is consistent,
    i.e.,
    $\Pr\big[ \htk = k \big] \to 1$,
    as long as $\theta_k = \omega(\theta_1/\sqrt{m})$.
\end{theorem}

\begin{proof}[Proof of \cref{cor:relative:dense:factors}]
    It follows from the assumptions and \cref{thm:relative:signflip:scale} that the signal and noise matrices satisfy the condition \cref{eq:relative:condition} of \cref{thm:relative:selection} with $a_n = \theta_1/\sqrt{m}$ and $b_n$ as given.
    Moreover, it follows from \cref{eq:relative:dense:factors:eps} that the $k$-th component satisfies \cref{eq:relative:observed:gap} with $\varepsilon_n = c\theta_k$.

    It now remains to verify that $\tau_n < \varepsilon_n$ eventually and that $a_n+b_n=o\{\min(\tau_n,\varepsilon_n-\tau_n)\}$.
    First, note that it follows from \cref{eq:relative:dense:factors:tau} that $\tau_n \leq b_n + \tau_n \lesssim \theta_1/\sqrt{m} = o(\theta_k)$.
    Since $\varepsilon_n = c\theta_k$, it follows that $\tau_n = o(\varepsilon_n)$ and so $\tau_n < \varepsilon_n$ eventually.
    Note next that by \cref{eq:relative:dense:factors:tau} we have that $a_n + b_n = \theta_1/\sqrt{m} + b_n = o(\tau_n)$.
    Moreover, $b_n \leq b_n + \tau_n = O(\theta_1/\sqrt{m})$, and so
    \begin{equation*}
        \frac{a_n + b_n}{\varepsilon_n - \tau_n}
        =
        \frac{\theta_1/\sqrt{m} + b_n}{c\theta_k - \tau_n}
        \lesssim
        \frac{\theta_1/\sqrt{m}}{c\theta_k - \tau_n}
        =
        \frac{(\theta_1/\sqrt{m})/\theta_k}{c - \tau_n/\theta_k}
        \to
        \frac{0}{c - 0}
        = 0
        .
    \end{equation*}
    Thus, we have $a_n+b_n=o\{\min(\tau_n,\varepsilon_n-\tau_n)\}$,
    which completes the proof.
\end{proof}

\Cref{cor:relative:dense:factors} establishes the consistency of FlipPA for bounded-rank incoherent components under the condition that $\theta_k = \omega(\theta_1/\sqrt{m})$, i.e., as long as $\theta_k \sqrt{m}/\theta \to \infty$.
In other words, a component is selected once its strength is asymptotically larger than $\theta_1/\sqrt{m}$, instead of the full strongest-signal scale of $\theta_1$.
Notably, this condition permits factor models with strong signals.
In particular, if $\bmX = \bmF\bmL^\top + \bmN$, then in the strong-factor setting of \citet{bai2002determining}, one has $\theta_1 \asymp \cdots \asymp \theta_k \asymp \sqrt{np}$ for each $j$.
Recent work has also considered weak loadings in factor models, for which the signal scale can satisfy $\|\bmS\| \ll \sqrt{np}$ \citep{bai2023approximate,choi2025high}.
The signflipped remnant of the strongest factor has size only $O_{\Pr}\{\theta_1/\sqrt m\}$ for dense factors, which is of the same order as the ordinary Gaussian noise edge when $n\asymp p$, and is much smaller than the full signal scale $\sqrt{np}$.
Hence, FlipPA does not require $\theta_1\ll\sqrt{np}$ to select strong factors.
In a mixed strong/weak model, stronger factors can shadow weaker factors only through the reduced scale $\theta_1/\sqrt m$.
Any perceptible component satisfying $\theta_j/\theta_1\gg1/\sqrt m$ is selected with probability tending to one; when $n\asymp p$, this is the claimed scale $\theta_j/\theta_1\gg1/\sqrt n$.
}

\section{Proof of \cref{thm:perm:homvar}}
\label{thm:perm:homvar:proof}

Note first that the homogenized noise matrix $\bbrN$
can be written as $\bbrN = \bmG \diag^{1/2}(\bmv)$,
where $\bmG \in \bbR^{n \times p}$ has entries $G_{ij} \overset{iid}{\sim} \clN(0, 1)$
and the empirical distribution function of $(v_1,\dots,v_p)$ converges to $H$.
Thus,
as is well known,
the empirical singular value distribution of $\bbrN/\sqrt{n}$ converges
to the generalized Mar{\v{c}}enko-Pastur distribution defined by $H$
\citep[see, e.g.,][Theorem 4.3]{bai2009spectral}.
Namely,
$\bbrN \bbrN^\top / n = (1/n) \bmG \diag(\bmv) \bmG^\top$
has a limiting spectral distribution
defined by the unique Stieltjes transform $m(z)$
that satisfies the equation
\begin{equation}
    \label{eq:marchenkopastur:stieltjes}
    z + \frac{1}{m(z)}
    =
    \gamma \int \frac{t dH(t)}{1 + t m(z)}
    ,
    \qquad z \in \bbC^+
    .
\end{equation}
So, it remains to show that
the empirical singular value distribution of $\bmN_\pi/\sqrt{n}$ also converges
to the generalized Mar{\v{c}}enko-Pastur distribution defined by $H$.
Namely, we need to show that $\bmN_\pi \bmN_\pi^\top / n$
has a limiting spectral distribution
defined by the unique Stieltjes transform $m(z)$
that satisfies \cref{eq:marchenkopastur:stieltjes}.
Since the random permutations $\bmpi_1,\dots,\bmpi_p$
induce dependence among the entries of $\bmN_\pi$,
doing so requires establishing the generalized Mar{\v{c}}enko-Pastur law
without assuming that the entries are independent.

The remainder of the proof consists of two parts.
First, \cref{proof:noise:baizhou} establishes general conditions
under which random matrices with dependent entries
still follow the generalized Mar{\v{c}}enko-Pastur law.
The conditions include settings
beyond what is needed here
and may be of independent interest.
Then, \cref{proof:noise:hetero:perm} completes the proof by
showing that $\bmN_\pi^\top$ satisfies these conditions.

\subsection{General conditions for the generalized Mar{\v{c}}enko-Pastur law}
\label{proof:noise:baizhou}

This section proves the following \lcnamecref{lemma:noise:baizhou},
which establishes the generalized Mar{\v{c}}enko-Pastur law
under relaxed independence conditions.
Specifically,
it allows for some dependence within each row
in the case where the rows all have isotropic covariances.
See \citet{hui2010lsd,wei2016tls,bryson2019marchenko}
for some related (but different) results.

\begin{lemma}
  [Generalized Mar{\v{c}}enko-Pastur law under relaxed independence conditions]
  \label{lemma:noise:baizhou}
  Let $\bmX \in \bbR^{n \times p}$
  have independent rows with zero mean entries,
  and suppose that $n,p \to \infty$ with $p/n \to \gamma > 0$.
  Suppose furthermore that
  \begin{enumerate}
    \item Each row $\bmx_k \in \bbR^{p}$ of $\bmX$ has an isotropic covariance of
      $\E(\bmx_k \bmx_k^\top) = \eta_k^2 \bmI_p$.
    \item The variances $\eta_1^2,\dots,\eta_n^2$ are uniformly bounded
      with empirical distribution
      converging to some deterministic distribution $H$.
    \item For any sequence of possibly complex-valued symmetric deterministic $p \times p$ matrices $(\bmA_p)_{p\ge 1}$
      with uniformly bounded spectral norms,
      and for every row $\bmx_k$,
      we have
      \begin{equation} \label{eq:var:quadratic}
        \var\left(\bmx_k^\top \bmA_p\bmx_k\right)
        \le C \|\bmA_p\|_F^2,
      \end{equation}
    where $C$ does not depend on the sequence $(\bmA_p)_{p\ge 1}$,
    and $\var(Z) = \E |Z-\E Z|^2$ denotes the variance
    for a complex-valued random variable.
  \end{enumerate}
  Then, with probability one,
  the empirical spectral distribution
  of $n^{-1} \bmX^\top \bmX$ converges weakly
  to the generalized Mar{\v{c}}enko-Pastur distribution,
  whose Stieltjes transform $m(z)$ satisfies:
  \begin{equation} \label{eq:homogenized:perm:stieltjes:yue}
    z+\frac{1}{m(z)}=\int\frac{t}{1+\gamma t m(z)}dH(t),~~~z\in\mathbb C^+.
  \end{equation}
\end{lemma}

The next section (\cref{proof:noise:hetero:perm})
shows that $\bmN_\pi^\top$ satisfies these conditions.
The remainder of this section proves \cref{lemma:noise:baizhou}
by carefully combining techniques used in the proofs of \cite[Theorem 1.1]{bai2008independence} and \cite[Theorem 4.3]{bai2009spectral}.

\begin{proof}[Proof of \cref{lemma:noise:baizhou}]
Let $m_n(z) = p^{-1} \tr ( n^{-1} \bmX^\top \bmX - z\bmI_p )^{-1}$
be the Stieltjes transform for the empirical spectral distribution of $n^{-1} \bmX^\top \bmX$.
Following the proof of Theorem 1.1 of \cite{bai2008independence}, we proceed in three steps:
\begin{enumerate}
  \item $m_n(z) - \E m_n(z) \to 0$, a.s.
  \item $\E m_n(z) \to m(z)$, which satisfies \cref{eq:homogenized:perm:stieltjes:yue}.
  \item \cref{eq:homogenized:perm:stieltjes:yue} has a unique solution in $\bbC^+$.
\end{enumerate}
Note that $m_n(z) = p^{-1} \tr\clB_n^{-1}$ where
\begin{align*}
  \clB_n &\coloneqq \bmB_n - z \bmI_p
  \in \bbC^{p \times p}
  , &
  \bmB_n &\coloneqq \frac{1}{n} \sum_{i=1}^n \bmx_i \bmx_i^\top
  \in \bbC^{p \times p}
  .
\end{align*}
For any $k\in\{1,\dots,n\}$, we also define 
\begin{align*}
  \clB_{k,n} &\coloneqq \bmB_{k,n} - z \bmI_p
  \in \bbC^{p \times p}
  , &
  \bmB_{k,n} &\coloneqq \frac{1}{n} \sum_{i \neq k} \bmx_i \bmx_i^\top
  \in \bbC^{p \times p}
  .
\end{align*}
Throughout the proof, for fixed $z\in\mathbb C^+$, we will write $z=\Re(z)+i\Im(z)=u+iv$, where $u\in\bbR,v>0$ are the real and imaginary parts of $z$. 
Also, we consider a constant $L>0$ such that for all $n$, $\eta_k^2 \le L$, $k=1,\ldots,n$.

\textbf{Step 1.} $m_n(z) - \E m_n(z) \to_{a.s.} 0$.

Using the notation $\bbE_k(\cdot) = \bbE(\cdot | \bmx_{k+1},\dots,\bmx_n)$, we have
\begin{align*}
  m_n(z) - \E m_n(z)
  &= \E_0 m_n(z) - \E_n m_n(z)
  = \sum_{k=1}^n (\E_{k-1} m_n(z) - \E_k m_n(z))
  \\
  &= \frac{1}{p} \sum_{k=1}^n (\E_{k-1}-\E_k)(\tr\clB_n^{-1} - \tr\clB_{k,n}^{-1})
  = \frac{1}{p} \sum_{k=1}^n (\E_{k-1}-\E_k) \nu_k
\end{align*}
where $\nu_k \coloneqq \tr\clB_n^{-1} - \tr\clB_{k,n}^{-1}$. By using Lemma 2.6 of \cite{silverstein1995empirical}, we have $|\nu_k| \leq v^{-1}$. 
Thus, $(\E_{k-1}-\E_k) \nu_k$ forms a bounded martingale difference sequence
and applying the Burkholder inequality \cite[Lemma 2.12]{bai2009spectral} yields,
for any $q\ge 2$ and for some $K_q>0$ depending only on $q$,
\begin{equation*}
  \E|m_n(z) - \E m_n(z)|^q
  \leq K_q p^{-q}
  \E\bigg(
    \sum_{k=1}^n |(\E_{k-1}-\E_k) \nu_k|^2
  \bigg)^{q/2}
  \leq K_q \bigg(\frac{2}{v}\bigg)^q p^{-q/2} \bigg(\frac{p}{n}\bigg)^{-q/2}
  .
\end{equation*}
By choosing $q > 2$, this implies $m_n(z) - \E m_n(z) \to_{a.s.} 0$ due to the Borel-Cantelli lemma.

\textbf{Step 2.} $\E m_n(z) \to m(z)$, which satisfies \cref{eq:homogenized:perm:stieltjes:yue}.

We define the scalars $K,\tilde K$
  by
\begin{equation*}
K=\frac{1}{n}\sum_{k=1}^n\frac{\eta_k^2}{1+n^{-1}\tr\clB_{k,n}^{-1}\eta_k^2},
~~~
\tilde{K}=\frac{1}{n}\sum_{k=1}^n\frac{\eta_k^2}{1+n^{-1}\E\tr\clB_{n}^{-1}\eta_k^2}.
\end{equation*}
Since $v>0$, it directly follows that $\Im K, \Im\tilde{K}<0$, 
hence $|(K - z)^{-1}|\leq 1/v$ and $|(\tilde K - z)^{-1}|\leq 1/v$. 

Now, 
by using the resolvent identity $\bmA^{-1} - \bmB^{-1} = -\bmA^{-1}(\bmA-\bmB)\bmB^{-1}$, holding for any two invertible square matrices of the same size, we have
\begin{align*}
  &(K - z)^{-1}\bmI_p - (\bmB_n - z\bmI_p)^{-1}
  \\&\qquad
  =
  \bigg\{
    \frac{1}{n}\sum_{k=1}^n (K - z)^{-1} \bmx_k \bmx_k^\top (\bmB_n - z\bmI_p)^{-1}
  \bigg\}
  - (K - z)^{-1}K (\bmB_n - z\bmI_p)^{-1}
  \\&\qquad
  = 
  \bigg\{
    \sum_{k=1}^n \frac{(K - z)^{-1}n^{-1} \bmx_k \bmx_k^\top \clB_{k,n}^{-1}}{1 + n^{-1}\bmx_k^\top\clB_{k,n}^{-1}\bmx_k}
  \bigg\}
  - (K - z)^{-1}K (\bmB_n - z\bmI_p)^{-1}
\end{align*}
where the last line uses that
\begin{align*}
  \bmx_k^\top \clB_n^{-1}
  &= \bmx_k^\top \clB_{k,n}^{-1}
  -
  \frac{\bmx_k^\top \clB_{k,n}^{-1} (n^{-1}\bmx_k\bmx_k^\top) \clB_{k,n}^{-1}}
    {1 + n^{-1}\bmx_k^\top\clB_{k,n}^{-1}\bmx_k}
  = 
  \frac{\bmx_k^\top \clB_{k,n}^{-1}}{1 + n^{-1}\bmx_k^\top\clB_{k,n}^{-1}\bmx_k}
  .
\end{align*}
Taking the trace and dividing by $p$ yields
\begin{align}
  &
  (K - z)^{-1} - \frac{1}{p} \tr \clB_n^{-1}
  = \frac{1}{p}
  \bigg\{
    \sum_{k=1}^n
      \frac{n^{-1} (K - z)^{-1} \bmx_k^\top \clB_{k,n}^{-1}  \bmx_k}
        {1 + n^{-1}\bmx_k^\top\clB_{k,n}^{-1}\bmx_k}
  \bigg\}
  - \frac{1}{p} (K - z)^{-1} K\tr [ \clB_n^{-1}]
  \nonumber\\
  &=\frac{1}{p}
  \bigg\{
    \sum_{k=1}^n
      \frac{n^{-1} (K - z)^{-1} \bmx_k^\top \clB_{k,n}^{-1}  \bmx_k}
        {1 + n^{-1}\bmx_k^\top\clB_{k,n}^{-1}\bmx_k}
  \bigg\}
  - \frac{1}{p} \bigg\{\sum_{k=1}^n\frac{n^{-1}\eta_k^2(K - z)^{-1} \tr [ \clB_n^{-1}]}{1+n^{-1}\tr\clB_{k,n}^{-1}\eta_k^2}\bigg\}
  \nonumber\\
  &= \frac{1}{p} \sum_{k=1}^n \frac{d_k}{1 + n^{-1}\bmx_k^\top\clB_{k,n}^{-1}\bmx_k}\label{kbdiff}
  ,
\end{align}
where for $k\in\{1,\dots,n\}$,
\begin{align*}
  &d_k\coloneqq n^{-1} (K - z)^{-1} \bmx_k^\top \clB_{k,n}^{-1}  \bmx_k- n^{-1}\eta_k^2 (K - z)^{-1} \tr [ \clB_n^{-1}] \left(\frac{1 + n^{-1}\bmx_k^\top\clB_{k,n}^{-1}\bmx_k}{1+n^{-1}\tr\clB_{k,n}^{-1}\eta_k^2}\right)\\
  &= d_{k1} + d_{k2} + d_{k3}
  ,
\end{align*}
with
\begin{align*}
  d_{k1} &\coloneqq
  n^{-1}\eta_k^2 (K - z)^{-1} \tr [ \clB_{k,n}^{-1}]- n^{-1} \eta_k^2(K - z)^{-1} \tr [ \clB_n^{-1}]
  , \\
  d_{k2} &\coloneqq
  n^{-1} (K - z)^{-1} \bmx_k^\top \clB_{k,n}^{-1}  \bmx_k- n^{-1}\eta_k^2 (K - z)^{-1} \tr [\clB_{k,n}^{-1}]
  , \\
  d_{k3} &\coloneqq
  n^{-1} (K - z)^{-1} \tr[\clB_n^{-1}]
  \cdot
  \eta_k^2 
\bigg(1 - \frac{1+n^{-1}\bmx_k^\top \clB_{k,n}^{-1}\bmx_k}{1+ n^{-1}\eta_k^2\tr\clB_{k,n}^{-1}}\bigg)
  .
\end{align*}
By using Lemma 2.6 of \cite{silverstein1995empirical} and the fact that $|(K - z)^{-1}|\leq 1/v$, we have 
\begin{align}\label{dk1}
|d_{k1}|\leq\frac{\eta_k^2|(K - z)^{-1}|}{nv}\leq\frac{L}{nv^2}.
\end{align}
Further, we have
$
|d_{k2}|\leq\frac{1}{nv}\left|\bmx_k^\top \clB_{k,n}^{-1} \bmx_k- \eta_k^2\tr \clB_{k,n}^{-1}\right|$.
Thus, from $\E(\bmx_k\bmx_k^\top)=\eta_k^2\bmI_p$, \cref{eq:var:quadratic}
by using that all singular values of $\clB_{k,n}^{-1}$ are bounded by $1/v$,
\begin{align}\label{dk2}
\E|d_{k2}|^2
&
\leq\frac{1}{n^2v^2}\E\left|\bmx_k^\top \clB_{k,n}^{-1} \bmx_k- \eta_k^2\tr \clB_{k,n}^{-1}\right|^2
\\& \nonumber
= 
\frac{1}{n^2v^2}\var(\bmx_k^\top \clB_{k,n}^{-1} \bmx_k)
= 
\frac{O(\|\clB_{k,n}^{-1}\|_F^2)}{n^2v^2} 
= 
\frac{O(p)}{n^2v^4}
= 
\frac{O(1)}{nv^4}.
\end{align}

Now recall from \cite[Corollary 3.1]{couillet2011random}, that if $m$ is the Stieltjes transform of a measure on $\mathbb R$, then for any $z\in\mathbb C^+$, we have
\begin{equation}\label{stb}
\left|\frac{1}{1+m(z)}\right|\leq\frac{|z|}{\Im(z)}.
\end{equation}

For $d_{k3}$, by using 
$|(K - z)^{-1}|\leq 1/v$, 
$|\tr\clB_n^{-1}|\leq p/v$, 
$\eta_k^2 \le L$,
and that \eqref{stb} holds for $n^{-1}\eta_k^2\tr\clB_{k,n}^{-1}$ as it is a Stieltjes transform, we have
\begin{align*}
|d_{k3}|\leq
\frac{p\eta_k^2}{n^2v^2}\left|\frac{\bmx_k^\top\clB_{k,n}^{-1}\bmx_k-\eta_k^2\tr\clB_{k,n}^{-1}}{1+n^{-1}\eta_k^2\tr\clB_{k,n}^{-1}}\right|
\leq\frac{p|z|L}{n^2v^3}\left|\bmx_k^\top \clB_{k,n}^{-1} \bmx_k- \eta_k^2\tr \clB_{k,n}^{-1}\right|.
\end{align*}
Thus,
\begin{align}\label{dk3}
\E|d_{k3}|^2
\leq\frac{p^2|z|^2L^2}{n^4v^6}\E\left|\bmx_k^\top \clB_{k,n}^{-1} \bmx_k- \eta_k^2\tr \clB_{k,n}^{-1}\right|^2
= 
\frac{O(|z|^2L^2)}{nv^8}.
\end{align}
From \eqref{dk1}, \eqref{dk2}, and \eqref{dk3}, 
since $z,L$ are fixed,
we have, uniformly over $k\in\{1,\dots,n\}$,
\begin{align*}
|\E d_k|^2=|\E d_{k1}+\E d_{k2}+\E d_{k3}|^2
\leq 
3(\E|d_{k1}|^2+\E|d_{k2}|^2+\E|d_{k3}|^2) = O(1/n).
\end{align*}
Hence, from \eqref{kbdiff}, by using that \eqref{stb} holds for $n^{-1}\bmx_k^\top\clB_{k,n}^{-1}\bmx_k$, as it is a Stieltjes transform, 
and since the bound for $d_k$ is uniform over $k$,
\begin{align}\label{kbn}
\left|\E\left((K - z)^{-1} -  \frac{1}{p}\tr \clB_n^{-1}\right)\right|
&
=\left|\frac{1}{p}\sum_{k=1}^n \E \frac{d_k}{1 + n^{-1}\bmx_k^\top\clB_{k,n}^{-1}\bmx_k}\right|
\\& \nonumber
\leq\frac{|z|}{pv}\sum_{k=1}^n\left|\E d_k\right|
= 
\frac{|z|}{pv}\cdot O(n \cdot 1/\sqrt{n})
\to0.
\end{align}
Next, 
we have
\begin{align*}
\left|\E[(\tilde K - z)^{-1}- (K - z)^{-1}]\right|
&
\leq\frac{\E|\tilde K-K|}{v^2}
\\&
=\frac{1}{n^2v^2}\E\left|\sum_{k=1}^n\eta_k^4\frac{\E\tr\clB_n^{-1}-\tr\clB_{k,n}^{-1}}{(1+n^{-1}\eta_k^2\tr\clB_{k,n}^{-1})(1+n^{-1}\eta_k^2\E\tr\clB_n^{-1})}\right|.
\end{align*}
Since \eqref{stb} holds for $n^{-1}\eta_k^2\tr\clB_{k,n}^{-1}$ and $n^{-1}\eta_k^2\tr\clB_{k,n}^{-1}$ as they are Stieltjes transforms, 
from $\left|\tr\clB_{k,n}^{-1}-\tr\clB_n^{-1}\right|\leq 1/v$, 
and from $\E|\tr\clB_n^{-1}-\E\tr\clB_n^{-1}|\to 0$ from Step 1,
we have
\begin{align*}
&
\frac{|z|^2L^4}{n^2v^4}\sum_{k=1}^n\E\left|\tr\clB_{k,n}^{-1}-\E\tr\clB_n^{-1}\right|
\\&\qquad
\leq\frac{|z|^2L^4}{n^2v^4}\sum_{k=1}^n\left(\E\left|\tr\clB_n^{-1}-\E\tr\clB_n^{-1}\right|+\E\left|\tr\clB_{k,n}^{-1}-\tr\clB_n^{-1}\right|\right)
\to0.
\end{align*}
Thus, from \eqref{kbn}, we reach
\begin{align}
\label{eq:fixpoint}
\E\left((\tilde K - z)^{-1} - \frac{1}{p} \tr \clB_n^{-1}\right)=\left(\frac{1}{n}\sum_{k=1}^n\frac{\eta_k^2}{1+n^{-1}\eta_k^2\E\tr\clB_n^{-1}}-z\right)^{-1}-\frac{1}{p}\E\tr\clB_n^{-1}\to0.
\end{align}
For each fixed $z\in\mathbb C^+$, $\E m_n(z)=p^{-1}\E\tr\clB_n^{-1}$ is a bounded sequence. 
Thus, for any subsequence $\{n'(a)\}_{a\ge 1}$ of the values taken by $n$, 
there is a subsubsequence $\{n''(b)\}_{b\ge 1}$ of  $\{n'(a)\}_{a\ge 1}$, such that $\E m_{n''(b)}(z)$ converges to a limit $m(z)$ as $b\to\infty$. 

Now, 
for all $t\in [0,L]$ and $n$ 
$$
\left|\frac{t}{1+n^{-1}t\E\tr\clB_n^{-1}}- \frac{t}{1+\gamma tm}\right|
= 
\left|\frac{\gamma t^2 (m- \E m_n(z))}
{(1+n^{-1}t m_n(z))(1+\gamma tm)}\right|
\le
\frac{\gamma L^2  |z|^2 |m- \E m_n(z)|}{v^2}.
$$
and thus the functions $f_n(t) = t/(1+\gamma tm_n(z))$ converge to $f(t) = t/(1+\gamma tm(z))$ uniformly over
$t\in [0,L]$, and over the sequence $\{n''(b)\}_{b\ge 1}$.
Then, from \cref{eq:fixpoint}, 
by using that the empirical distribution of $(\eta_k^2)_{k=1,\ldots,n}$ converges to the distribution $H$ on $[0,L]$ (on any sequence, in particular on $\{n''(b)\}_{b\ge 1}$), we have, on the sequence $\{n''(b)\}_{b\ge 1}$
$$
\frac{1}{n}\sum_{k=1}^n\frac{\eta_k^2}{1+n^{-1}\eta_k^2\E\tr\clB_n^{-1}}
\to 
\int\frac{t}{1+\gamma tm}dH(t).
$$
Further, we have $\Im(\int\frac{t}{1+\gamma tm}dH(t))\le 0$ and $\Im(z)>0$.
Hence
$m$ 
satisfies the equation
\begin{equation}
\label{eq:final:fixpoint}
\left(\int\frac{t}{1+\gamma tm}dH(t)-z\right)^{-1}=m.
\end{equation}
We will show in the next step that the solution to \cref{eq:final:fixpoint} is unique.
Thus, $\E m_n(z)$ converges to a limit which is the unique solution to \cref{eq:final:fixpoint}. Combining Step 1, we have $m_n(z)\to_{a.s.} m(z)$ for any fixed $z\in\mathbb C^+$. Finally, applying a standard argument based on Vitali's convergence theorem (e.g., see the proof of Theorem 2.9 of \cite{bai2009spectral}) yields $m_n(z)\to_{a.s.} m(z)$ for all $z\in\mathbb C^+$, 
where for all $z\in\mathbb C^+$, $m(z)$ is the unique solution to \cref{eq:final:fixpoint}.

\textbf{Step 3.} Show \cref{eq:homogenized:perm:stieltjes:yue} has a unique solution in $\mathbb C^+$. This step follows immediately as a special case of \cite[Proof of Theorem 4.3, Step 3]{bai2009spectral}.
\end{proof}

\subsection{Showing that \texorpdfstring{$\bmN_\pi^\top$}{N\_pi\^{}T} satisfies the conditions of \cref{lemma:noise:baizhou}}
\label{proof:noise:hetero:perm}

This section completes the proof of \cref{thm:perm:homvar}
by applying \cref{lemma:noise:baizhou}.
Since $\bmN_\pi$ has independent columns,
we apply the \lcnamecref{lemma:noise:baizhou} to $\bmN_\pi^\top$,
which has independent rows.
For this, we must show that $\bmN_\pi^\top$ satisfies
the three conditions of \cref{lemma:noise:baizhou}.
The first two conditions are straightforward to verify in this case:
\begin{enumerate}
    \item The covariance matrix of
    the $j$-th row of $\bmN_\pi^\top$ (i.e., the $j$-th column of $\bmN_\pi$)
    is
    \begin{equation*}
        \bbE \left[ (\bmN_\pi)_{:j}(\bmN_\pi)_{:j}^\top \right]
        =
        \bbE \left[ (\bmpi_j \bmN_{:j})(\bmpi_j \bmN_{:j})^\top \right]
        =
        \bbE \left[ \bmpi_j \bmN_{:j} \bmN_{:j}^\top \bmpi_j^\top \right]
        .
    \end{equation*}
    Since $\bmpi_j$ and $\bmN$ are independent
    and
    $\bbE \left[ \bmN_{:j} \bmN_{:j}^\top \right] = 
    \diag\left(
        \bbE \left[ |N_{1j}|^2 \right],
        \dots,
        \bbE \left[ |N_{nj}|^2 \right]
    \right)$,
    \begin{align*}
        &
        \bbE \left[ \bmpi_j \bmN_{:j} \bmN_{:j}^\top \bmpi_j^\top \right]
        =
        \bbE \left[
            \bmpi_j \,
            \bbE \left[ \bmN_{:j} \bmN_{:j}^\top \right]
            \bmpi_j^\top
        \right]
        \\& \qquad
        =
        \bbE \left[
            \bmpi_j
            \diag\left(
                \bbE \left( |N_{1j}|^2 \right),
                \dots,
                \bbE \left( |N_{nj}|^2 \right)
            \right)
            \bmpi_j^\top
        \right]
        \\& \qquad
        =
        \bbE \left[
            \sum_{m=1}^n
            \bbE \left( |N_{mj}|^2 \right) \cdot (\bmpi_j)_{:m}(\bmpi_j)_{:m}^\top
        \right]
        =
        \sum_{m=1}^n
        \bbE \left( |N_{mj}|^2 \right)
        \cdot
        \bbE \left[ (\bmpi_j)_{:m}(\bmpi_j)_{:m}^\top \right]
        .
    \end{align*}
    Finally,
    $\bbE \left[ (\bmpi_j)_{:m}(\bmpi_j)_{:m}^\top \right] = (1/n)\bmI_n$
    since $\bmpi_j$ is uniform over all permutations,
    so we have
    \begin{align*}
        \bbE \left[ (\bmN_\pi)_{:j}(\bmN_\pi)_{:j}^\top \right]
        &
        =
        \sum_{m=1}^n
        \bbE \left( |N_{mj}|^2 \right)
        \cdot
        \bbE \left[ (\bmpi_j)_{:m}(\bmpi_j)_{:m}^\top \right]
        \\&
        =
        \sum_{m=1}^n \left[
            \bbE \left( |N_{mj}|^2 \right)
            \cdot
            \frac{1}{n} \bmI_n
        \right]
        =
        \left[ \frac{1}{n} \sum_{m=1}^n \bbE |N_{mj}|^2 \right] \bmI_n
        =
        v_j \bmI_n
        .
    \end{align*}

    \item The row variances $v_1,\dots,v_p$
    have an empirical distribution converging to $H$ by assumption.
    Moreover, since the entries of $\bmN$ have uniformly bounded fourth moments,
    it follows that
    \begin{equation*}
        \max_{j=1,\dots,p} v_j
        =
        \max_{j=1,\dots,p} \left[
            \frac{1}{n} \sum_{m=1}^n \bbE |N_{mj}|^2
        \right]
        \leq
        \max_{\substack{j=1,\dots,p\\m=1,\dots,n}} \bbE |N_{mj}|^2
        \leq
        \max_{\substack{j=1,\dots,p\\m=1,\dots,n}} \bbE |N_{mj}|^4
        ,
    \end{equation*}
    is also uniformly bounded,
    i.e., $v_1,\dots,v_p$ are uniformly bounded.
\end{enumerate}

\noindent
It now remains to establish the third condition,
i.e., to bound the variance of the quadratic form \cref{eq:var:quadratic}.
This condition controls the amount of dependence
and is more involved.

\begin{enumerate}[start=3]
\item
For convenience,
let $\bmx_k = (\bmN_\pi)_{:k} \in \bbR^n$ denote the $k$-th row of $\bmN_\pi^\top$,
let $T_{ik} = \sqrt{\bbE |N_{ik}|^2}$ denote the standard deviation of $N_{ik}$,
and let $\eta_k = \sqrt{v_k}$ denote the standard deviation for the $k$-th row of $\bmN_\pi^\top$.
Now, note that
\begin{align*}
\var(\bmx_k^\top \bmA\bmx_k)
&
= 
\bbE |\bmx_k^\top \bmA\bmx_k -  \bbE \bmx_k^\top \bmA\bmx_k|^2
=
\bbE (\bmx_k^\top \bmA\bmx_k -  \bbE \bmx_k^\top \bmA\bmx_k) (\bmx_k^\top \bmA\bmx_k -  \bbE \bmx_k^\top \bmA\bmx_k)^*
\\&
=
\bbE [
\bmx_k^\top \bmA\bmx_k \bmx_k^\top \bmA^*\bmx_k] 
- \bbE \bmx_k^\top \bmA\bmx_k \cdot\bbE \bmx_k^\top \bmA^*\bmx_k.
\end{align*} 
where these expectations are all with respect to both
the randomness in the noise $\bmN$
and the random permutations $\bmpi$.
By using the independence of the coordinates of $\bmx_k$, as well as the symmetry of $\bmA$, we find
\begin{align*}
&\bbE [\bmx_k^\top \bmA\bmx_k \bmx_k^\top \bmA^*\bmx_k] 
=\sum_{i,j,l,m}A_{ij}A_{lm}^*\bbE (x_{ki}x_{kj}x_{kl}x_{km})\\
&=\left(\sum_{i=j=l=m}+\sum_{i=j\neq l=m}+\sum_{i=l\neq j=m}+\sum_{i=m\neq j=l}\right)A_{ij}A_{lm}^*\bbE (x_{ki}x_{kj}x_{kl}x_{km})\\
&=\sum_i|A_{ii}|^2\bbE (x_{ki}^4)
+\sum_{i\neq l}
(A_{ii}A_{ll}^* +
2|A_{il}|^2)\bbE (x_{ki}^2x_{kl}^2).
\end{align*}
For $i\neq l$, we have
\begin{equation*}
\bbE (x_{ki}^2x_{kl}^2)=\frac{\sum_{i\neq j}T_{ik}^2T_{lk}^2}{n(n-1)}
=\frac{(\sum_iT_{ik}^2)^2-\sum_iT_{ik}^4}{n(n-1)}
=\frac{n^2\eta_k^4-\sum_iT_{ik}^4}{n(n-1)}=:M.
\end{equation*}
Thus,
\begin{align*}
\bbE \left[|\bmx_k^\top \bmA\bmx_k|^2\right]
=&\sum_i|A_{ii}|^2\bbE (x_{ki}^4)
+
\sum_{i\neq l}
(A_{ii}A_{ll}^* +
|A_{il}|^2+A_{il}A_{li}^*)
M\\
=&\sum_i|A_{ii}|^2(\bbE (x_{ki}^4)-M)
+
\left(|\sum_iA_{ii}|^2+
\sum_{i\neq l}
2|A_{il}|^2\right)
M.
\end{align*}
Similarly, by using $\bbE (x_{ki}^2)=n^{-1}\sum_iT_{ik}^2=\eta_k^2$ we have
\begin{align*}
\bbE \bmx_k^\top \bmA\bmx_k \cdot\bbE \bmx_k^\top \bmA^*\bmx_k
&=\sum_i|A_{ii}|^2\left[\bbE (x_{ki}^2)\right]^2
+\sum_{i\neq l}A_{ii}A_{ll}^*\bbE (x_{ki}^2)\bbE (x_{kl}^2)\\
&=\sum_i\eta_k^4|A_{ii}|^2+\sum_{i\neq l}\eta_k^4A_{ii}A_{ll}^*
= \eta_k^4|\sum_iA_{ii}|^2.
\end{align*}
By putting these together, we have 
\begin{align*}
\var(\bmx_k^\top \bmA\bmx_k)
&=\sum_i|A_{ii}|^2\left[\bbE (x_{ki}^4)-\eta_k^4-M\right]
+|\sum_iA_{ii}|^2\left[M-\eta_k^4\right]
+\sum_{i\neq j}2|A_{ij}|^2M.
\end{align*}
For the first term, by using the assumptions that the entries of $\bmN$ have uniformly bounded fourth moments and the entries of $\eta$ are uniformly bounded, we have
\begin{align*}
\sum_i|A_{ii}|^2\left[\bbE (x_{ki}^4)-\eta_k^4-M\right]\lesssim \sum_i|A_{ii}|^2.
\end{align*}
For the second and last terms,
\begin{align*}
|M-\eta_k^4|
&=
|\frac{\eta_k^4-\sum_iT_{ik}^4}{n(n-1)}|
= O(1/n)
, &
\left|\sum_{i\neq j}
2|A_{ij}|^2
M\right|
&\lesssim \sum_{i\neq j}|A_{ij}|^2
.
\end{align*}
Thus,
\begin{equation*}
\var(\bmx_k^\top \bmA\bmx_k)\lesssim 
\sum_i|A_{ii}|^2+|\sum_{i}A_{ii}|^2/n+\sum_{i\neq j}|A_{ij}|^2
=\|\bmA\|_F^2+|\tr \bmA|^2/n.
\end{equation*}
Moreover, $|\tr \bmA|^2 = |\sum_{i=1}^n a_{ii}|^2 \le n\sum_{i=1}^n  |a_{ii}|^2 \le  n \|\bmA\|_F^2$. Hence, $\var(\bmx_k^\top \bmA\bmx_k)\lesssim \|\bmA\|_F^2$.
\end{enumerate}

\noindent
Thus, we have shown that $\bmN_\pi^\top$ satisfies
all the conditions of \cref{lemma:noise:baizhou}
and so we have from \cref{lemma:noise:baizhou}
that $p^{-1} \bmN_\pi \bmN_\pi^\top$ has a limiting spectral distribution
defined by the unique Stieltjes transform
that satisfies \cref{eq:homogenized:perm:stieltjes:yue}.
From properties of the Stieltjes transform,
it follows
that $\bmN_\pi \bmN_\pi^\top / n$ has a limiting spectral distribution
defined by the unique Stieltjes transform
that satisfies \cref{eq:marchenkopastur:stieltjes}.
Namely,
the empirical singular value distributions of $\bmN_\pi/\sqrt{n}$ and $\bbrN$
both converge to the same generalized Mar{\v{c}}enko-Pastur distribution.
This finishes the proof.
\qed


\section{Theoretical Properties of BlockFlipPA}
\label{sec:blockflippa:consistency}

This \lcnamecref{sec:blockflippa:consistency}
provides theoretical guarantees for BlockFlipPA
(introduced in \cref{sec:sim:block:dep})
that are analogous to those provided for FlipPA in \cref{sec:theory}.
Before we state the results,
we introduce the following notation
for the $(i,j)$-th $b_1 \times b_2$ block of a matrix $\bmA \in \bbR^{n \times p}$:
\begin{equation*}
    \bmA_{[i],[j]} = \bmA_{ (i-1) b_1 + 1 : i b_1, (j-1) b_2 + 1 : j b_2}
    .
\end{equation*}
For simplicity,
we will assume that $b_1$ divides $n$ and $b_2$ divides $p$.
We also introduce the following blockwise analogue of \cref{assump:noise:sym:ind}:

\begin{condition}[Noise with independent symmetric blocks]
  \label{assump:blockwise:noise:sym:ind}
  The noise matrix $\bmN$
  has independent blocks
  with symmetric distributions,
  i.e., $ N_{[i],[j]} =_d -N_{[i],[j]}$ for all $i, j$.
\end{condition}

\subsection{Analogue of \cref{prop:type-I}}

Under \cref{assump:blockwise:noise:sym:ind},
we immediately have the following analogue of \cref{prop:type-I}:

\begin{proposition}[Type~I error control for BlockFlipPA]
  \label{prop:blockflippa:type-I}
  Suppose the signal-plus-noise model~\cref{eq:signal:plus:noise}
  satisfies \cref{assump:blockwise:noise:sym:ind}.
  Then the type~I error rate of BlockFlipPA,
  for both upper-edge and pairwise comparison methods,
  is bounded above as $
  \Pr_{H_0}\bigl[\htk > 0\bigr]
  \leq 1 - \lfloor \paquant T \rfloor/(T+1)
  $.
\end{proposition}

\Cref{prop:blockflippa:type-I}
can be proved in nearly the same way as \cref{prop:type-I};
see its proof in \cref{thm:false:positive:proof}.
One need only replace all the entrywise signflip matrices
with the blockwise signflip matrices defined as
$\bmR = \btlR \otimes \bm1_{b_1 \times b_2}$, where
$\tlR_{ij} \overset{iid}{\sim}\pm1$ with probability 1/2.

\subsection{Analogue of \cref{thm:consistency}}

We now obtain the following analogue to \cref{thm:consistency},
which provided a general sufficient condition for FlipPA consistency
(of which
\cref{thm:consistency:component,thm:consistency:component:rates}
were special cases).
This \lcnamecref{thm:blockflippa:consistency}
provides a corresponding condition for BlockFlipPA consistency.

\begin{theorem}[Asymptotic consistency of BlockFlipPA: general condition]
  \label{thm:blockflippa:consistency}
  Suppose the signal-plus-noise model \cref{eq:signal:plus:noise}
  satisfies \cref{assump:blockwise:noise:sym:ind,assump:noise:upper:edge,assump:signal:perceptible}
  with fixed block sizes $b_1$ and $b_2$ as $n,p \to \infty$.
  Then BlockFlipPA using the upper-edge comparison method
  with threshold $\tau \in (0,\ep)$
  is consistent,
  i.e.,
  $\Pr\bigl[ \htk = k \bigr] \to 1$,
  as long as
  $\bbE \| \btlS \|_{2,\infty} \to 0$,
  $\bbE \| \btlS^{\top} \|_{2,\infty} \to 0$, and
  $\min\big[ \bbE \, \rho_2(\btlS), \bbE \, \rho_\infty(\btlS) \big] \to 0$,
  where
  $\btlS$ is the matrix of signal block norms
  defined as
  \begin{equation*}
    \btlS
    =
    \begin{bmatrix}
    \|\bmS_{[1],[1]}\| & \cdots & \|\bmS_{[1],[p/b_2]}\| \\
    \vdots & & \vdots \\
    \|\bmS_{[n/b_1],[1]}\| & \cdots & \|\bmS_{[n/b_1],[p/b_2]}\|
    \end{bmatrix}
    \in
    \bbR^{n/b_1 \times p/b_2}
    .
  \end{equation*}
\end{theorem}

Note that \cref{thm:blockflippa:consistency}
is identical to \cref{thm:consistency},
except that entrywise independence (\cref{assump:noise:sym:ind})
has been replaced with blockwise independence (\cref{assump:blockwise:noise:sym:ind})
and entrywise signal delocalization
has been replaced with blockwise signal delocalization.

\Cref{thm:blockflippa:consistency}
can be proved in nearly the same way as \cref{thm:consistency};
see its proof in \cref{pf:thm:consistency}.
For the most part,
one need only replace all the entrywise signflip matrices
with the blockwise signflip matrices defined as
$\bmR = \btlR \otimes \bm1_{b_1 \times b_2}$, where
$\tlR_{ij} \overset{iid}{\sim}\pm1$ with probability 1/2.
The main change required is in the analysis of signal destruction,
where \cref{lem:opnorm:bound} is used
to obtain \cref{eq:consistency:flipdataflipnoise}.
\Cref{lem:opnorm:bound} does not apply directly
to the blockwise signflip matrices used in BlockFlipPA;
the following \lcnamecref{lem:blockwise:opnorm:bound}
provides a blockwise variant.
Using this variant in place of \cref{lem:opnorm:bound}
completes the proof of \cref{thm:blockflippa:consistency}.

\begin{lemma}[Operator norms of blockwise signflipped matrices]
    \label{lem:blockwise:opnorm:bound}
    Let $\bmA \in \bbR^{n \times p}$ be arbitrary
    and let $\bmR = \btlR \otimes \bm1_{b_1 \times b_2}$
    with $\btlR \sim \operatorname{Unif}(\{-1, 1\}^{n/b_1 \times p/b_2})$
    be a blockwise Rademacher random matrix.
    Then
    $\bbE \|\bmA \circ \bmR\|
     \lesssim
     \|\btlA\|_{2,\infty}
     +
     \|\btlA^\top\|_{2,\infty}
     +
     \min\big[ \rho_2(\btlA), \rho_\infty(\btlA) \big]$,
    where
    \begin{equation*}
        \btlA
        =
        \begin{bmatrix}
            \|\bmA_{[1],[1]}\| & \cdots & \|\bmA_{[1],[p/b_2]}\| \\
            \vdots & & \vdots \\
            \|\bmA_{[n/b_1],[1]}\| & \cdots & \|\bmA_{[n/b_1],[p/b_2]}\|
        \end{bmatrix}
        \in
        \bbR^{n/b_1 \times p/b_2}
        .
    \end{equation*}
\end{lemma}

\begin{proof}[Proof of \cref{lem:blockwise:opnorm:bound}]
    Let $\tln = n/b_1$ and $\tlp = p/b_2$,
    and recall that
    \begin{equation*}
        \|\bmA \circ \bmR\|
        =
        \max_{\substack{
            \bmu \in \bbR^n : \|\bmu\|_2 = 1 \\
            \bmz \in \bbR^p : \|\bmz\|_2 = 1
        }}
        \bmu^\top (\bmA \circ \bmR) \bmz
        .
    \end{equation*}
    Now, note that
    \begin{align*}
        \bmu^\top (\bmA \circ \bmR) \bmz
        &
        =
        \begin{bmatrix}
            \bmu_{[1]} \\ \vdots \\ \bmu_{[\tln]}
        \end{bmatrix}^\top
        \begin{bmatrix}
            \tlR_{1,1} \bmA_{[1],[1]} & \cdots & \tlR_{1,\tlp} \bmA_{[1],[\tlp]} \\
            \vdots & & \vdots \\
            \tlR_{\tln,1} \bmA_{[\tln],[1]} & \cdots & \tlR_{\tln,\tlp} \bmA_{[\tln],[\tlp]}
        \end{bmatrix}
        \begin{bmatrix}
            \bmz_{[1]} \\ \vdots \\ \bmz_{[\tlp]}
        \end{bmatrix}
        \\&
        =
        \sum_{i=1}^{\tln}
        \sum_{j=1}^{\tlp}
        \bmu_{[i]}^\top \tlR_{i,j} \bmA_{[i],[j]} \bmz_{[j]}
        =
        \sum_{i=1}^{\tln}
        \sum_{j=1}^{\tlp}
        \tlR_{i,j} \bmu_{[i]}^\top \bmA_{[i],[j]} \bmz_{[j]}
        \\&
        \leq
        \sum_{i=1}^{\tln}
        \sum_{j=1}^{\tlp}
        \tlR_{i,j} \|\bmu_{[i]}\|_2 \|\bmA_{[i],[j]}\| \|\bmz_{[j]}\|_2
        ,
    \end{align*}
    where
    $\bmu_{[1]},\dots,\bmu_{[\tln]} \in \bbR^{b_1}$
    and
    $\bmz_{[1]},\dots,\bmz_{[\tlp]} \in \bbR^{b_2}$
    are the blocks of $\bmu$ and $\bmz$
    corresponding to the $b_1 \times b_2$ blocks of $\bmA \circ \bmR$.
    Thus, we have that
    \begin{align*}
        \|\bmA \circ \bmR\|
        &
        \leq
        \max_{\substack{
            \bmu \in \bbR^n : \|\bmu\|_2 = 1 \\
            \bmz \in \bbR^p : \|\bmz\|_2 = 1
        }}
        \sum_{i=1}^{\tln}
        \sum_{j=1}^{\tlp}
        \tlR_{i,j} \|\bmu_{[i]}\|_2 \|\bmA_{[i],[j]}\| \|\bmz_{[j]}\|_2
        \\&
        =
        \max_{\substack{
            \bmu \in \bbR^n : \|\bmu\|_2 = 1 \\
            \bmz \in \bbR^p : \|\bmz\|_2 = 1
        }}
        \begin{bmatrix}
            \|\bmu_{[1]}\|_2 \\ \vdots \\ \|\bmu_{[\tln]}\|_2
        \end{bmatrix}^\top
        \underbrace{
        \begin{bmatrix}
            \tlR_{1,1} \|\bmA_{[1],[1]}\| & \cdots & \tlR_{1,\tlp} \|\bmA_{[1],[\tlp]}\| \\
            \vdots & & \vdots \\
            \tlR_{\tln,1} \|\bmA_{[\tln],[1]}\| & \cdots & \tlR_{\tln,\tlp} \|\bmA_{[\tln],[\tlp]}\|
        \end{bmatrix}
        }_{=\btlA \circ \btlR}
        \begin{bmatrix}
            \|\bmz_{[1]}\|_2 \\ \vdots \\ \|\bmz_{[\tlp]}\|_2
        \end{bmatrix}
        .
    \end{align*}
    Note next that
    for any $\bmu \in \bbR^n$ such that $\|\bmu\|_2 = 1$
    and any $\bmz \in \bbR^p$ such that $\|\bmz\|_2 = 1$,
    \begin{align*}
        \left\|
            \begin{bmatrix}
                \|\bmu_{[1]}\|_2 \\ \vdots \\ \|\bmu_{[\tln]}\|_2
            \end{bmatrix}
        \right\|_2
        &
        =
        \sqrt{\|\bmu_{[1]}\|_2^2 + \cdots + \|\bmu_{[\tln]}\|_2^2}
        =
        \sqrt{\|\bmu\|_2^2}
        =
        \|\bmu\|_2
        =
        1
        , \\
        \left\|
            \begin{bmatrix}
                \|\bmz_{[1]}\|_2 \\ \vdots \\ \|\bmz_{[\tlp]}\|_2
            \end{bmatrix}
        \right\|_2
        &
        =
        \sqrt{\|\bmz_{[1]}\|_2^2 + \cdots + \|\bmz_{[\tlp]}\|_2^2}
        =
        \sqrt{\|\bmz\|_2^2}
        =
        \|\bmz\|_2
        =
        1
        ,
    \end{align*}
    and so it follows that
    \begin{equation*}
        \|\bmA \circ \bmR\|
        \leq
        \max_{\substack{
            \btlu \in \bbR^{\tln} : \|\btlu\|_2 = 1 \\
            \btlz \in \bbR^{\tlp} : \|\btlz\|_2 = 1
        }}
        \btlu^\top
        (\btlA \circ \btlR)
        \btlz
        =
        \|\btlA \circ \btlR\|
        .
    \end{equation*}
    Taking expectations and applying \cref{lem:opnorm:bound}
    completes the proof.
\end{proof}

\subsection{Analogue of \cref{thm:consistency:component}}

We finally obtain the following analogue to \cref{thm:consistency:component},
which provided a sufficient condition for FlipPA consistency
in terms of the delocalization of the signal components.
This \lcnamecref{thm:blockflippa:consistency:component} provides a corresponding condition for BlockFlipPA consistency.

\begin{theorem}[Asymptotic consistency of BlockFlipPA]
  \label{thm:blockflippa:consistency:component}
  Suppose the signal-plus-noise model~\cref{eq:signal:plus:noise}
  satisfies \cref{assump:blockwise:noise:sym:ind,assump:noise:upper:edge,assump:signal:perceptible}
  with fixed block sizes $b_1$ and $b_2$ as $n,p \to \infty$.
  Then BlockFlipPA using the upper-edge comparison method
  with threshold $\tau \in (0,\ep)$
  is consistent,
  i.e.,
  $\Pr\big[ \htk = k \big] \to 1$,
  as long as
  the signal components are delocalized as follows:
  \begin{equation}
    \label{eq:blockflippa:consistency:component}
    \bbE \left\{
      \sum_{i = 1}^{k}
      \theta_i \|\bmu_i\|_2 \|\bmz_i\|_2
      \cdot
      \left[
        \frac{\|\bmu_i\|_\infty/\|\bmu_i\|_2 + \|\bmz_i\|_\infty/\|\bmz_i\|_2}{2}
      \right]
    \right\}
    \to 0
    .
  \end{equation}
\end{theorem}

Note that \cref{thm:blockflippa:consistency:component}
is identical to \cref{thm:consistency:component},
except that entrywise independence (\cref{assump:noise:sym:ind})
has been replaced with blockwise independence (\cref{assump:blockwise:noise:sym:ind}).
Notably,
the delocalization condition \cref{eq:blockflippa:consistency:component} for BlockFlipPA
is identical to the corresponding condition \cref{eq:consistency:component} for FlipPA.
Consequently,
an analogous form of \cref{thm:consistency:component:rates} follows immediately.
We conclude this \lcnamecref{sec:blockflippa:consistency}
with a brief proof of \cref{thm:blockflippa:consistency:component}.

\begin{proof}[Proof of \cref{thm:blockflippa:consistency:component}]
\Cref{thm:blockflippa:consistency:component}
can be proved in a very similar way as \cref{thm:consistency:component};
see its proof in \cref{thm:consistency:component:proof}.
The main change is that here we build on \cref{thm:blockflippa:consistency}
in place of \cref{thm:consistency},
and so we must show that the signal matrix
$\bmS = \sum_{i = 1}^{k} \theta_i \bmu_i \bmz_i^\top$
satisfies the conditions of \cref{thm:blockflippa:consistency}
instead of the conditions of \cref{thm:consistency}.
Specifically, we need to show that
\begin{equation*}
    \bbE \| \btlS \|_{2,\infty} \to 0
    , \quad
    \bbE \| \btlS^{\top} \|_{2,\infty} \to 0
    , \quad
    \text{and} \quad
    \min\big[ \bbE \, \rho_2(\btlS), \bbE \, \rho_\infty(\btlS) \big] \to 0
    ,
\end{equation*}
where
$\btlS$ is the matrix of signal block norms
defined as
\begin{equation*}
    \btlS
    =
    \begin{bmatrix}
    \|\bmS_{[1],[1]}\| & \cdots & \|\bmS_{[1],[p/b_2]}\| \\
    \vdots & & \vdots \\
    \|\bmS_{[n/b_1],[1]}\| & \cdots & \|\bmS_{[n/b_1],[p/b_2]}\|
    \end{bmatrix}
    \in
    \bbR^{n/b_1 \times p/b_2}
    .
\end{equation*}
We begin with the first two conditions.
Note first that
\begin{align*}
    \|\btlS\|_{2,\infty}^2
    &
    =
    \max_{i=1,\dots,n/b_1} \|\btlS_{i:}\|_2^2
    =
    \max_{i=1,\dots,n/b_1}
    \|\bmS_{[i],[1]}\|^2 + \cdots + \|\bmS_{[i],[p/b_2]}\|^2
    \\&
    \leq
    \max_{i=1,\dots,n/b_1}
    \|\bmS_{[i],[1]}\|_F^2 + \cdots + \|\bmS_{[i],[p/b_2]}\|_F^2
    =
    \max_{i=1,\dots,n/b_1}
    \|\bmS_{[i],:}\|_F^2
    \\&
    \leq
    \max_{i=1,\dots,n/b_1}
    \left[
        b_1
        \cdot
        \max_{j=1,\dots,b_1}
        \|(\bmS_{[i],:})_{j:}\|_2^2
    \right]
    =
    b_1
    \cdot
    \max_{i=1,\dots,n}
    \|\bmS_{i,:}\|_2^2
    =
    b_1 \|\bmS\|_{2,\infty}^2
    .
\end{align*}
Thus,
\begin{equation*}
    \bbE \|\btlS\|_{2,\infty}
    \leq
    \sqrt{b_1} \bbE \|\bmS\|_{2,\infty}
    ,
\end{equation*}
and likewise
\begin{equation*}
    \bbE \|\btlS^\top\|_{2,\infty}
    \leq
    \sqrt{b_2} \bbE \|\bmS^\top\|_{2,\infty}
    .
\end{equation*}
In the proof of \cref{thm:consistency},
we showed that the condition \cref{eq:blockflippa:consistency:component} implies that
$\bbE \|\bmS\|_{2,\infty} \to 0$
and
$\bbE \|\bmS^\top\|_{2,\infty} \to 0$.
As a result,
we immediately have
\begin{align}
    \label{eq:blockflippa:consistency:component:proof:results:1:2}
    \bbE \| \btlS \|_{2,\infty} &\to 0
    , &
    \bbE \| \btlS^{\top} \|_{2,\infty} &\to 0
    ,
\end{align}
and are done with the first two conditions of \cref{thm:blockflippa:consistency}.

To show the final condition of \cref{thm:blockflippa:consistency},
we will show that $\bbE \|\btlS\|_{4,4} \to 0$ then apply \cref{thm:decay:suff},
as was done for $\bmS$ in the proof of \cref{thm:consistency}.
To show this, note that
\begin{align*}
    \|\btlS\|_{4,4}^4
    &
    =
    \sum_{i=1}^{n/b_1}
    \sum_{j=1}^{p/b_2}
    \|\bmS_{[i],[j]}\|^4
    \leq
    \sum_{i=1}^{n/b_1}
    \sum_{j=1}^{p/b_2}
    \|\bmS_{[i],[j]}\|_F^4
    \leq
    \sum_{i=1}^{n/b_1}
    \sum_{j=1}^{p/b_2}
    (\sqrt[4]{b_1 b_2}\|\bmS_{[i],[j]}\|_{4,4})^4
    \\&
    =
    b_1 b_2
    \sum_{i=1}^{n/b_1}
    \sum_{j=1}^{p/b_2}
    \|\bmS_{[i],[j]}\|_{4,4}^4
    =
    b_1 b_2
    \sum_{i=1}^{n/b_1}
    \sum_{j=1}^{p/b_2}
    \left[
        \sum_{\tli=1}^{b_1}
        \sum_{\tlj=1}^{b_2}
        (\bmS_{[i],[j]})_{\tli,\tlj}^4
    \right]
    =
    b_1 b_2
    \|\bmS\|_{4,4}^4
    ,
\end{align*}
and so
\begin{equation*}
    \bbE \|\btlS\|_{4,4}
    \leq
    \sqrt[4]{b_1 b_2}
    \cdot
    \bbE \|\bmS\|_{4,4}
    .
\end{equation*}
In the proof of \cref{thm:consistency},
we showed that the condition \cref{eq:blockflippa:consistency:component} implies that
$\bbE \|\bmS\|_{4,4} \to 0$,
so we immediately have
$\bbE \|\btlS\|_{4,4} \to 0$.
Applying \cref{thm:decay:suff} then yields
\begin{equation}
    \label{eq:blockflippa:consistency:component:proof:results:3}
    \min\big[ \bbE \, \rho_2(\btlS), \bbE \, \rho_\infty(\btlS) \big] \to 0
    ,
\end{equation}
so we are done with the final condition of \cref{thm:blockflippa:consistency}.

Combining \cref{eq:blockflippa:consistency:component:proof:results:1:2,eq:blockflippa:consistency:component:proof:results:3}
with \cref{thm:blockflippa:consistency} completes the proof of \cref{thm:blockflippa:consistency:component}.
\end{proof}


\section{Simulations with \texorpdfstring{$n \ll p$}{n << p}}
\label{sec:sim:highdim}

This \lcnamecref{sec:sim:highdim} repeats the experiments from
\cref{sec:sim:hom,sec:sim:block,sec:sim:block:dep}
for a higher dimensional setting with $n \ll p$,
specifically with $n = 60$ and $p = 5000$.
We consider the same methods as before,
with the exception of BEMA,
which took too long to run.

\begin{figure} \centering
  \includegraphics[scale=0.15]{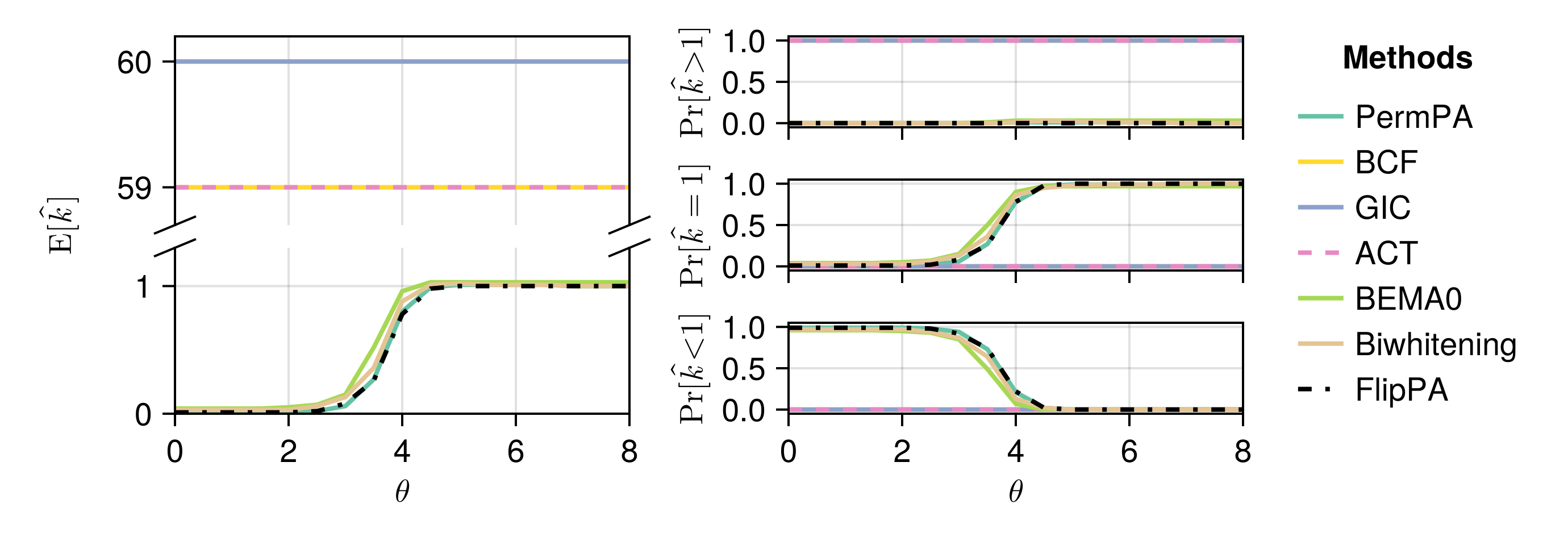}
  \caption{Higher-dimensional analogue of \cref{fig:sim:hom}
    from \cref{sec:sim:hom},
    where here $n = 60$ and $p = 5000$.
    Performance across $100$ runs
    is shown
    for each method
    for a rank-one signal in homogeneous noise
    from \cref{sec:sim:hom},
    where the signal strength $\theta$ increases
    from zero (buried in the noise)
    to eight (above the noise).
    The left plot shows the average selected rank across the runs $\bbE[\htk]$;
    the second column of plots shows what proportion
    of runs resulted in
    over-estimation $\Pr[\htk > 1]$,
    correct estimation $\Pr[\htk = 1]$,
    and
    under-estimation $\Pr[\htk < 1]$.}
  \label{fig:sim:highdim:hom}
\end{figure}

\Cref{fig:sim:highdim:hom} provides the analogue
of \cref{fig:sim:hom}
from \cref{sec:sim:hom}.
PermPA, BEMA0, Biwhitening, and FlipPA
were again highly effective
once $\theta$ was large enough,
i.e., once the signal rose above the noise,
and they estimated the rank correctly in most of the runs.
BCF, GIC, and ACT
dramatically over-estimated the rank across the entire range in this experiment;
they selected all or nearly all the components.
For $\theta=0$,
BEMA0 rejected the pure noise null hypothesis
in 4 of the 100 trials
(achieving an empirical type~I error rate of 4\%);
Biwhitening rejected the null
in 3 of the 100 trials
(achieving an empirical type~I error rate of 3\%);
PermPA and FlipPA both rejected the null
in 1 of the 100 trials
(achieving empirical type~I error rates of 1\%);
and the remaining methods (BCF, GIC, and ACT)
all rejected the null
in all of the 100 trials
(achieving empirical type~I error rates of 100\%).

\begin{figure} \centering
  \includegraphics[scale=0.15]{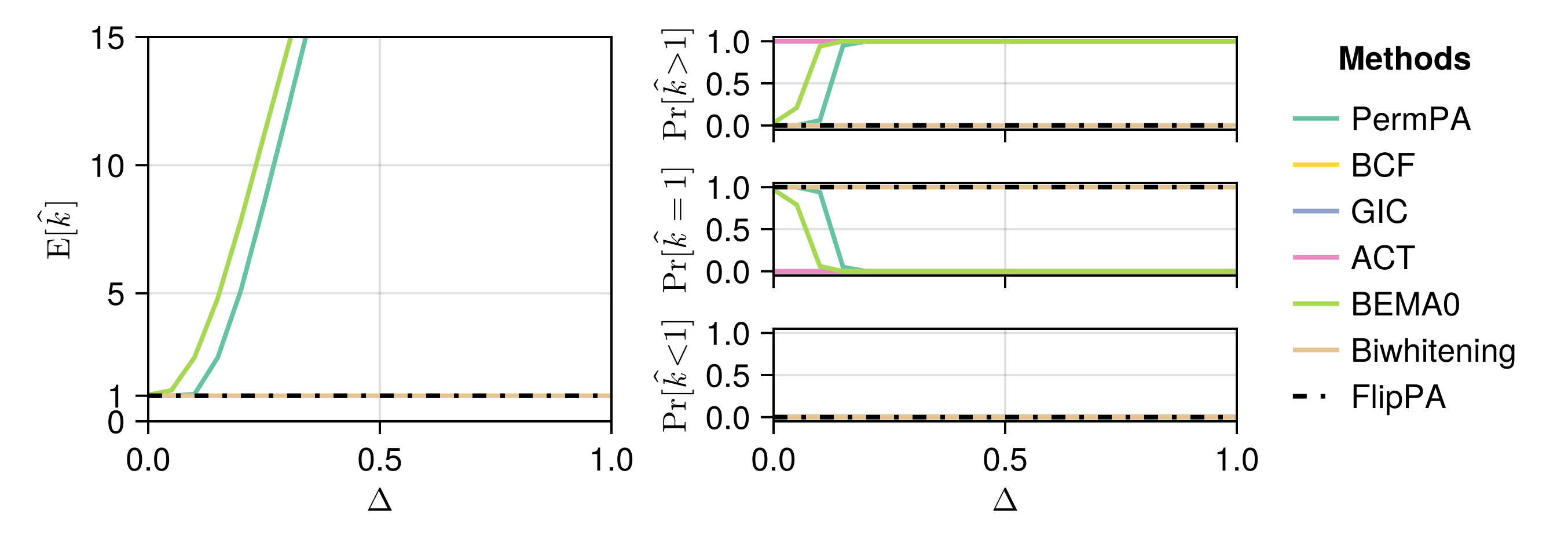}
  \caption{Higher-dimensional analogue of \cref{fig:sim:block}
    from \cref{sec:sim:block},
    where here $n = 60$ and $p = 5000$.
    Performance across $100$ runs
    is shown for each method for
    a rank-one signal in noise having a
    block-structured noise variance profile
    \cref{eq:sim:block:noise},
    where the amount of heteroscedasticity increases
    as $\Delta$ increases
    from zero (homoscedastic noise)
    to one (maximal heteroscedasticity).
    The left plot shows the average selected rank across the runs $\bbE[\htk]$;
    the second column of plots shows what proportion
    of runs resulted in
    over-estimation $\Pr[\htk > 1]$,
    correct estimation $\Pr[\htk = 1]$,
    and
    under-estimation $\Pr[\htk < 1]$.
    BCF, GIC, and ACT do not appear in the left plot
    because their averages are close to $60$ for the entire sweep.}
  \label{fig:sim:highdim:block}
\end{figure}

\Cref{fig:sim:highdim:block} provides the analogue
of \cref{fig:sim:block}
from \cref{sec:sim:block},
where here $\theta = 12$, $n_1 = n_2 = 30$, $p_1 = 4000$, and $p_2 = 1000$.
With the exceptions of BCF, GIC, and ACT
(all of which dramatically over-estimated the rank again),
the remaining methods (PermPA, BEMA0, Biwhitening, and FlipPA)
performed similarly to
the lower dimensional setting in \cref{fig:sim:block}.
They were highly effective when $\Delta$ was small,
but only FlipPA and Biwhitening remained effective
as $\Delta$ grew and the noise became more heteroscedastic.
PermPA and BEMA0 again overestimated the rank as $\Delta$ grew;
notably, their performance
degraded more rapidly here than in \cref{fig:sim:block}.

\begin{figure} \centering
  \includegraphics[scale=0.15]{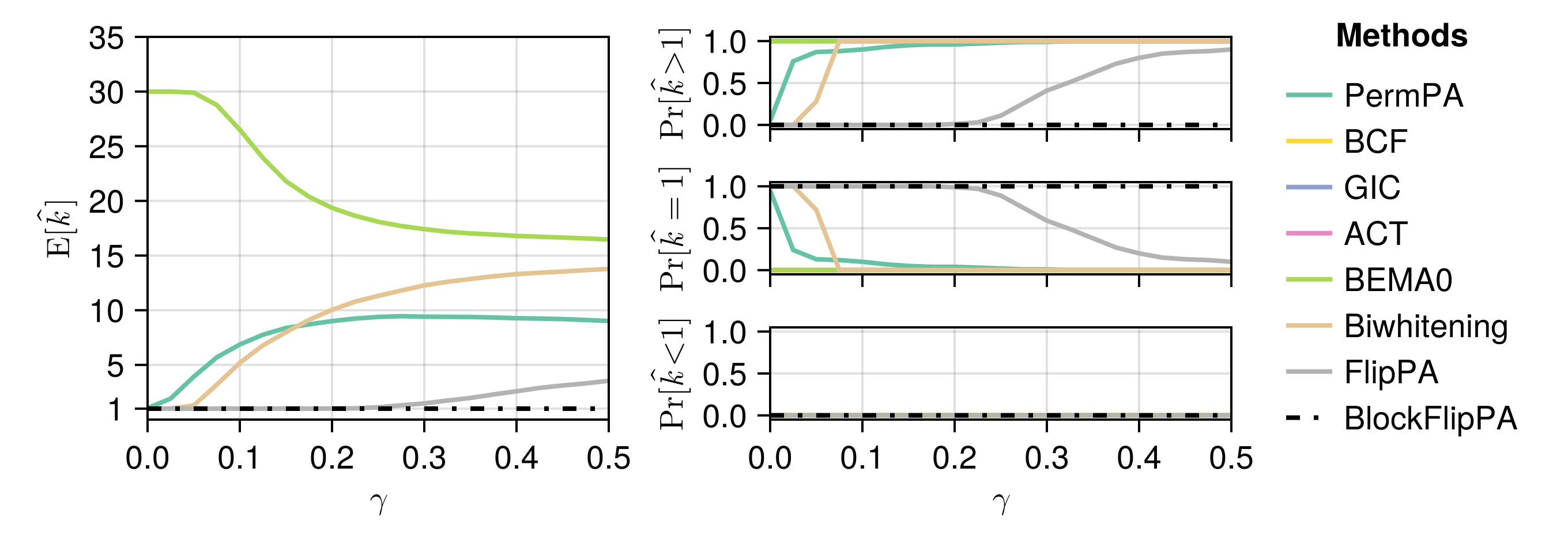}
  \caption{Higher-dimensional analogue of \cref{fig:sim:block:dep}
    from \cref{sec:sim:block:dep},
    where here $n = 60$ and $p = 5000$.
    Performance across $100$ runs
    is shown for each method for
    a rank-one signal in noise having a
    block-structured noise variance profile
    with blockwise dependence,
    where the amount of dependence increases
    as $\gamma$ increases
    from zero (independent entries)
    to $1/2$ (increasing blockwise dependence).
    The left plot shows the average selected rank across the runs $\bbE[\htk]$;
    the second column of plots shows what proportion
    of runs resulted in
    over-estimation $\Pr[\htk > 1]$,
    correct estimation $\Pr[\htk = 1]$,
    and
    under-estimation $\Pr[\htk < 1]$.
    BCF, GIC, and ACT do not appear in the left plot
    because their averages are all above $35$ for the entire sweep;
    in fact BCF and GIC are close to $60$ for the entire sweep.}
  \label{fig:sim:highdim:block:dep}
\end{figure}

\Cref{fig:sim:highdim:block:dep} provides the analogue
of \cref{fig:sim:block:dep}
from \cref{sec:sim:block:dep},
where here $\theta = 48$, $n_1 = n_2 = 30$, $p_1 = 4000$, $p_2 = 1000$, $b_1 = 3$, and $b_2 = 250$.
As in the lower dimensional setting of \cref{fig:sim:block:dep},
a few methods (PermPA, Biwhitening, FlipPA, and BlockFlipPA) were effective when $\gamma$ was small,
but only BlockFlipPA remained effective across the entire sweep.
The rest dramatically over-estimated the rank across the entire sweep.
Interestingly, FlipPA was significantly more robust to the growing dependence here
than it was in the lower dimensional setting in \cref{fig:sim:block:dep}.

\begin{figure} \centering
  \includegraphics[scale=0.15]{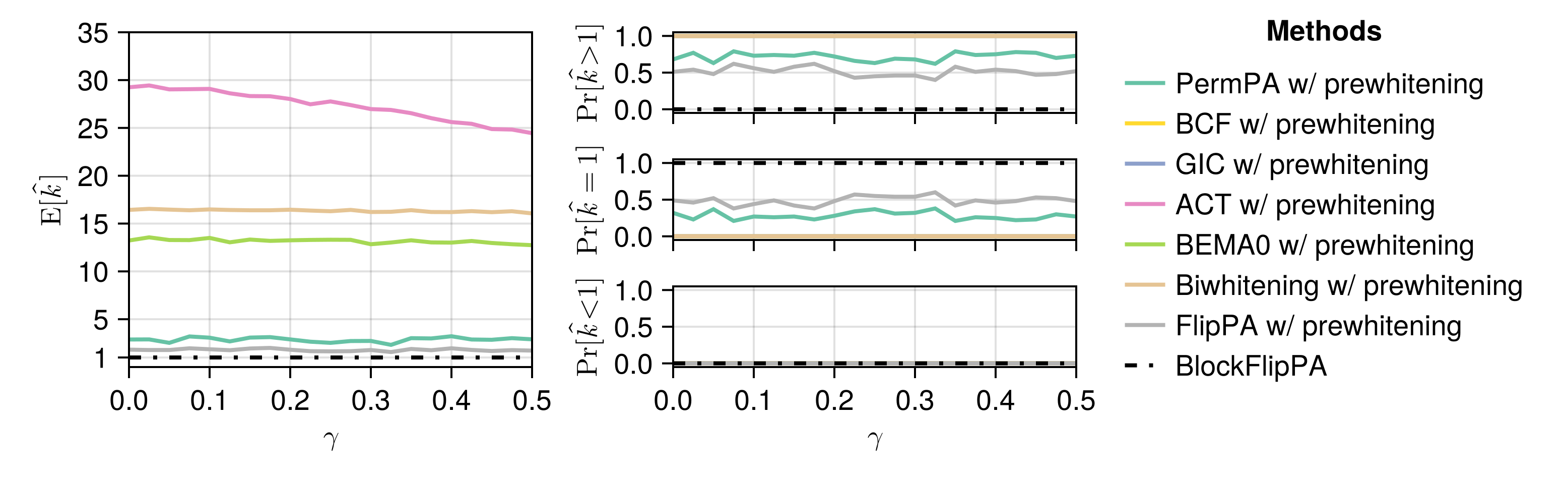}
  \caption{Higher-dimensional analogue of \cref{fig:sim:block:dep:prewhiten:mask}
    from \cref{sec:sim:block:dep},
    where here $n = 60$ and $p = 5000$.
    Performance is shown
    for the same setting as \cref{fig:sim:highdim:block:dep}
    across $100$ runs of each method
    with prewhitening using the estimates $\bhtSigma_1,\bhtSigma_2$.
    The left plot shows the average selected rank across the runs $\bbE[\htk]$;
    the second column of plots shows what proportion
    of runs resulted in
    over-estimation $\Pr[\htk > 1]$,
    correct estimation $\Pr[\htk = 1]$,
    and
    under-estimation $\Pr[\htk < 1]$.
    BCF and GIC do not appear in the left plot
    because their averages are close to $60$ for the entire sweep.}
  \label{fig:sim:highdim:block:dep:prewhiten:mask}
\end{figure}

\Cref{fig:sim:highdim:block:dep:prewhiten:mask} provides the analogue
of \cref{fig:sim:block:dep:prewhiten:mask}
from \cref{sec:sim:block:dep},
where here $\theta = 48$, $n_1 = n_2 = 30$, $p_1 = 4000$, $p_2 = 1000$, $b_1 = 3$, and $b_2 = 250$.
In contrast to the lower dimensional setting of \cref{fig:sim:block:dep},
prewhitening did not significantly improve the performance of most of the methods.
BCF, GIC, ACT, and BEMA0 continued to dramatically over-estimate the rank.
Biwhitening actually performed worse with prewhitening than without,
perhaps due to poor estimation of the prewhitening matrices
in this higher dimensional setting.
Interestingly,
the performance of PermPA and FlipPA improved with prewhitening for large values of $\gamma$,
but FlipPA performed worse with prewhitening for small values of $\gamma$.
Only BlockFlipPA was effective across the sweep.


\section{Simulations studying strong signal shadowing}
\label{sec:sim:shadowing}

This \lcnamecref{sec:sim:shadowing} investigates
the parallel analysis phenomenon known as ``shadowing'',
where strong signals can cause weak signals to be missed.
Roughly speaking,
insufficient destruction of strong signals
(e.g., by permutations in PermPA or by signflipping in FlipPA)
can result in an inflated estimate of the noise floor,
which in turn can cause weak signals
(below the inflated noise floor)
to be missed.
To study this phenomenon,
we \edit{begin by considering} the setting of \cref{fig:sim:hom} from \cref{sec:sim:hom}
but now with a rank-$2$ signal.
Namely, the $n \times p$ data matrix $\bmX$ is
generated as $\bmX = \theta_1 \bmu_1 \bmz_1^\top + \theta_2 \bmu_2 \bmz_2^\top + \bmN$,
where $N_{ij} \overset{iid}{\sim} \clN(0, v/n)$,
$\bmu_1, \bmu_2 \in \bbR^n$ and $\bmz_1, \bmz_2 \in \bbR^p$
are drawn uniformly from the respective unit spheres,
$\theta_1$ and $\theta_2$ are swept from zero to twenty,
and we take $v = 1$ without loss of generality.
As $\theta_1$ and $\theta_2$ increase,
the corresponding signals transition from being buried in the noise
to rising above it (roughly around a value of one)
then to finally being well above the noise.

\begin{figure} \centering
  \includegraphics[scale=0.15]{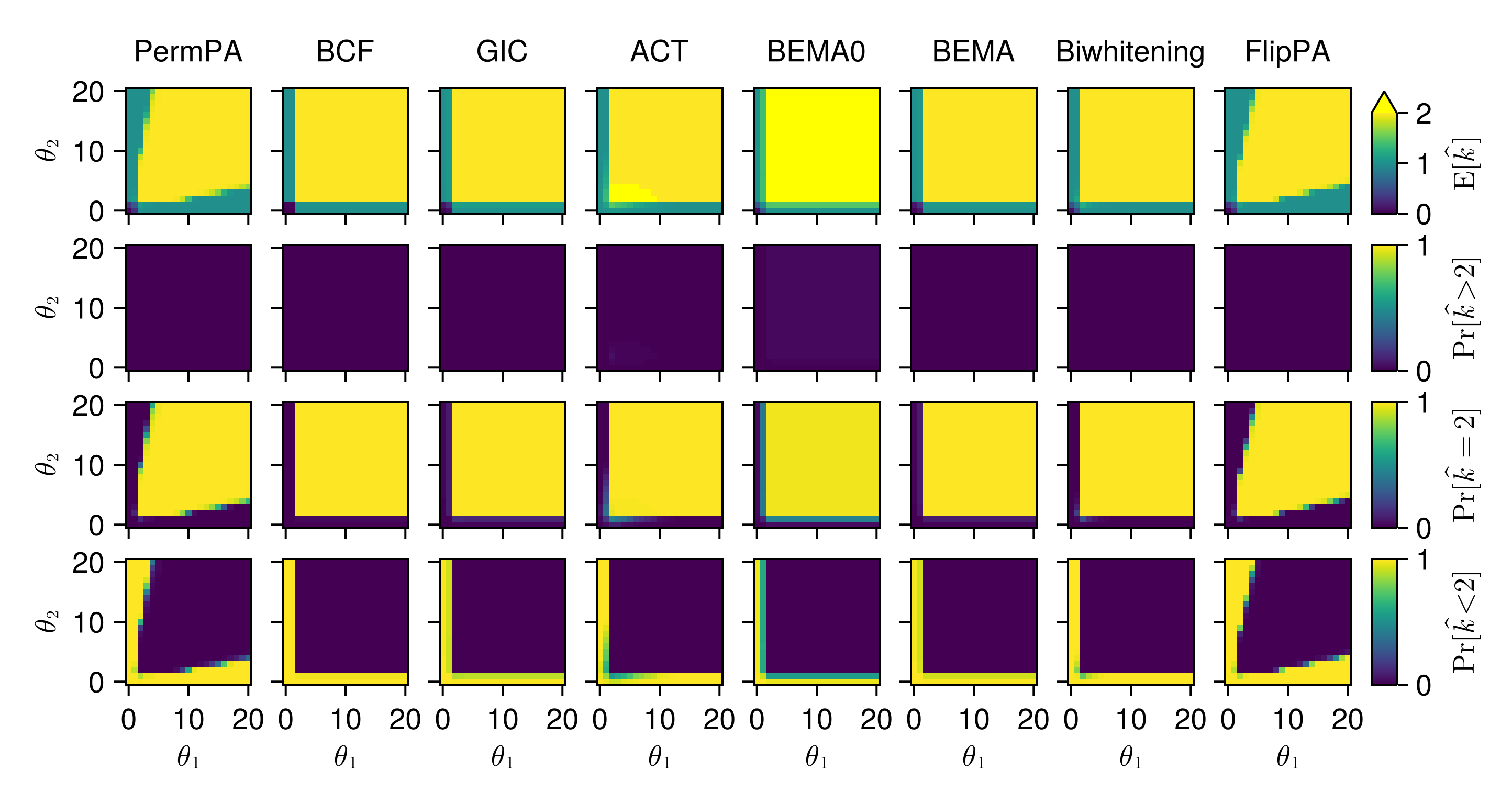}
  \caption{Rank-two analogue of \cref{fig:sim:hom}
    from \cref{sec:sim:hom}.
    Performance across $100$ runs
    is shown
    for each method,
    where the signal strengths $\theta_1$ and $\theta_2$ range
    from zero (buried in the noise)
    to twenty (well above the noise).
    The first row of heatmaps shows the average selected rank across the runs $\bbE[\htk]$;
    the remaining rows show what proportion
    of runs resulted in
    over-estimation $\Pr[\htk > 2]$,
    correct estimation $\Pr[\htk = 2]$,
    and
    under-estimation $\Pr[\htk < 2]$.}
  \label{fig:sim:shadowing}
\end{figure}

\Cref{fig:sim:shadowing} shows the resulting performance of each method
across $100$ runs.
The first row of heatmaps shows the average selected rank $\bbE\bigl[\htk\bigr]$ across the runs;
the remaining rows show the proportion $\Pr\bigl[\htk > 2\bigr]$
of runs resulting in over-estimation,
correct estimation (i.e., $\Pr\bigl[\htk = 2\bigr]$),
and under-estimation (i.e., $\Pr\bigl[\htk < 2\bigr]$).
As before, all the methods were highly effective
when both $\theta_1$ and $\theta_2$ were large,
i.e., when both signals rose above the noise
and neither dominated the other.
When either $\theta_1$ or $\theta_2$ was small,
the corresponding signal was buried in the noise
and none of the methods correctly found it.
The final regime,
where both signals rise above the noise
and one dominates the other,
is where shadowing occurs.
As expected,
both parallel analysis methods (PermPA and FlipPA) exhibited shadowing
in this regime,
i.e., they only identified one of the two signals.
The remaining methods did not suffer from shadowing;
they correctly identified both signals
even when one was much stronger than the other.
Combining these methods with FlipPA
to combine their resilience to shadowing
with the strengths of FlipPA (type~I error control, robustness to heterogeneous noise, etc.)
is an interesting direction for future work.

\edit{
Notably,
shadowing only occurred when the stronger signal
was over four times stronger than the weaker signal.
When both of the signals were strong,
both of the parallel analysis methods correctly identified them.
This matches the findings of the theoretical analysis in \cref{sec:relative:strength}.
To further investigate this mixed strong/weak factor regime,
we now consider a more focused version of the preceding sweep,
where one component is fixed at a strong-factor scale
and the other component is swept from buried to clearly detectable.
Specifically,
we again take $n=500$, $p=300$, and $N_{ij}\overset{iid}{\sim}\clN(0,1/n)$,
but now fix $\theta_1=20$ and sweep $\theta_2$ from zero to ten.
This means that in the factor model scaling $\sqrt{n}\bmX$ of \cite{bai2002determining},
the first signal has a strength of $\theta_1\sqrt{n}\approx447$,
which is comparable to $\sqrt{np}\approx387$.
This places us in the strong-factor regime for this factor model.

\begin{figure} \centering
  \includegraphics[scale=0.15]{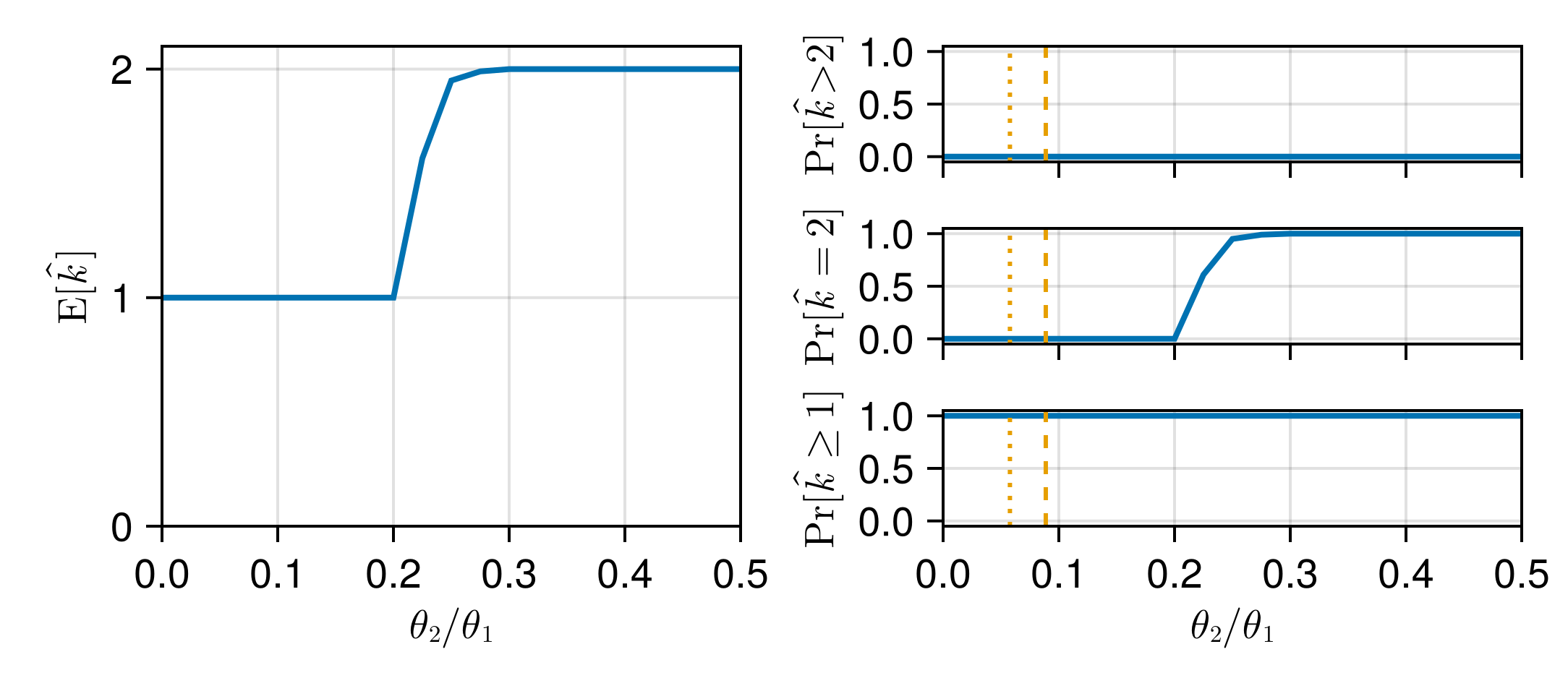}
  \caption{\edit{Focused experiment in the mixed strong/weak factor regime.
    The data are generated as in \cref{fig:sim:shadowing},
    with $\theta_1=20$ fixed and $\theta_2$ swept from zero to ten.
    The performance of FlipPA across $100$ runs is shown as a function of the relative strength $\theta_2/\theta_1$.
    The dotted vertical line marks $1/\sqrt{\min(n,p)}$,
    and the dashed vertical line marks
    $(1+\sqrt{p/n})/\theta_1$, the approximate Gaussian noise-edge scale
    relative to the strong factor.}}
  \label{fig:sim:mixed:strong:weak}
\end{figure}

\Cref{fig:sim:mixed:strong:weak} shows the performance of FlipPA across $100$ runs.
Note first that it selected the strong factor
in all $100$ runs throughout the sweep, even when the second factor was absent
or buried in the noise.
Moreover, there was no over-selection in any of the runs.
The weaker factor was selected increasingly often once its relative strength
was large enough compared with the empirical-null scale induced by the
signflipped strong factor:
FlipPA selected both factors
in $61\%$ of the runs at $\theta_2/\theta_1=0.225$,
in $95\%$ of the runs at $\theta_2/\theta_1=0.25$,
in $99\%$ of the runs at $\theta_2/\theta_1=0.275$,
and in all the runs once $\theta_2/\theta_1\geq0.3$.
This supports the relative-strength interpretation (which was analyzed in \cref{sec:relative:strength}):
strong dense factors do not prevent FlipPA from selecting strong components,
and weaker components are recovered once they rise sufficiently above
the combined noise and signflipped-signal scale.

\begin{figure} \centering
  \includegraphics[scale=0.15]{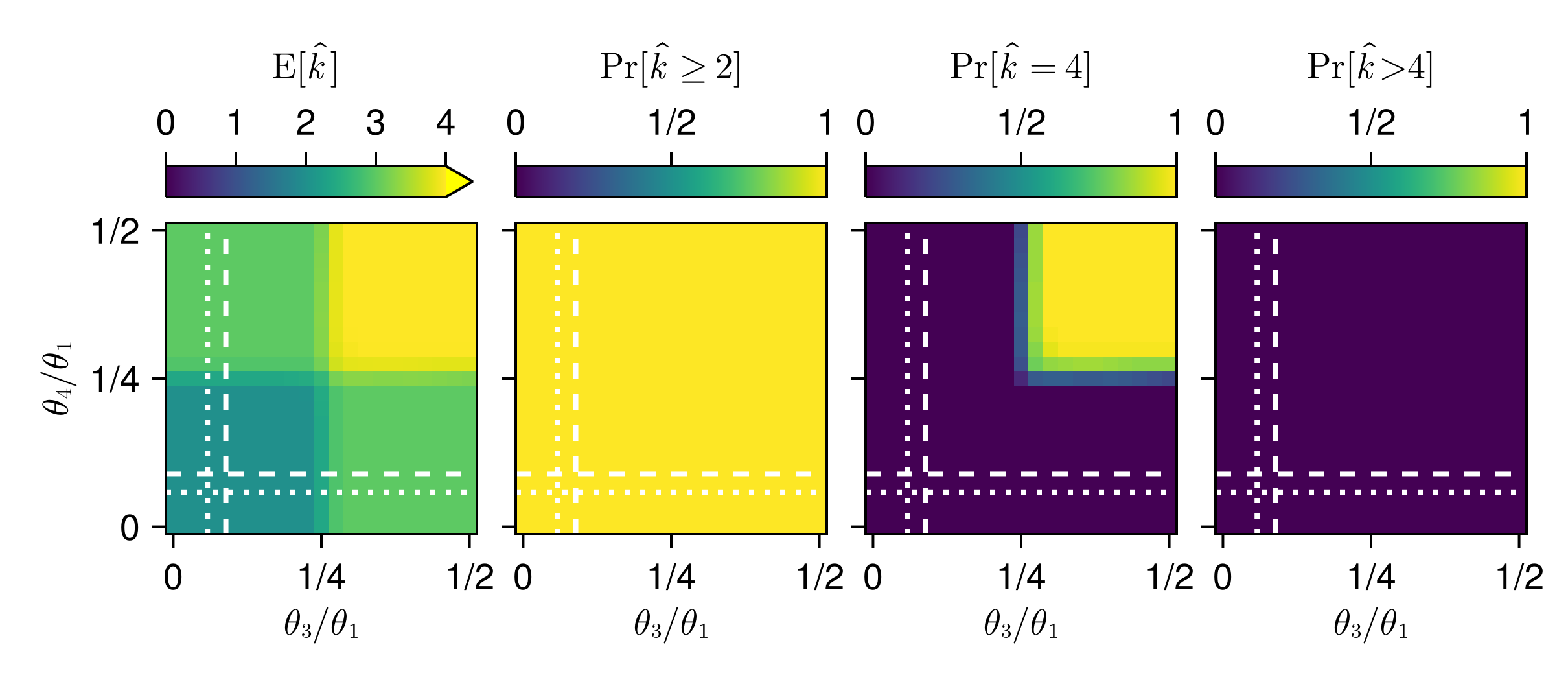}
  \caption{\edit{Experiment with two-strong and two-weak factors.
    The data are generated as in \cref{fig:sim:shadowing},
    but with four factors for which $\theta_1=\theta_2=20$ is fixed
    and $\theta_3,\theta_4$ are both swept from zero to ten.
    The performance of FlipPA across $100$ runs is shown as a function of the relative strengths $\theta_3/\theta_1$ and $\theta_4/\theta_1$.
    The dotted reference lines mark
    $1/\sqrt{\min(n,p)}$,
    and the dashed reference lines mark
    $(1+\sqrt{p/n})/\theta_1$.}}
  \label{fig:sim:mixed:two:strong:two:weak}
\end{figure}

To ensure that the above conclusions are not merely an artifact of having only one strong and one weak factor, \cref{fig:sim:mixed:two:strong:two:weak} repeats the experiment but with two strong factors
and two weak factors.
The first two signal strengths are fixed at $\theta_1=\theta_2=20$,
while the final two signal strengths $\theta_3$ and $\theta_4$ are both swept from zero to ten.
Here we observe the same relative-strength
behavior in the presence of multiple strong factors as we did with only one.
When either weak factor was absent or too small,
FlipPA typically selected the two strong factors and,
when present, the larger weak factor,
but did not select all four.
Once both weak factors were above the signflipped strong-factor scale,
the probability of selecting all four factors rose sharply:
on the grid points with both $\theta_3$ and $\theta_4$ above $6.5$
(i.e., with $\theta_j/\theta_1\geq0.325$),
all the runs correctly selected all four components.
Overall, we find that the above conclusions carry over to scenarios with multiple strong factors and multiple weak factors.
}


\section{Simulations with asymmetric noise}
\label{sec:sim:asymmetric}

This \lcnamecref{sec:sim:asymmetric}
investigates the impact of asymmetry in the noise distribution.
As discussed in \cref{rem:noise:sym},
FlipPA can achieve consistency even
if the noise entries have asymmetric distributions.
One only needs the signflipped noise $\bmR \circ \bmN$
to share the same
asymptotic upper-edge $\brsigma$ as $\bmN$ (from \cref{assump:noise:upper:edge}),
which often occurs for large random matrices
due to universality.
This \lcnamecref{sec:sim:asymmetric}
investigates the nonasymptotic behavior of FlipPA
under asymmetric noise
via experiments with Bernoulli, Poisson and Skew Normal noise.

\begin{figure} \centering
  \begin{subfigure}{\linewidth} \centering
    \includegraphics[scale=0.14, trim=0 50 0 0, clip]{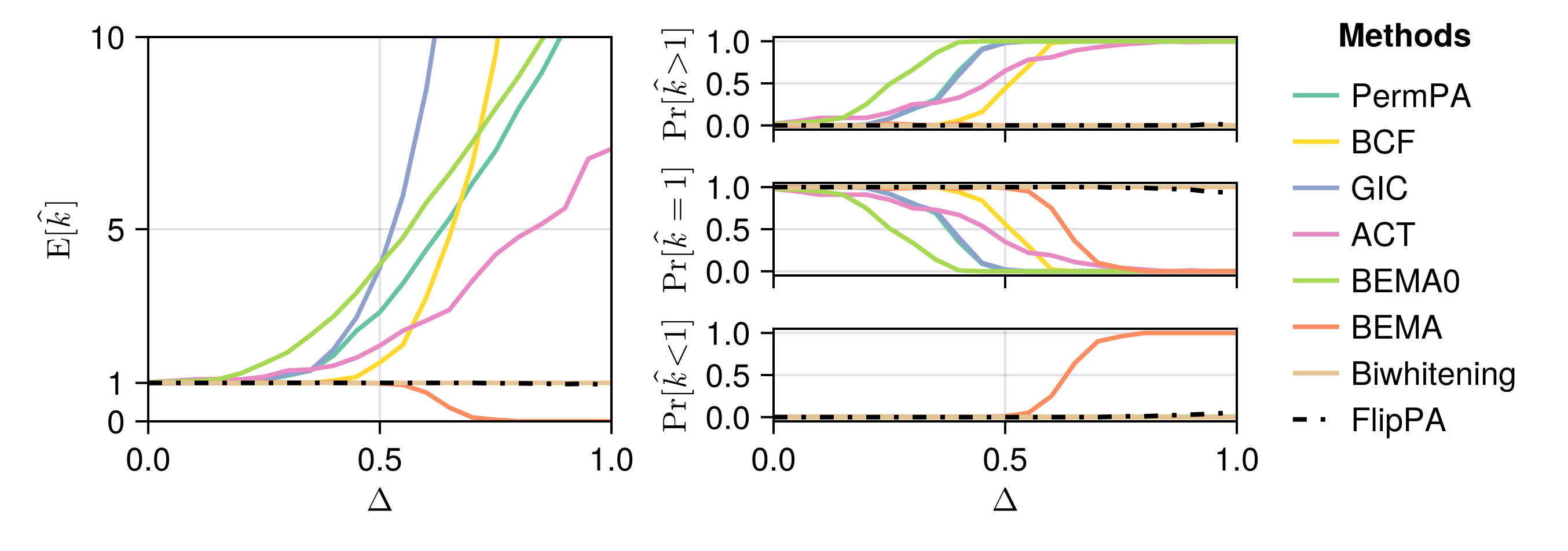}
    \caption{Centered Bernoulli noise \cref{eq:sim:centered:bernoulli}}
    \label{fig:sim:block:centered:bernoulli}
  \end{subfigure}
  \begin{subfigure}{\linewidth} \centering
    \includegraphics[scale=0.14, trim=0 50 0 0, clip]{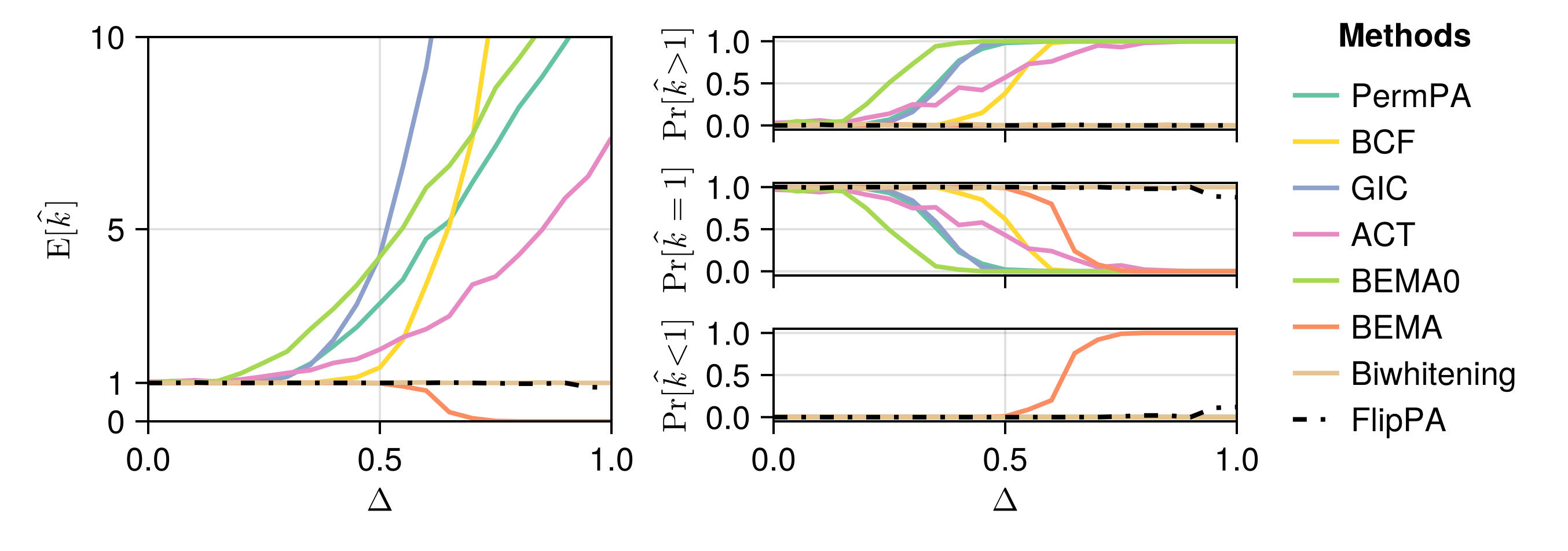}
    \caption{Centered Poisson noise \cref{eq:sim:centered:poisson}}
    \label{fig:sim:block:centered:poisson}
  \end{subfigure}
  \begin{subfigure}{\linewidth} \centering
    \includegraphics[scale=0.14, trim=0 50 0 0, clip]{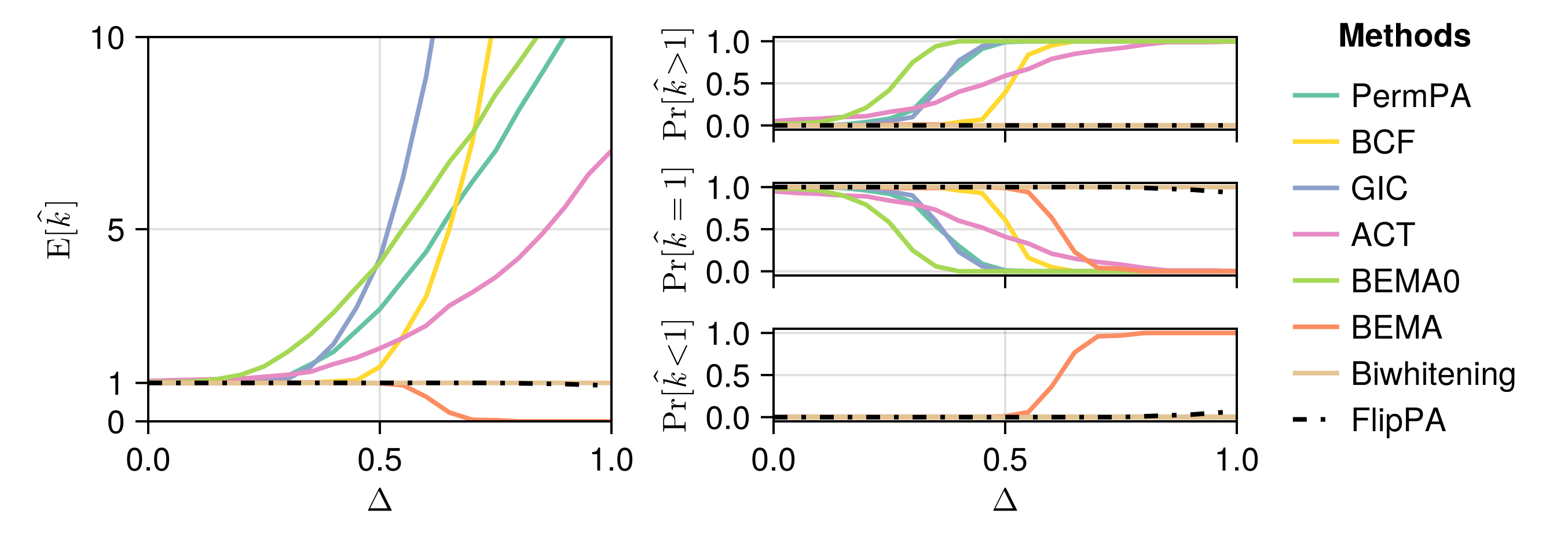}
    \caption{Skew Normal noise \cref{eq:sim:skew:normal} with $\alpha = 2$}
    \label{fig:sim:block:skew:normal:2}
  \end{subfigure}
  \begin{subfigure}{\linewidth} \centering
    \includegraphics[scale=0.14, trim=0 50 0 0, clip]{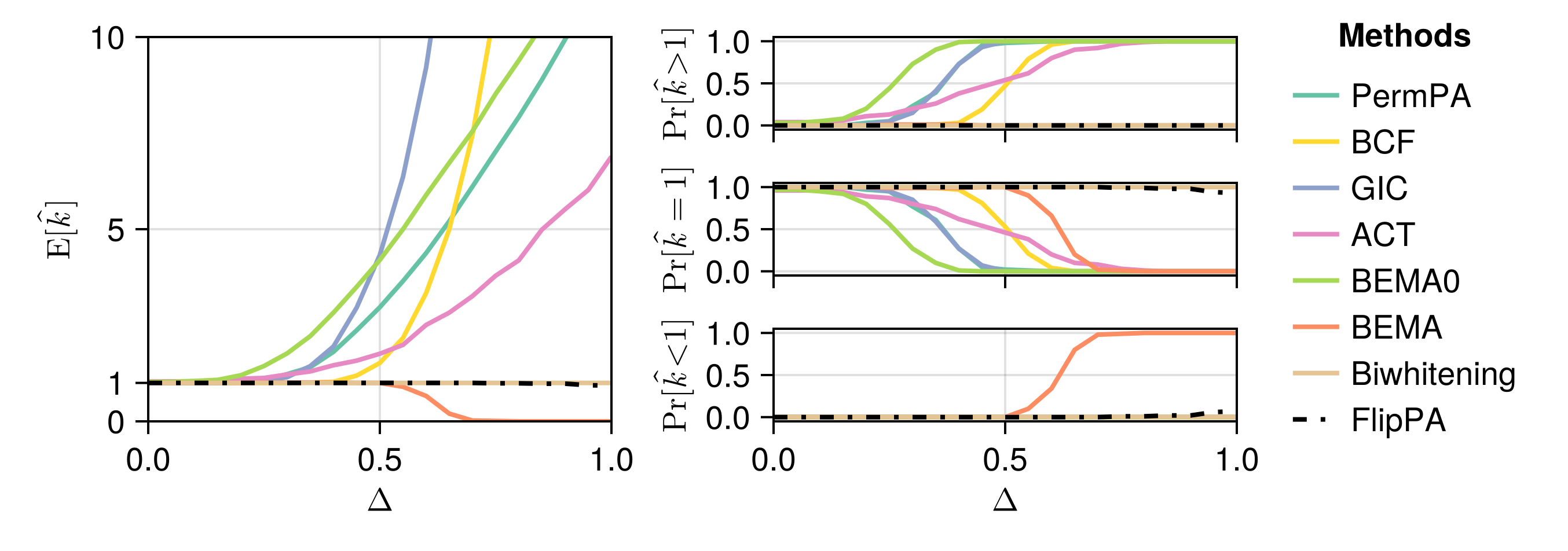}
    \caption{Skew Normal noise \cref{eq:sim:skew:normal} with $\alpha = 10$}
    \label{fig:sim:block:skew:normal:10}
  \end{subfigure}
  \caption{Analogue of \cref{fig:sim:block}
    from \cref{sec:sim:block}
    for asymmetric noise distributions.
    The performance across $100$ runs is shown for each method,
    where the noise has a block-structured variance profile
    as in \cref{eq:sim:block:noise}
    with the amount of heteroscedasticity parameterized by $\Delta$.
    The performance of all the methods is nearly identical across the distributions.}
  \label{fig:sim:block:asymmetric}
\end{figure}

\begin{figure} \centering
  \includegraphics[scale=0.14, trim=0 50 0 0, clip]{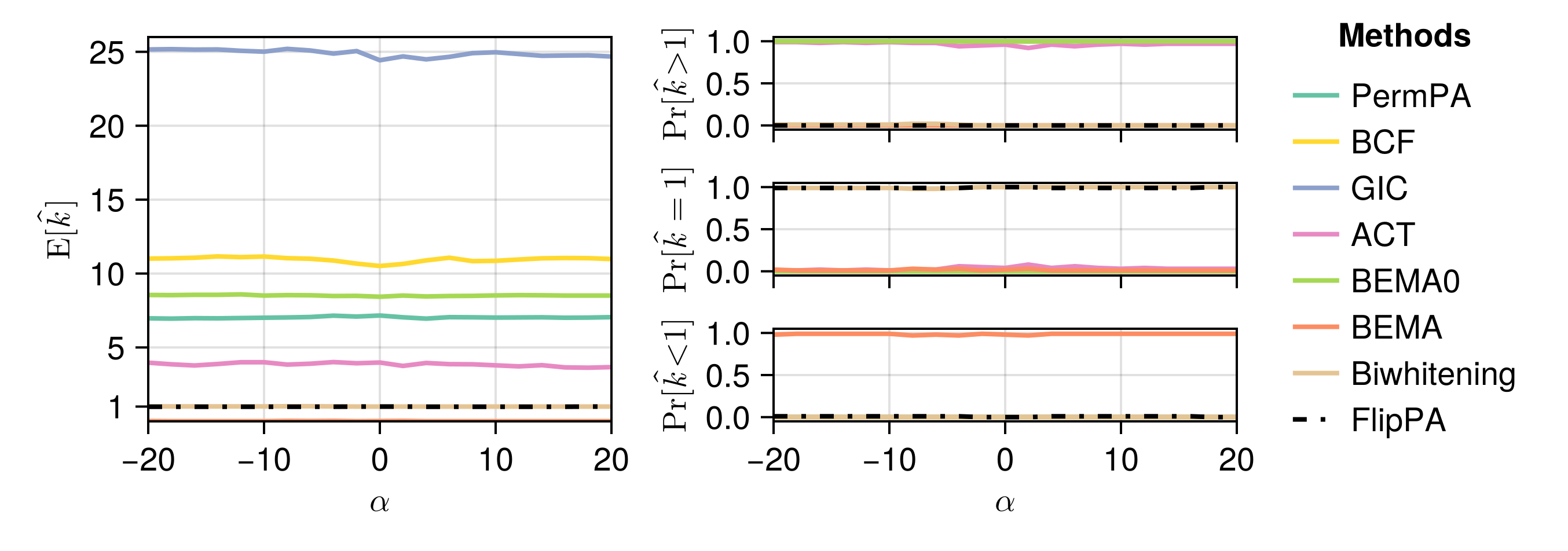}
  \caption{Performance across $100$ runs of each method for
    a rank-one signal in Skew Normal noise \cref{eq:sim:skew:normal},
    where we sweep the shape parameter $\alpha$,
    which determines the skewness of the distribution.
    The noise variance profile is
    as in \cref{eq:sim:block:noise},
    where here we set $\Delta = 0.75$.}
  \label{fig:sim:block:skew:sweep}
\end{figure}

\Cref{fig:sim:block:asymmetric} repeats the experiment
of \cref{fig:sim:block} from \cref{sec:sim:block},
where the noise variance profile is as in \cref{eq:sim:block:noise}
but now with
centered Bernoulli, Poisson, and Skew Normal noise distributions,
all of which are asymmetric distributions.
\Cref{fig:sim:block:centered:bernoulli}
considers (centered and scaled) Bernoulli noise
generated as follows
\begin{equation}
  \label{eq:sim:centered:bernoulli}
  N_{ij}
  \overset{ind}{\sim}
  \sqrt{\frac{8}{n}}
  \Big( \operatorname{Bernoulli}(P_{ij}) - P_{ij} \Big)
  ,
\end{equation}
where
\begin{equation*}
    P_{ij}
    =
    \frac{1}{2}
    +
    \frac{1}{2}
    \sqrt{1 - \frac{V_{ij}}{2}}
    .
\end{equation*}
This produces centered Bernoulli noise entries
with means of $\bbE N_{ij} = \sqrt{8/n}(P_{ij} - P_{ij}) = 0$,
variances of
\begin{equation*}
  \var(N_{ij})
  =
  \frac{8}{n} P_{ij} (1 - P_{ij})
  =
  \frac{8}{n}
  \Bigg(
    \frac{1}{2}
    +
    \frac{1}{2}
    \sqrt{1 - \frac{V_{ij}}{2}}
  \Bigg)
  \Bigg(
    \frac{1}{2}
    -
    \frac{1}{2}
    \sqrt{1 - \frac{V_{ij}}{2}}
  \Bigg)
  =
  \frac{V_{ij}}{n}
  ,
\end{equation*}
and skewness given by
\begin{equation*}
  \skewness(N_{ij})
  =
  \frac{1 - 2 P_{ij}}{\sqrt{P_{ij}(1-P_{ij})}}
  =
  -2\sqrt{2/V_{ij} - 1}
  .
\end{equation*}
\Cref{fig:sim:block:centered:poisson}
considers (centered and scaled) Poisson noise
generated as follows
\begin{equation}
  \label{eq:sim:centered:poisson}
  N_{ij}
  \overset{ind}{\sim}
  \frac{1}{\sqrt{n}}
  \Big( \operatorname{Poisson}(V_{ij}) - V_{ij} \Big)
  ,
\end{equation}
which produces centered Poisson noise entries
with means of $\bbE N_{ij} = (1/\sqrt{n}) (V_{ij} - V_{ij}) = 0$,
variances of $\var(N_{ij}) = V_{ij}/n$,
and skewness given by $\skewness(N_{ij}) = 1/\sqrt{V_{ij}}$.
\Cref{fig:sim:block:skew:normal:2,fig:sim:block:skew:normal:10}
consider Skew Normal noise
generated as follows
\begin{equation}
  \label{eq:sim:skew:normal}
  N_{ij} \overset{ind}{\sim} \operatorname{SkewNormal}\Big(
    \xi \sqrt{V_{ij} / n},
    \omega \sqrt{V_{ij} / n},
    \alpha
  \Big)
  ,
\end{equation}
where
$\xi = -\omega \delta \sqrt{2 / \pi}$,
$\omega = 1/\sqrt{1 - 2\delta^2 / \pi}$,
$\delta = \alpha / \sqrt{1 + \alpha^2}$,
and we take $\alpha = 2$ in \cref{fig:sim:block:skew:normal:2}
and $\alpha = 10$ in \cref{fig:sim:block:skew:normal:10}
(note that setting $\alpha = 0$ would produce Gaussian noise entries).
This produces Skew Normal noise entries
with means of
\begin{equation*}
  \bbE N_{ij}
  =
  \xi
  \sqrt{V_{ij} / n}
  +
  \omega
  \sqrt{V_{ij} / n}
  \delta
  \sqrt{2 / \pi}
  =
  (\xi + \omega \delta \sqrt{2 / \pi})
  \sqrt{V_{ij} / n}
  =
  0
  ,
\end{equation*}
variances of
\begin{equation*}
  \var(N_{ij})
  =
  \Big( \omega \sqrt{V_{ij} / n} \Big)^2
  \left( 1 - \frac{2\delta^2}{\pi} \right)
  =
  \left[
    \omega^2
    \left( 1 - \frac{2\delta^2}{\pi} \right)
  \right]
  \frac{V_{ij}}{n}
  =
  \frac{V_{ij}}{n}
  ,
\end{equation*}
and skewness given by
\begin{equation*}
  \skewness(N_{ij})
  =
  \frac{4 - \pi}{2}
  \frac{(\delta \sqrt{2/\pi})^3}{(1-2\delta^2/\pi)^{3/2}}
  .
\end{equation*}
In all four cases,
all the methods performed similarly
to the Gaussian setting considered in \cref{fig:sim:block}.
All the methods were effective for small $\Delta$
(where the noise is close to homoscedastic)
but only FlipPA and Biwhitening
correctly estimated the rank across
the entire sweep.

To study the impact of skewness,
\cref{fig:sim:block:skew:sweep}
considers the Skew Normal noise \cref{eq:sim:skew:normal}
used in \cref{fig:sim:block:skew:normal:2,fig:sim:block:skew:normal:10}
but varies the shape parameter $\alpha$
instead of the heterogeneity parameter $\Delta$,
which is set to $\Delta = 0.75$.
As $\alpha$ is swept from $-20$ to $20$,
the noise entries go from being left-skewed
to right-skewed.
Notably, the performance of the methods
did not vary with the skewness.

\begin{figure} \centering
  \begin{subfigure}{\linewidth} \centering
    \includegraphics[scale=0.13, trim=0 50 0 0, clip]{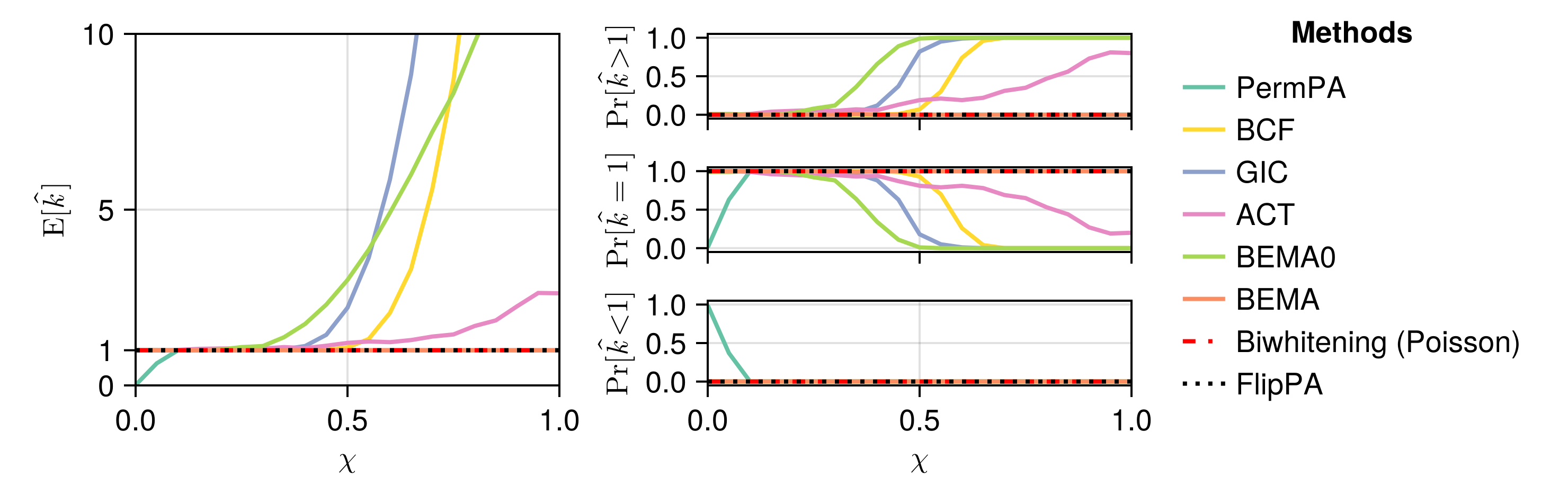}
    \caption{Bernoulli data \cref{eq:sim:rankone:bernoulli}
      with rank-one parameter matrix}
    \label{fig:sim:rankone:bernoulli}
  \end{subfigure}
  \begin{subfigure}{\linewidth} \centering
    \includegraphics[scale=0.13, trim=0 50 0 0, clip]{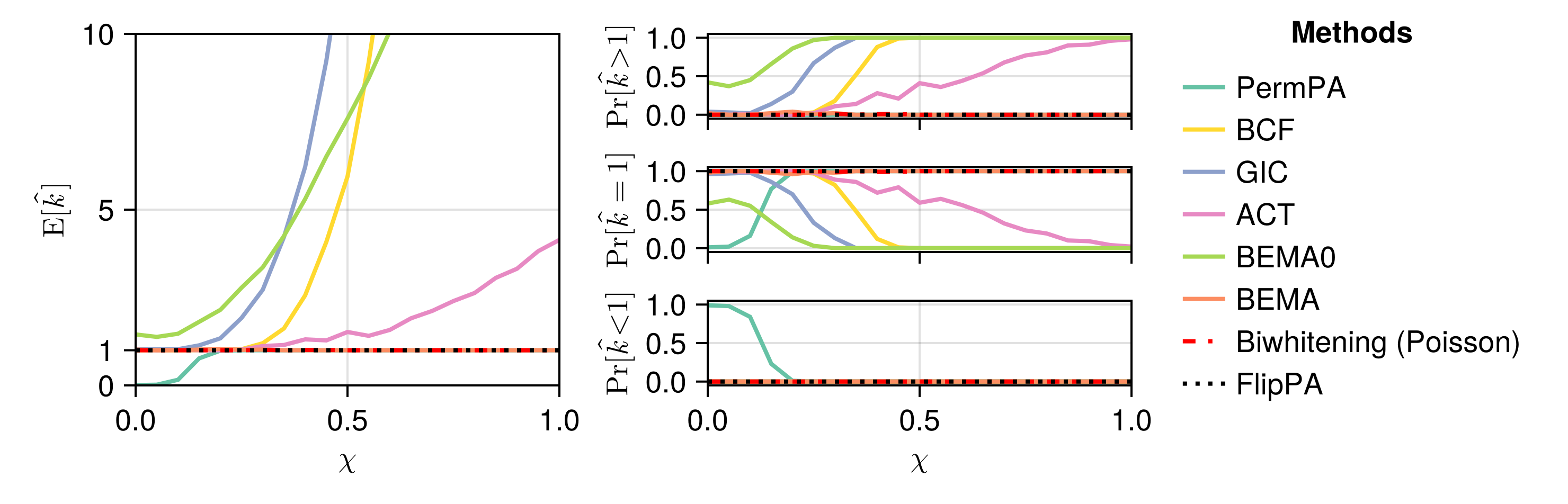}
    \caption{Poisson data \cref{eq:sim:rankone:poisson}
      with rank-one parameter matrix}
    \label{fig:sim:rankone:poisson}
  \end{subfigure}
  \caption{Performance across $100$ runs of each method for
    Bernoulli and Poisson data
    with a rank-one parameter matrix $\bmS$,
    as defined in \cref{eq:sim:rankone:bernoulli,eq:sim:rankone:poisson,eq:sim:rankone:factors}.
    As $\chi$ increases,
    $\bmS$ generally becomes more heterogeneous,
    which in turn produces more heteroscedasticity in the noise.}
  \label{fig:sim:rankone:bernoulli:poisson}
\end{figure}

Finally,
\cref{fig:sim:rankone:bernoulli:poisson} considers
Bernoulli and Poisson data
rather than the
rank-one signal plus
centered Bernoulli noise \cref{eq:sim:centered:bernoulli}
or centered Poisson noise \cref{eq:sim:centered:poisson}
used in \cref{fig:sim:block:centered:bernoulli,fig:sim:block:centered:poisson};
this matches how Bernoulli and Poisson noise
often arises in practice.
In particular,
\cref{fig:sim:rankone:bernoulli}
generates the $n \times p$ data matrix $\bmX$ as
\begin{equation}
  \label{eq:sim:rankone:bernoulli}
  X_{ij}
  \overset{ind}{\sim}
  \operatorname{Bernoulli}(S_{ij})
  \quad \text{with} \quad
  \bmS = \frac{1}{2} \cdot \frac{\bmu \bmz^\top}{\max(\bmu \bmz^\top)}
  \in \bbR^{n \times p}
  ,
\end{equation}
where
\begin{equation}
  \label{eq:sim:rankone:factors}
  u_i \overset{iid}{\sim} \operatorname{LogUniform}(e^{-\chi},e^{\chi})
  \quad \text{and} \quad
  z_i \overset{iid}{\sim} \operatorname{LogUniform}(e^{-\chi},e^{\chi})
  ,
\end{equation}
and we sweep $\chi$ from zero to one.
The entries of the noise matrix $\bmN = \bmX - \bmS$
are centered Bernoulli random variables
with means of $\bbE N_{ij} = \bbE X_{ij} - S_{ij} = S_{ij} - S_{ij} = 0$
and variances of $\var(N_{ij}) = \var(X_{ij}) = S_{ij}(1 - S_{ij})$.
Note that the scaling in $\bmS$ keeps
the maximum Bernoulli probability in the data fixed at $1/2$,
while
$\chi$ controls the degree of heterogeneity
in $\bmS$,
which consequently controls the degree of heterogeneity
in the noise variances.
\Cref{fig:sim:rankone:poisson}
considers a Poisson variant,
where the $n \times p$ data matrix $\bmX$ is generated as
\begin{equation}
  \label{eq:sim:rankone:poisson}
  X_{ij}
  \overset{ind}{\sim}
  \operatorname{Poisson}(S_{ij})
  \quad \text{with} \quad
  \bmS = \lambda \cdot \frac{\bmu \bmz^\top}{\operatorname{mean}(\bmu \bmz^\top)}
  \in \bbR^{n \times p}
  ,
\end{equation}
where $\lambda = 0.1$,
$\bmu$ and $\bmz$ are as above,
and we again sweep $\chi$ from zero to one.
The entries of the corresponding noise matrix $\bmN = \bmX - \bmS$
are centered Poisson random variables here
with means of $\bbE N_{ij} = \bbE X_{ij} - S_{ij} = S_{ij} - S_{ij} = 0$
and variances of $\var(N_{ij}) = \var(X_{ij}) = S_{ij}$.
The scaling of $\bmS$ in \cref{eq:sim:rankone:poisson}
keeps the average Poisson parameter in the data fixed at $\lambda = 0.1$,
while $\chi$ again controls the degree of heterogeneity
in the signal and in the noise variances.

The methods behaved somewhat similarly for
the Bernoulli data \cref{eq:sim:rankone:bernoulli} in \cref{fig:sim:rankone:bernoulli}
as for
the Poisson data \cref{eq:sim:rankone:poisson} in \cref{fig:sim:rankone:poisson}.
In both cases,
BCF, GIC, ACT, and BEMA0 all performed fairly well when $\chi$ was small
(though BEMA0 slightly overestimated the rank in the Poisson case);
the noise is close to homoscedastic in this regime.
They all overestimated the rank as $\chi$ grew.
Interestingly, PermPA
correctly estimated the rank
after around $\chi = 0.1$ in \cref{fig:sim:rankone:bernoulli}
and after around $\chi = 0.2$ in \cref{fig:sim:rankone:poisson},
even though permutations do not preserve
the heteroscedastic noise.
This seems to have been because permutations
were also not very effective at destroying the signal here
(e.g., $\bmS = \lambda \bm1_n \bm1_p^\top$ when $\chi = 0$,
which is unchanged by permutations).
Indeed, this is why PermPA underestimated the rank for small $\chi$.
Only BEMA, Biwhitening (Poisson),%
\footnote{The data in these experiments were nonnegative,
so we used the variant of Biwhitening designed for Poisson data
\citep[Algorithm~1.1]{landa2021brt:arxiv:v2}.
The variant that estimates a quadratic variance function
\citep[Algorithm~5.1]{landa2021brt:arxiv:v2}
performed similarly,
as did
the variant used in \cref{sec:simulations}
(which can handle the negative-valued data considered in that \lcnamecref{sec:simulations}).}
and FlipPA were effective
throughout the sweep.

Overall,
the experiments in this \lcnamecref{sec:sim:asymmetric}
indicate that FlipPA is robust to asymmetry in the noise distribution.
Indeed, all the methods considered here
generally performed similarly under asymmetric noise
as they did under the symmetric noise
used in \cref{sec:simulations}.
It should be noted, however,
that some common asymmetric settings
such as the Bernoulli and Poisson data settings
considered in \cref{fig:sim:rankone:bernoulli,fig:sim:rankone:poisson}
can introduce other important features beyond asymmetry.
For example,
as discussed above,
permutations were not very effective at destroying
these nonnegative signals.
These settings can also have a very strong component
that corresponds roughly to the mean of the data.
As discussed in \cref{sec:sim:shadowing},
strong factors like this can cause weak factors to be missed
in parallel analysis methods.
This can sometimes be addressed
by running FlipPA after centering or standardizing the data
(e.g., as done in Section~3.2 of \cite{landa2021brt:arxiv:v2})
to reduce this component.


\section{Power vs. type~I error for FlipPA}
\label{sec:sim:roc}

As usual,
the choice of the quantile $q$ in FlipPA affects both
the resulting type~I error (characterized in \cref{prop:type-I})
and the power.
In particular, $q$ tunes the trade-off between these two performance metrics:
increasing $q$ makes FlipPA select more conservatively
which reduces the type~I error but at the cost of also reducing the power.
However, it turns out that this trade-off can essentially vanish
when $n$ and $p$ are sufficiently large.

\begin{figure} \centering
  \includegraphics[scale=0.15]{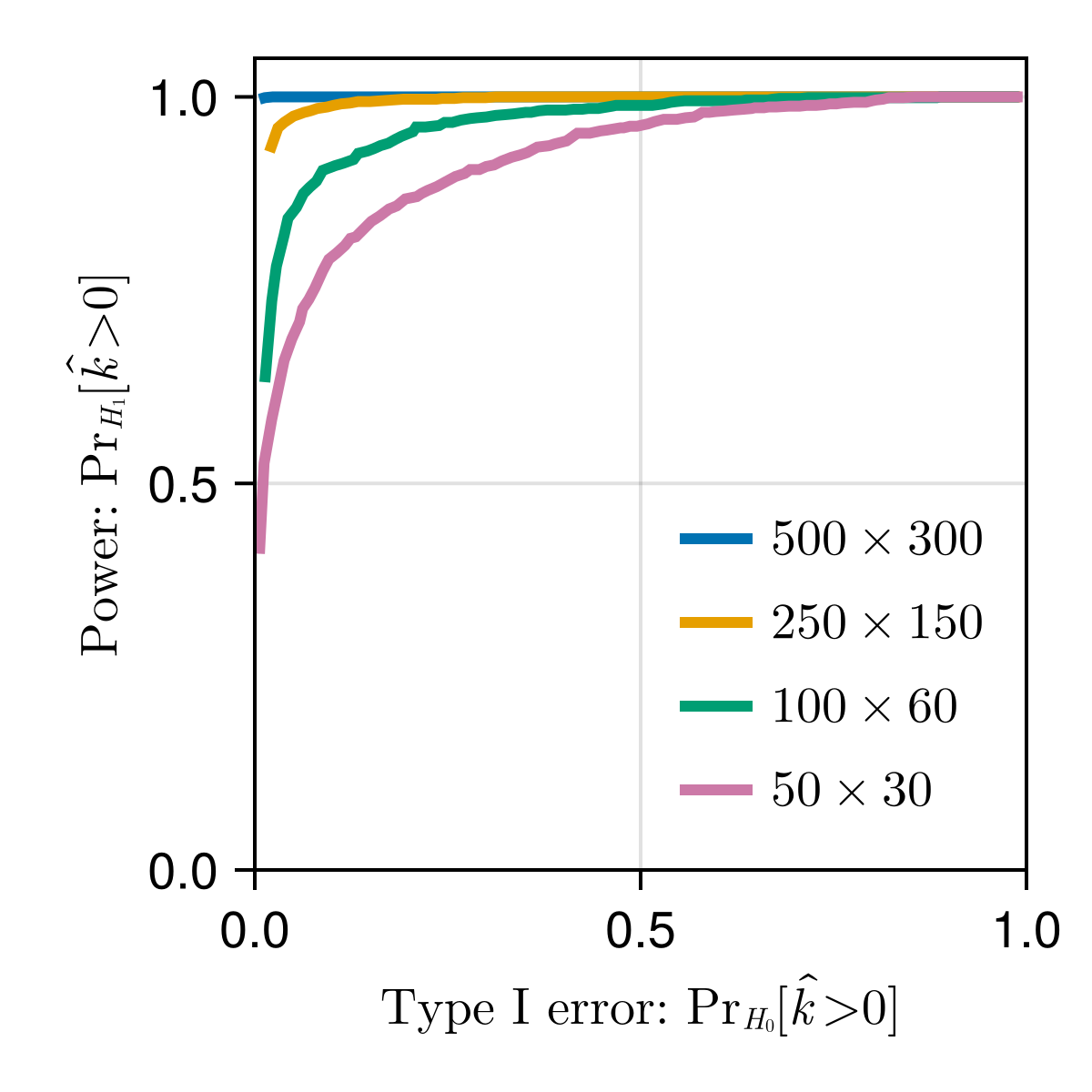}
  \caption{Trade-off curves showing the power ($\Pr_{H_1}\bigl[\htk > 0 \bigr]$) of FlipPA vs. its type~I error ($\Pr_{H_0}\bigl[\htk > 0 \bigr]$) for various matrix sizes $n \times p$.
  For sufficiently large $n$ and $p$,
  the trade-off essentially vanishes
  and choosing $q = 1$ simultaneously achieves power close to one
  with type~I error close to zero.}
  \label{fig:sim:roc}
\end{figure}

\Cref{fig:sim:roc} illustrates this phenomenon for various matrix sizes.
In particular,
it plots the power ($\Pr_{H_1}\bigl[\htk > 0 \bigr]$) of FlipPA
vs. its type~I error ($\Pr_{H_0}\bigl[\htk > 0 \bigr]$)
as $q$ is swept from zero to one.
Here, we consider the null hypothesis of no signal ($k=0$)
against the alternative of a rank-one signal ($k=1$)
of $\bmS = \theta \bmu \bmz^\top$
where
$\theta = 1.5$,
and
$\bmu \in \bbR^n$ and $\bmz \in \bbR^p$
are drawn uniformly from the respective unit spheres.
In both cases, the noise matrix is generated
with a block-structured noise variance profile
as follows:
\begin{equation}
  N_{ij} \overset{ind}{\sim} \clN(0, V_{ij}/n)
  \quad \text{where} \quad
  \bmV
  =
  \begin{bmatrix}
    0.25 \cdot \bm1_{n_1 \times p_1} &
    1          \cdot \bm1_{n_1 \times p_2} \\
    1.75 \cdot \bm1_{n_2 \times p_1} &
    1          \cdot \bm1_{n_2 \times p_2}
  \end{bmatrix}
  ,
\end{equation}
where $n_1 = n_2 = 0.5 n$, $p_1 = 0.8 p$, $p_2 = 0.2 p$,
and we consider $n \times p$ matrices of sizes
$50 \times 30$, $100 \times 60$, $250 \times 150$, and $500 \times 300$.
Essentially, this is the setting of \cref{fig:sim:block} with $\Delta = 0.75$
for varying matrix sizes where the block sizes are scaled accordingly.
The power and type~I error are computed empirically from $1000$ runs,
where we set $T = 100$ in FlipPA.

For the small $50 \times 30$ matrix,
we observe the expected trade-off.
Small type~I error (close to zero) corresponds to low power (under $1/2$);
this is the extreme where $q = 1$.
Decreasing $q$ towards zero increases the power towards one
but at the cost of a corresponding increase in type~I error.
However, this trade-off reduces as the matrix grows larger
and has essentially vanished once the matrix is $500 \times 300$.
Consequently,
$q = 1$ seems to be a reasonable default choice for the quantile in FlipPA
when the data is large;
for smaller data, one may consider reducing $q$ to increase power.
Indeed, as shown in \cref{sec:consistency},
FlipPA is asymptotically consistent under suitable conditions.


\section{Weak signals and bound-based thresholds}
\label{sec:weak:signals}

A general approach to estimating the rank is to select all the singular values
that rise above a threshold
set by using a high-probability upper-bound for the operator norm of the noise.
For example,
if the entries of the noise are bounded in absolute value by one,
the threshold can be set at $(2 + \eta)\sqrt{\max(n,p)}$ for some $\eta \geq 0$,
as is done by the universal singular value thresholding method
of \citet{chatterjee2015meb}.
Such approaches can be highly effective
in ``strong signal'' regimes,
where the singular values of the signal
diverge away from those of the noise.
Roughly speaking,
using this bound-based threshold
ensures that the noise singular values are not selected,
while the diverging nature of the signal singular values in this regime
ensures that they eventually rise above the threshold.

However, that is not the case for the ``weak signal'' regime
we have focused on,
where weak ``emergent'' factors
can produce signal singular values
that are of the same order as those of the noise.
In this regime, the threshold must be carefully calibrated
to match the operator norm of the noise,
and not simply upper bound it.
This requires properly accounting
for both
the aspect ratio of the data (i.e., the ratio of $n$ and $p$)
and the heterogeneity of the noise variances.
Indeed, the bound-based threshold above (which accounts for neither)
can be overly conservative in this regime.

\begin{figure} \centering
    \begin{subfigure}[b]{0.5\linewidth} \centering
        \includegraphics[scale=0.15]{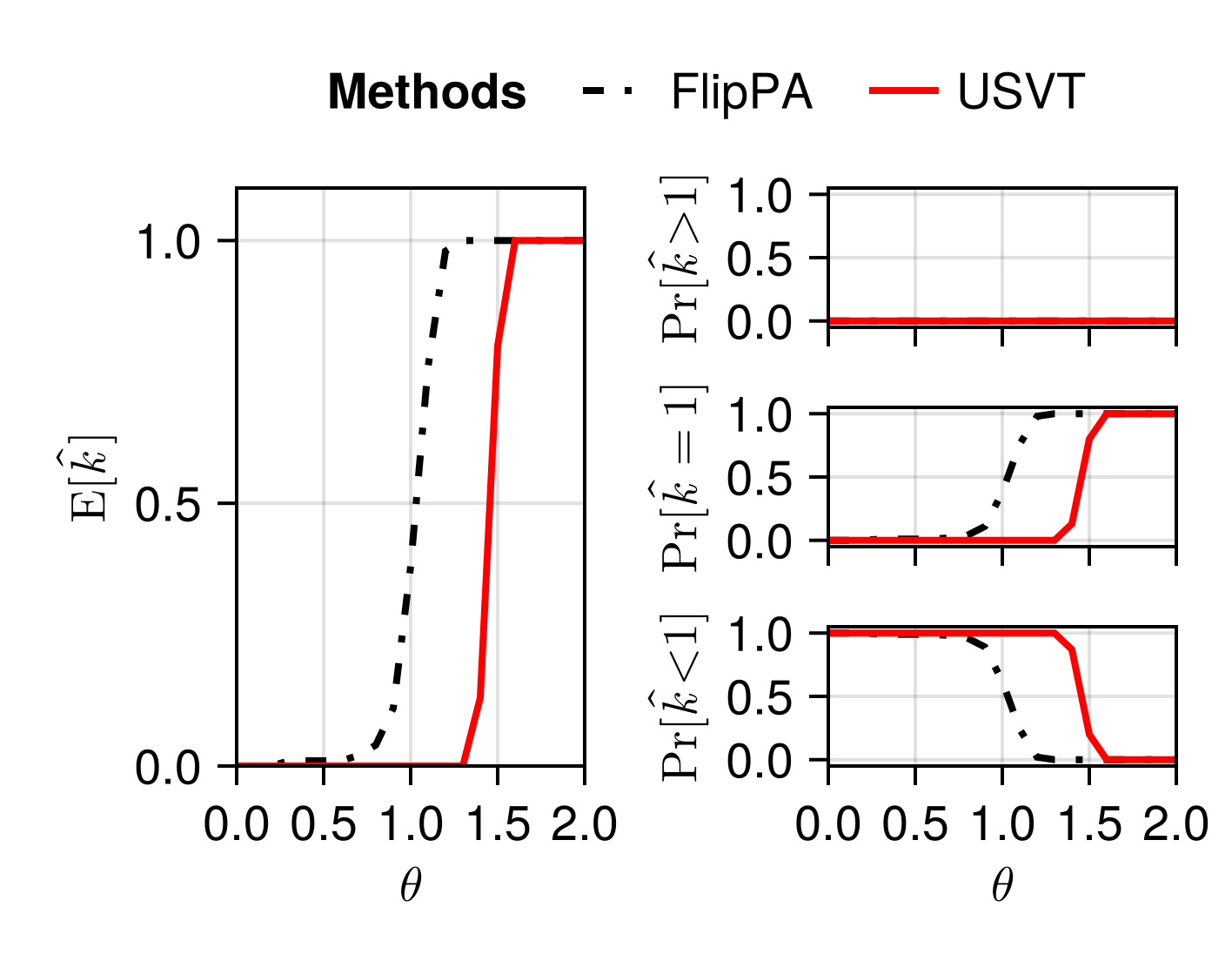}
        \caption{Comparison in the setting of \cref{fig:sim:hom}.}
        \label{fig:sim:weak:signals:hom}
    \end{subfigure}%
    \begin{subfigure}[b]{0.5\linewidth} \centering
        \includegraphics[scale=0.15]{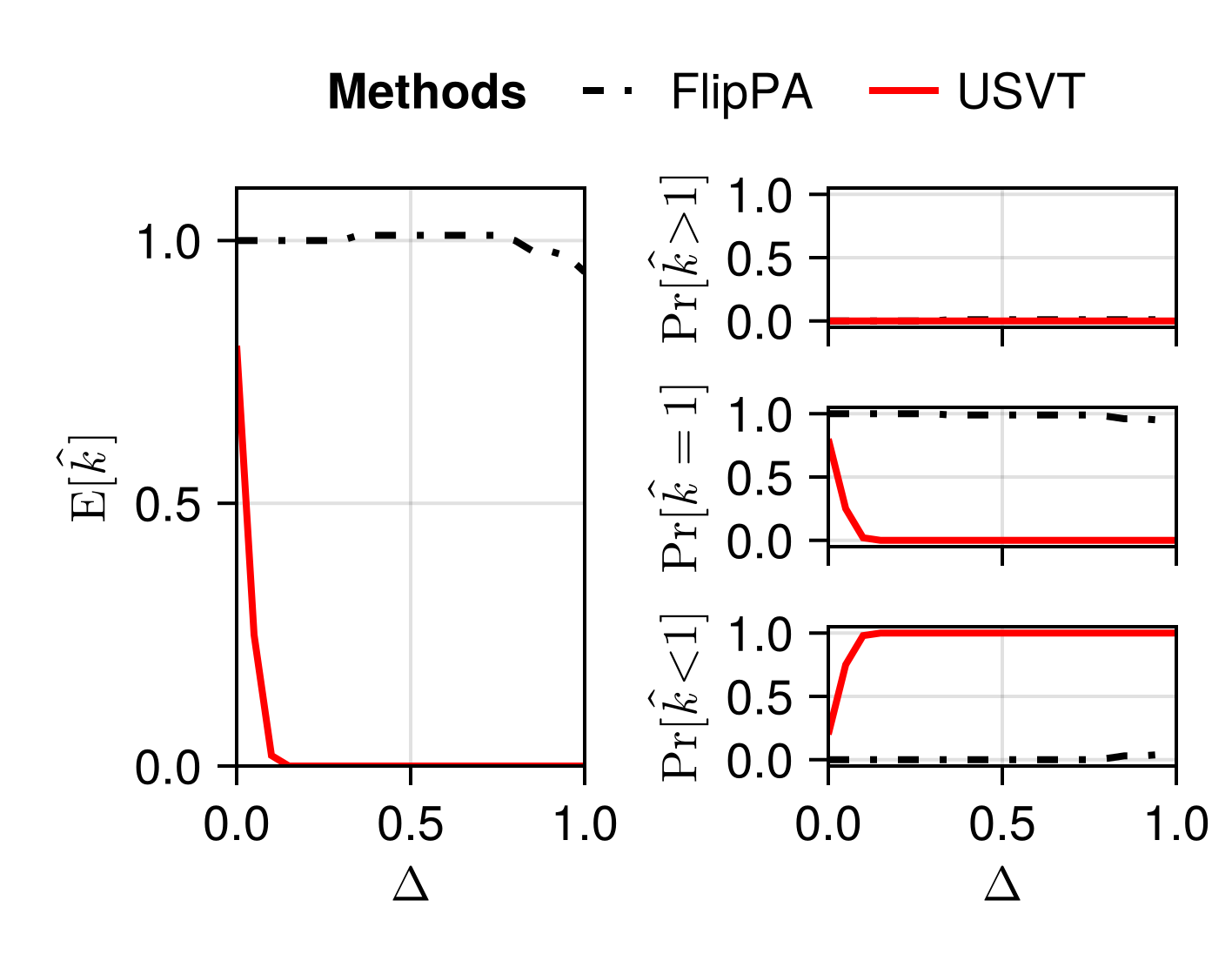}
        \caption{Comparison in the setting of \cref{fig:sim:block}.}
        \label{fig:sim:weak:signals:block}
    \end{subfigure}
    \caption{Comparison of FlipPA with
        universal singular value thresholding (USVT)
        from \citet{chatterjee2015meb},
        where the threshold is chosen using the true maximum variance.}
    \label{fig:sim:weak:signals}
\end{figure}

To illustrate this point,
here we
repeat the simulations of \cref{fig:sim:hom,fig:sim:block}
with the addition of universal singular value thresholding \citep{chatterjee2015meb}.
Since our noise is not bounded in absolute value by one,
we used the threshold of $(2 + \eta)\sqrt{v_{\max}\max(n,p)}$
suggested by \citet[Section~1.3]{chatterjee2015meb},
where $v_{\max}$ is the maximum variance of the noise entries.
\Cref{fig:sim:weak:signals} shows the performance
of universal singular value thresholding (USVT) with this threshold
compared with FlipPA
(we omit the other methods here to aid readability).

\Cref{fig:sim:weak:signals:hom} repeats the experiment of \cref{fig:sim:hom},
which sweeps the signal strength $\theta$ of a rank-one signal
in homogeneous noise;
see \cref{sec:sim:hom} for the detailed description.
Initially,
the signal is buried in the noise and neither method finds it.
Around $\theta = 1$, the signal rises above the noise
and is identified by FlipPA.
However, the bound-based threshold used by USVT
over-estimates the operator norm of the noise
so it does not find the signal until around $\theta = 1.5$.
\Cref{fig:sim:weak:signals:block}
repeats the experiment of \cref{fig:sim:block},
which considers heterogeneous noise with an increasing degree of heterogeneity $\Delta$,
where the rank-one signal has a constant strength;
see \cref{sec:sim:block} for the detailed description.
In this experiment,
the signal rises above the noise
and FlipPA finds it throughout the sweep.
Once again, the bound-based threshold used by USVT
over-estimates the operator norm of the noise,
resulting in under-selection.
Indeed,
this gap appears to grow as the degree of the heterogeneity $\Delta$ increases
(and the maximum variance used by USVT becomes less representative),
resulting in worse performance for USVT.
In contrast,
FlipPA uses a data-driven threshold that adapts to the noise heterogeneity
and correctly estimates the rank.


\section{Selection and preprocessing for astronomy data}
\label{sec:exp:quasar:spectra:preprocess}

Here we describe the selection and preprocessing steps
taken to produce the dataset used in \cref{sec:exp:quasar:spectra}
and illustrated in \cref{fig:exp:quasar:spectra:data};
it is essentially the same as the steps in \cite[Section~SM5]{hong2023owp}
except that here we do not remove spectra with highly heterogeneous variances.
Specifically,
the dataset was formed as follows:
\begin{enumerate}
  \item Select spectra from DR16Q that satisfy the following conditions:
  \begin{itemize}
    \item $\verb|SURVEY| = \verb|"eboss"|$,
    \item $\verb|PLATEQUALITY| = \verb|"good"|$,
    \item redshift: $2.0 < \verb|Z| < 2.1$,
    \item $\verb|BAL_PROB| < 0.2$,
    \item the rest frame wavelengths measured cover the range 1280--1820
      without missing entries
      (i.e., without entries for which $\verb|IVAR| = 0$).
  \end{itemize}
  \item Form a data vector $\bmy_i \in \bbR^p$
    and an associated variance vector $\bmv_i \in \bbR^p$
    for each of the $\tln = 10052$ selected spectra
    by linearly interpolating \verb|FLUX| and $1 \oslash \verb|IVAR|$
    to obtain entries for the $p=1081$ rest frame wavelengths
    $\verb|LAMREST| = (1280, 1280.5, \dots, 1820) \in \bbR^p$.
  \item Center each spectrum, i.e., $\bmy_i \gets \bmy_i - (1/p) \bm1_{p \times p} \bmy_i$ for $i = 1,\dots,\tln$.
  \item Normalize each spectrum so that its average value
    for rest frame wavelengths inside the range 1525--1575 is $\pm 1$,
    i.e., for $i = 1,\dots,\tln$,
  \begin{enumerate}
    \item $\sigma_i \gets |\verb|mean|(\bmy_i(1525 < \verb|LAMREST| < 1575))|$,
    \item $\bmy_i \gets \bmy_i / \sigma_i$,
    \item $\bmv_i \gets \bmv_i / \sigma_i^2$.
  \end{enumerate}
  \item Sort the spectra in increasing order
    of their average noise variance profiles $(1/p) \bm1_{p}^\top \bmv_i$.
  \item Select the $n=5000$ final (i.e., noisiest) spectra that remain
    after first dropping the 80 final (i.e., noisiest) spectra.
  \item Stack the selected vectors:
  $\bmY = [\bmy_1^\top; \dots; \bmy_n^\top] \in \bbR^{n \times p}$,
  $\bmV = [\bmv_1^\top; \dots; \bmv_n^\top] \in \bbR^{n \times p}$.
\end{enumerate}


\section{SDSS Acknowledgements}
\label{ack:sdss}

The example quasar spectra were provided by the Sloan Digital Sky Survey
\cite{ahumada2020t1d,lyke2020tsd}.
Funding for the Sloan Digital Sky 
Survey IV has been provided by the 
Alfred P. Sloan Foundation, the U.S. 
Department of Energy Office of 
Science, and the Participating 
Institutions. 

SDSS-IV acknowledges support and 
resources from the Center for High 
Performance Computing  at the 
University of Utah. The SDSS 
website is www.sdss4.org.

SDSS-IV is managed by the 
Astrophysical Research Consortium 
for the Participating Institutions 
of the SDSS Collaboration including 
the Brazilian Participation Group, 
the Carnegie Institution for Science, 
Carnegie Mellon University, Center for 
Astrophysics | Harvard \& 
Smithsonian, the Chilean Participation 
Group, the French Participation Group, 
Instituto de Astrof\'isica de 
Canarias, The Johns Hopkins 
University, Kavli Institute for the 
Physics and Mathematics of the 
Universe (IPMU) / University of 
Tokyo, the Korean Participation Group, 
Lawrence Berkeley National Laboratory, 
Leibniz Institut f\"ur Astrophysik 
Potsdam (AIP),  Max-Planck-Institut 
f\"ur Astronomie (MPIA Heidelberg), 
Max-Planck-Institut f\"ur 
Astrophysik (MPA Garching), 
Max-Planck-Institut f\"ur 
Extraterrestrische Physik (MPE), 
National Astronomical Observatories of 
China, New Mexico State University, 
New York University, University of 
Notre Dame, Observat\'ario 
Nacional / MCTI, The Ohio State 
University, Pennsylvania State 
University, Shanghai 
Astronomical Observatory, United 
Kingdom Participation Group, 
Universidad Nacional Aut\'onoma 
de M\'exico, University of Arizona, 
University of Colorado Boulder, 
University of Oxford, University of 
Portsmouth, University of Utah, 
University of Virginia, University 
of Washington, University of 
Wisconsin, Vanderbilt University, 
and Yale University.


\edit{
\section{AI Acknowledgements}
\label{ack:ai}

GPT-5.5 was used to aid in the development of conjectures and initial proofs for the theoretical results in \cref{sec:relative:strength}, as well as for the development of prototype scripts for some of the experiments in \cref{sec:sim:shadowing}.
Moreover, it was used to help develop initial versions of some of the corresponding discussions, and to catch a few typos and inconsistencies throughout the paper.
Claude Sonnet 4.6 was also used to catch a few typos and inconsistencies throughout the paper.
All the content has been rigorously revised by the authors.
}

\end{document}